# Local recognition of reflection graphs on Coxeter groups

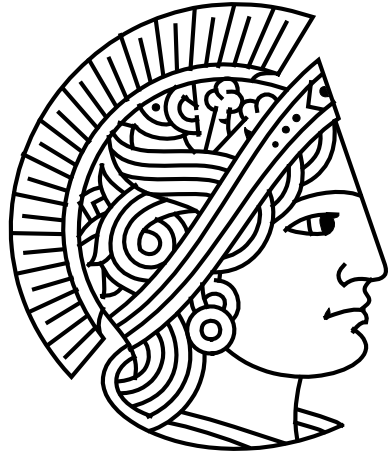

**Diploma Thesis**

*31 Mar 2008*
(minor corrections, 13 May 2008)

Armin Straub

Technische Universität Darmstadt
Department of Mathematics
*Algebra, Geometry, Functional Analysis*

supervised by
PD dr. Ralf Gramlich

# Contents

# Notations

# 1 Introduction

## 1.1 Overview

Given a graph $\Gamma$ one may ask to which extent it is determined by its local graphs that is by the induced subgraphs on the vertices adjacent to a particular vertex. If all these local graphs are isomorphic to a single graph $\Lambda$ then $\Gamma$ is locally homogeneous and we say that $\Gamma$ is locally $\Lambda$. A classification of the graphs which are locally a given graph $\Lambda$ is called a local recognition result. The local recognition of graphs has been studied extensively in the literature, for instance in [BH77], [HS85], [Gra02], [Gra04], [CCG05], [Alt07], and a particularly guiding example is the local recognition of the Kneser graphs studied in [Hal80] and [Hal87].

In this thesis we are interested in the local recognition of Weyl graphs $\mathbb{W}(M)$ which are commuting graphs on the reflections of irreducible Coxeter groups. If the corresponding Dynkin diagram $M$ is simply laced then we identify $\mathbb{W}(M)$ with locally cotriangular graphs as studied and locally recognized in [HS85]. On the other hand, if $M$ is not simply laced then the Weyl graph $\mathbb{W}(M)$ is a bichromatic graph. We prove that the Weyl graphs $\mathbb{W}(B_n)$ and $\mathbb{W}(C_n)$ are locally recognizable, that is that they are characterized among connected graphs by their local graphs. On the contrary, the Weyl graph $\mathbb{W}(F_4)$ is not locally recognizable. We therefore study bichromatic graphs which are locally like $\mathbb{W}(F_4)$, and in the sequel obtain several characterizations of $\mathbb{W}(F_4)$ as one of two tightest bichromatic graphs which are locally like $\mathbb{W}(F_4)$.

In the last section, we turn to group theoretical applications of local recognition results for graphs. The paradigmatic and guiding result exemplified in [GLS94, Theorem 27.1] is the characterization of the symmetric groups by means of the structure of its transposition centralizers. We give a similar characterization for the Coxeter group of type $F_4$ in terms of its reflection centralizers. The interest in such local recognition results for groups stems from the classification of finite simple groups and the fact that the majority of finite simple groups arises from (possibly twisted) Chevalley groups. These can be defined similar to Coxeter groups as groups generated by $SL(2, q)$ subgroups subject to certain relations. Local recognition results for instance for Chevalley groups of type $A_n$, $n \geqslant 8$, based on graph theoretical results have been studied in [Gra02], [Gra04]. Recently, Ralf Gramlich and Kristina Altmann proved a local recognition result for Chevalley groups of type $A_7$ and $E_6$ based on results of [Alt07] and making use of the local recognition of graphs that are locally $\mathbb{W}(A_5)$, see [Gra08]. We hope that this thesis can help to approach a similar recognition result for Chevalley groups of type $F_4$.

### 1.1.1 Organization

This thesis is organized in such a way that later sections depend on earlier ones. The introductory section mainly contains the necessary background on groups in general, graphs as well as Coxeter groups. These parts may as well be skipped if the reader feels sufficiently familiar with the presented content.

Examples often contain definitions and basic observations that are used later on, and therefore form a logically relevant part of this thesis. Remarks, on the other hand, usually contain comments that are meant to supplement or to motivate definitions or results. To meet this goal some of the remarks certainly are handwavy in nature. They may, however, be omitted without impact on the logical coherence of the thesis.

At some points we use the computer algebra systems GAP, see [GAP07], and SAGE, see [SAG07], for computations. These computations are deferred to the appendix.

### 1.1.2 General notations

Throughout this thesis, the letters $X$, $Y$ are used to denote sets. Likewise, $G$, $H$ denote groups, and $\Gamma$, $\Lambda$ denote graphs. $n$ represents a natural number, $p$ a prime, and $q$ a prime power. $(W, S)$ is used to denote a Coxeter system, $\Phi$ a root system, and $M$ a diagram.

Indices often are not quantified and are meant to be chosen in such a way that the resulting expression is defined. For instance, suppose that the elements $x_1, x_2, ..., x_n$ have been defined. Then $x_i$ is used to refer to one of the vertices $x_i$ with the implicit assumption that $i \in \{1, 2, ..., n\}$.

We write $x \triangleq y$ to mean that $x$ is defined to be $y$. The maximum of two numbers $x$, $y$ is denoted by $x \vee y$. Moreover, we denote the $n$-element set $\{1, 2, ..., n\}$ by $[n]$, and the $k$-subsets of a set $X$ by

$$\binom{X}{k} \triangleq \{A \subset X \colon |A| = k\}.$$

The disjoint union of two sets $A, B$ is written as $A \sqcup B$.

Finally, we write $A = B$ instead of $A \cong B$ if $A$ and $B$ can be identified in a canonical way. For instance, we write $\mathrm{Aut}(K_n) = \mathrm{Sym}_n$.

### 1.1.3 Acknowledgments

I would like to thank my thesis advisor Ralf Gramlich for his constant and helpful support throughout the writing of this thesis, for always nurturing my enthusiasm, for sharing his enviable insights, for borrowing several articles and books, for numerous discussions, and for reading and carefully commenting on early drafts of this thesis. I would also like to thank Hendrik Van Maldeghem for a fruitful discussion and for suggesting the approach outlined in Remark 3.14. Finally, I'm grateful to Ralf Gramlich and Jonathan I. Hall for making available to me a note of theirs that served as the starting point for this thesis.

## 1.2 Group theory

We assume that the reader is familiar with the basic notions of group theory. For a gentle introduction we refer to [DF04]. A thorough treatment can also be found in the classical texts [Bou89] and [Bou90]. This section is mainly used to fix notation and to introduce some of the groups that will be of importance later on.

Let $G$ be a group. We write $H \leqslant G$ to mean that $H$ is a subgroup of $G$. If $H \leqslant G$ we denote the set of (left) cosets of $H$ in $G$ by $G/H$. These cosets form a group, called the *quotient group*, with multiplication given by $(gH)(hH) = (gh)H$ if and only if $H$ is *normal* in $G$. Recall that $H$ is defined to be normal in $G$ if $gHg^{-1} = H$ for all $g \in G$. A group is said to be *simple* if it contains no nontrivial proper normal subgroups. We write $H \trianglelefteq G$ to denote that $H$ is normal in $G$. If $H \trianglelefteq G$ then $G/H$ also denotes the quotient group.

Let $G$ be a group, and $X \subseteq G$ a subset of elements. We write $\langle X \rangle$ for the subgroup of $G$ generated by $X$, that is the smallest subgroup of $G$ containing $X$. If $Y \subseteq G$ is another subset we will sometimes write $\langle X, Y \rangle$ to mean $\langle X \cup Y \rangle$. Now, let $X$ be any set and denote with $F(X)$ the free group on $X$. Given $R \subseteq F(X)$, let $N$ be the smallest normal subgroup of $F(X)$ containing $R$. We define

$$\langle X \colon R \rangle \triangleq \frac{F(X)}{N},$$

and say that $\langle X \colon R \rangle$ is the group generated by $X$ with respect to the *relations* $R$. In practice, we often use the notation

$$\langle X \colon r_1 = r_2 = ... = r_n = 1 \rangle \triangleq \langle X \colon \{r_1, r_2, ..., r_n\} \rangle.$$

If $G$ is a group and $G \cong \langle X \colon R \rangle$ for some $X, R$ then we say that $\langle X \colon R \rangle$ is a *presentation* for $G$. Note that the groups generated by $X$ in which the relations $R$ hold are precisely the quotient groups of $\langle X \colon R \rangle$. More on presentations of groups can be found for instance in [Bou89, §7].

### 1.2.1 Examples of groups

**Example 1.1.** The *cyclic group* of order $n$ is denoted by $\mathbb{Z}/n$.

**Example 1.2.** Let $X$ be a finite set. The *symmetric group* $\mathrm{Sym}_X$ is the group of all permutations of $X$. We abridge $\mathrm{Sym}_n$ for $\mathrm{Sym}_{[n]}$. Note that $\mathrm{Sym}_X \cong \mathrm{Sym}_n$ where $n$ is the cardinality of $X$. This explains why we will mostly work with $\mathrm{Sym}_n$. Permutations of $\mathrm{Sym}_n$ will be denoted in the usual *cycle* notation: the transposition interchanging $i \neq j \in [n]$ is written as $(i\ j)$ and likewise the cycle $(i_1\ i_2\ ...\ i_m)$ refers to the permutation $\sigma$ defined by mapping $\sigma(i_1) = i_2$, $\sigma(i_2) = i_3$, ..., $\sigma(i_m) = i_1$ and fixing every other $j \in [n]$. Since they refer to the same permutation we regard cyclic shifts of a cycle as the same. For instance, $(i_2\ i_3\ ...\ i_m\ i_1) = (i_1\ i_2\ ...\ i_m)$. The cycle $(i_1\ i_2\ ...\ i_m)$ is said to be of length $m$ and is therefore also called an $m$-cycle. Note that $m$-cycles have order $m$ in $\mathrm{Sym}_n$. Juxtaposition of cycles like $(i_1\ i_2\ ...\ i_m)(j_1\ j_2\ ...\ j_k)$ refer to the composition $(i_1\ i_2\ ...\ i_m) \circ (j_1\ j_2\ ...\ j_k)$ of the associated permutations. For instance, $(1\ 2)(2\ 3) = (1\ 2\ 3)$. The support of a cycle $(i_1\ i_2\ ...\ i_m)$ is the set $\{i_1, i_2, ..., i_m\}$, and two cycles are said to be *disjoint* if their supports are disjoint. Notice that disjoint cycles commute. It is easy to see that every permutation $\sigma \in \mathrm{Sym}_n$ can be written as a product of disjoint cycles and that this decomposition is unique up to reordering the cycles. Two permutations are conjugate in $\mathrm{Sym}_n$ if and only if the number of $k$-cycles in their disjoint cycle decomposition is equal for $k \in [n]$.

$\mathrm{Sym}_n$ is generated by its transpositions. In fact, it is already generated by the $n-1$ transpositions $(1\ 2), (2\ 3), ..., (n-1\ n)$ as will become important later on.

**Example 1.3.** Since $\mathrm{Sym}_n$ is generated by transpositions we can write every permutation as a product of transpositions. We say that a permutation is *even* if it can be written as a product of an even number of transpositions. In fact, the even permutations form a group. This is the *alternating group* $\mathrm{Alt}_n$.

**Example 1.4.** The *dihedral group* $\mathrm{Dih}_n$ is the group of symmetries of a regular $n$-gon. $\mathrm{Dih}_n$ is a group of order $2n$ consisting of the $n$ rotations and $n$ reflections.

**Example 1.5.** Let $\mathbb{F}$ be a field. The *general linear group* $\mathrm{GL}(n, \mathbb{F})$ is the group of $\mathbb{F}$-linear automorphisms of the vector space $\mathbb{F}^n$ over $\mathbb{F}$. We identify $\mathrm{GL}(n, \mathbb{F})$ with the $n \times n$ invertible matrices with entries in $\mathbb{F}$. As usual, $\det\colon \mathrm{GL}(n, \mathbb{F}) \to \mathbb{F}^\times$ denotes the determinant homomorphism. Here, $\mathbb{F}^\times = \mathbb{F}\setminus\{0\}$ are the invertible elements. The *special linear group* $\mathrm{SL}(n, \mathbb{F})$ is the kernel of $\det$, that is $\mathrm{SL}(n, \mathbb{F}) \leqslant \mathrm{GL}(n, \mathbb{F})$ consists of those $n \times n$ matrices over $\mathbb{F}$ which have determinant 1. The *projective special linear group* is the group

$$\mathrm{PSL}(n, \mathbb{F}) \triangleq \frac{\mathrm{SL}(n, \mathbb{F})}{Z(\mathrm{SL}(n, \mathbb{F}))}$$

where $Z(\mathrm{SL}(n, \mathbb{F}))$ is the center of $\mathrm{SL}(n, \mathbb{F})$ made up by the scalar multiples of the identity matrix.

Recall that for any prime power $q$ there exists (up to isomorphism) exactly one field with $q$ elements, see for instance [Bou90, §12]. We denote this field by $\mathbb{F}_q$, and usually abbreviate $\mathrm{GL}(n, \mathbb{F}_q)$ as $\mathrm{GL}(n, q)$. We do the same for $\mathrm{SL}(n, \mathbb{F}_q)$ and $\mathrm{PSL}(n, \mathbb{F}_q)$.

More on general linear groups can be found for instance in [Rot95, Chapter 8].

**Remark 1.6.** Consider the $n$-dimensional vector space $\mathbb{F}_q^n$. There are $q^n - 1$ nonzero vectors in $\mathbb{F}_q^n$ each of which determines a 1-dimensional subspace. Since every 1-dimensional subspace contains $q - 1$ nonzero vectors we count

$$[n]_q \triangleq \frac{q^n - 1}{q - 1} = 1 + q + ... + q^{n-1}$$

1-dimensional subspaces of $\mathbb{F}_q^n$. $[n]_q$ is said to be a *q-analog* of the number $n$ because if we let $q \to 1$ then $[n]_q \to n$. With respect to the number of elements, an $n$-element set therefore behaves like an $n$-dimensional vector space over the non-existing field $\mathbb{F}_1$. In fact, this vague analogy goes much further. For instance, the number of $k$-dimensional subspaces of $\mathbb{F}_q^n$ is given by the $q$-binomial coefficient

$$\begin{bmatrix} n \\ k \end{bmatrix}_q \triangleq \frac{[n]_q!}{[n-k]_q!\,[k]_q!}$$

where $[n]_q!$ is the $q$-factorial defined as $[n]_q [n-1]_q \cdots [1]_q$. Two proofs and a lot more beautiful counting can be found in [And76, 13.1]. Again, if we let $q$ approach 1 then

$$\lim_{q \to 1} \begin{bmatrix} n \\ k \end{bmatrix}_q = \binom{n}{k},$$

that is we get the number of $k$-element subsets of an $n$-element set.

We therefore informally say that $n$-dimensional vector spaces over the finite field $\mathbb{F}_q$ are a $q$-analog of $n$-element sets (actually, the analogy usually is most visible between projective spaces $\mathbb{P}_n(\mathbb{F}_q)$ and $n+1$-elements sets). Likewise, there are striking similarities between the general linear group $\mathrm{GL}(n, q)$ and the symmetric group $\mathrm{Sym}_n$. Just to give an example, let us determine the order of $\mathrm{GL}(n, q)$ which is equivalent to counting invertible $n \times n$ matrices over $\mathbb{F}_q$. For the first column of such a matrix we may choose any nonzero vector, for the second column we may then choose any vector which is not a multiple of the first column, for the third column we may choose any vector not in the span of the first two columns, and so on. We therefore count

$$\begin{aligned} |\mathrm{GL}(n,q)| &= (q^n - 1)(q^n - q) \cdots (q^n - q^{n-1}) \\ &= \frac{q^n - 1}{q - 1} \frac{q^{n-1} - 1}{q - 1} \cdots \frac{q - 1}{q - 1} q^{\binom{n}{2}} (q-1)^n \\ &= [n]_q!\, q^{\binom{n}{2}} (q-1)^n. \end{aligned}$$

Observe the occurrence of $[n]_q!$ which is a $q$-analog of $n!$, the order of $\mathrm{Sym}_n$.

In this text, we shall encounter at least two further examples of $q$-analogs. Namely, we shall mention a $q$-analog of Kneser graphs, see Remark 2.15, and we will briefly introduce Chevalley groups as $q$-analogs of Coxeter groups, see Remark 3.25.

**Example 1.7.** Given a bilinear form $B \colon \mathbb{F}^n \times \mathbb{F}^n \to \mathbb{F}$ we say that a matrix $M \in \mathrm{GL}(n, \mathbb{F})$ preserves $B$ if $B(x, y) = B(Mx, My)$ for all $x, y \in \mathbb{F}^n$. The matrices preserving $B$ form a subgroup of $\mathrm{GL}(n, \mathbb{F})$. The *orthogonal group* $O(n, \mathbb{F})$ is the group of matrices preserving the bilinear form $\mathbb{F}^n \times \mathbb{F}^n \to \mathbb{F}$ defined by

$$(x, y) \mapsto x_1 y_1 + x_2 y_2 + \ldots + x_n y_n.$$

If $\mathbb{F} = \mathbb{R}$ then the orthogonal group $O(n, \mathbb{R})$ consists precisely of the linear isometries of $\mathbb{R}^n$.

**Example 1.8.** A bilinear form $B$ is said to be *alternating* if $B(x, x) = 0$ for all $x \in \mathbb{F}^n$. A *symplectic form* is an alternating bilinear form. $B: \mathbb{F}^{2n} \times \mathbb{F}^{2n} \to \mathbb{F}$ defined by

$$B(x, y) = x_1 y_2 - x_2 y_1 + x_3 y_4 - \ldots + x_{2n-1} y_{2n} - x_{2n} y_{2n-1}$$

is a nondegenerate symplectic form referred to as the *standard symplectic form* on $\mathbb{F}^{2n}$. The group preserving this form is called the *symplectic group* and is denoted by $\mathrm{Sp}(2n, \mathbb{F})$. Again, we usually write $\mathrm{Sp}(2n, q)$ for $\mathrm{Sp}(2n, \mathbb{F}_q)$. As it turns out, the symplectic group $\mathrm{Sp}(2n, \mathbb{F})$ is a subgroup of the special linear group $\mathrm{SL}(2n, \mathbb{F})$. Some more details on symplectic groups can be found for instance in [Gar97, 8.1].

### 1.2.2 Group properties

Let $G$ be a group. The *center* of $G$ defined as

$$Z(G) \triangleq \{g \in G : (\forall h \in G) gh = hg\}$$

is a normal subgroup. We will sometimes abbreviate the quotient $G/Z(G)$ as $G/Z$.

**Proposition 1.9.** *Let $G$ and $H$ be groups.*
- $Z(G \times H) = Z(G) \times Z(H)$.
- *If $Z(G) \leqslant H$ then $Z(G) \leqslant Z(H)$ provided that $H \leqslant G$.*

**Proof.** The first claim is obvious. If $H \leqslant G$ then $Z(G) \cap H \leqslant Z(H)$ and the second claim follows. □

The *commutator* of $g, h \in G$ is the element

$$[g, h] \triangleq g h g^{-1} h^{-1}.$$

Notice that $[g, h] = 1$ if and only $g$ and $h$ commute. Analogously, for subsets $X, Y \subseteq G$ we define

$$[X, Y] \triangleq \langle [x, y] : x \in X, y \in Y \rangle.$$

The *commutator subgroup* of $G$ is the normal subgroup $G' \triangleq [G, G]$. It is the smallest normal subgroup of $G$ such that the quotient group is abelian. A group $G$ is said to be *perfect* if $G' = G$.

**Example 1.10.** Every nonabelian simple group $G$ is perfect. This is just a consequence of the fact that $G$ contains no nontrivial normal subgroup other than itself.

**Example 1.11.** For $n \geqslant 5$ the alternating groups $\mathrm{Alt}_n$ are simple, see for instance [DF04, 4.6.24], and therefore perfect. Note that every commutator in the symmetric group $\mathrm{Sym}_n$ is even. The commutator subgroup of $\mathrm{Sym}_n$ is therefore $\mathrm{Alt}_n$ when $n \geqslant 5$. In fact, one checks that for all $n$ the commutator subgroup of $\mathrm{Sym}_n$ is $\mathrm{Alt}_n$.

**Example 1.12.** The commutator subgroup of $\mathrm{GL}(n,q)$ is $\mathrm{SL}(n,q)$ except when $(n,q)=(2,2)$, see [GLS94, A.1].

**Remark 1.13.** The groups $\mathrm{GL}(n,\mathbb{F})$ and certain subgroups including $\mathrm{SL}(n,\mathbb{F})$ and $\mathrm{Sp}(2n,\mathbb{F})$ are often referred to as the *classical groups*, a term coined by Hermann Weyl. These groups are of particular interest because of the following fact proved by Leonard E. Dickson and Jean Dieudonné. If $G$ is a classical group then up to some exceptions the group $G'/Z(G')$ is simple. For further reading on classical groups the enjoyable overview [Cur67] which predates the final classification of the finite simple groups is recommended to the reader.

### 1.2.3 Group actions

Let $X$ be a set, possibly with additional structure. By this we essentially mean that the group of automorphisms of $X$, that is the group of all bijective mappings $X \to X$ preserving the structure on $X$, may be a subgroup of the group $\mathrm{Sym}_X$ of all bijective mappings on $X$. We denote the automorphism group of $X$ by $\mathrm{Aut}(X)$. A (left) *group action* of a group $G$ on $X$ is a homomorphism $a\colon G \to \mathrm{Aut}(X)$. We usually write

$$g \cdot x \triangleq a(g)(x)$$

for $g \in G$ and $x \in X$ if there is no possible confusion about the action $a$.

**Example 1.14.** Let $X$ be the vector space $\mathbb{F}^n$ over some field $\mathbb{F}$. In this case $\mathrm{Aut}(X)$ is the general linear group $\mathrm{GL}(n,\mathbb{F})$. Fix a basis $B = \{\varepsilon_1, \varepsilon_2, \ldots, \varepsilon_n\}$ of $X$. Let $g \in \mathrm{Sym}_B$ and set

$$g \cdot (c_1 \varepsilon_1 + c_2 \varepsilon_2 + \ldots + c_n \varepsilon_n) = c_1\, g(\varepsilon_1) + c_2\, g(\varepsilon_2) + \ldots + c_n\, g(\varepsilon_n)$$

for $c_i \in \mathbb{F}$. This defines an action of $\mathrm{Sym}_B$ on $X$.

The *kernel* of the group action is the kernel of the homomorphism $G \to \mathrm{Aut}\, X$. An action is said to be *faithful* if this kernel is trivial that is if $G \to \mathrm{Aut}\, X$ is injective. Observe that if $N$ is the kernel of the action of $G$ on $X$ then $G/N$ acts faithfully on $X$. The *stabilizer* of $x \in X$ in $G$ is the subgroup

$$C_G(x) \triangleq \{g \in G\colon g \cdot x = x\}.$$

On the other hand, the *orbit* of $x$ is the set

$$G \cdot x \triangleq \{g \cdot x \colon g \in G\}.$$

Note that $X$ is partitioned by the orbits $G \cdot x$, $x \in X$. We say that the action of $G$ on $X$ is *transitive* if there is only a single orbit or, equivalently, if for every $x, y \in X$ there exists $g \in G$ such that $g \cdot x = y$.

The next result is elementary and easily checked.

**Lemma 1.15.** *(**Orbit-stabilizer formula**) Let $G$ act on $X$, and $x \in X$. Then the mapping $G/C_G(x) \to G \cdot x$ given by $gC_G(x) \mapsto g \cdot x$ is a bijection. In particular, if $G$ is finite then*

$$|G| = |C_G(x)||G \cdot x|. \qquad \square$$

Given a group action of $G$ on $X$ there are several induced actions which we regard as natural. For instance, any subgroup of $G$ acts on $X$ as well, $G$ acts on $X \times X$ given by $g \cdot (x, y) = (g \cdot x, g \cdot y)$, or $G$ acts on the $k$-subsets of $X$ by $g \cdot \{x_1, ..., x_k\} = \{g \cdot x_1, ..., g \cdot x_k\}$.

**Example 1.16.** Every group $G$ acts on itself by *conjugation* that is by $g \cdot h = ghg^{-1}$. Since this notation is potentially confusing we adopt the notation

$$h^g \triangleq g^{-1}hg.$$

Another reason for doing so is that this notation is used in the computer algebra system GAP which we will use for some calculations. For a subset $H$ of $G$ we analogously write $H^g$ to mean $g^{-1}Hg$. Likewise, we denote the conjugacy class of $h$ in $G$ as

$$h^G \triangleq \{h^g : g \in G\}.$$

Note that the stabilizer $C_G(h)$ is the *centralizer* of $h$ in $G$, that is the subgroup of $G$ of all elements commuting with $h$.

The ubiquity of the conjugation action is demonstrated in the following easily verified fact about the relation of the stabilizer subgroups with respect to arbitrary group actions.

**Proposition 1.17.** *Let $G$ act on $X$, $g \in G$, and $x \in X$. Then*

$$C_G(g \cdot x) = g C_G(x) g^{-1}. \qquad \square$$

The following observations are elementary but will prove helpful later on.

**Proposition 1.18.** *Let $G$ be a group and $x, y \in G$ be involutions.*
- *$(xy)^2 = 1$ if and only if $[x, y] = 1$.*
- *$(xy)^3 = 1$ if and only if $x^y = y^x$.*
- *$(xy)^4 = 1$ if and only if $[x, x^y] = 1$.*
- *$(xy)^n = 1$ for odd $n$ implies that $x$ and $y$ are conjugate in $G$.* $\qquad \square$

### 1.2.4 Group products

**Definition 1.19.** *Let $G, H$ be groups and $\alpha$ be an action of $H$ on $G$. The* semidirect product *$G \rtimes_\alpha H$ is the group with elements $G \times H$ and multiplication given by*

$$(g_1, h_1)(g_2, h_2) = (g_1 \, \alpha(h_1)(g_2), h_1 h_2).$$

**Example 1.20.** *Let $s, t$ be two distinct involutions of a group $G$. Then $\langle s, t \rangle$ is a dihedral group. Furthermore, $\langle s, t \rangle = \langle st \rangle \rtimes \langle s \rangle$.*

**Example 1.21.** *Let $\mathbb{F}$ be a field and $\mathrm{Aut}(\mathbb{F})$ its automorphism group. Note that $\mathrm{Aut}(\mathbb{F})$ acts on the general linear group $\mathrm{GL}(n, \mathbb{F})$ componentwise. The corresponding semidirect product*

$$\Gamma L(n, \mathbb{F}) \triangleq \mathrm{GL}(n, \mathbb{F}) \rtimes \mathrm{Aut}(\mathbb{F})$$

*is called the* semilinear group. *Likewise, the* special semilinear group *is*

$$\Sigma L(n, \mathbb{F}) \triangleq \mathrm{SL}(n, \mathbb{F}) \rtimes \mathrm{Aut}(\mathbb{F}).$$

Let $X$ be a set, and $G$ a group. We write $G^{(X)}$ for the direct sum of copies of $G$ indexed by $X$. To be precise, $G^{(X)}$ are the functions $X \to G$ which take only finitely many nontrivial values. Of course, the group structure of $G^{(X)}$ is defined componentwise.

**Definition 1.22.** *Let $G, H$ be groups, and let $H$ act on the set $X$. The* wreath product *$G \wr (H, X)$ is the semidirect product $G^{(X)} \rtimes_\alpha H$ where $\alpha \colon H \to \mathrm{Aut}(G^{(X)})$ is defined by*

$$\alpha(h)(f) = f \circ h^{-1}.$$

*We abbreviate $G \wr H$ if there is no confusion about the set $X$.*

The following two lemmata describe a useful criterion for recognizing when a group is a direct respectively semidirect product of two of its subgroups.

**Lemma 1.23.** *$G$ is a direct product of two subgroups $H, K \leqslant G$ if and only if $H, K \trianglelefteq G$, $HK = G$, and $H \cap K = 1$.*

**Proof.** See for instance [DF04, 5.4.9]. □

**Lemma 1.24.** *$G$ is a semidirect product of two subgroups $H, K \leqslant G$ if and only if $H \trianglelefteq G$, $HK = G$, and $H \cap K = 1$.*

**Proof.** See for instance [DF04, 5.5.12]. □

The conditions $HK = G$ and $H \cap K = 1$ say that every element of $G$ is a unique product of an element in $H$ with an element in $K$. More on semidirect products of groups can be found for instance in [DF04, Chap. 5] or [Bou89, §6].

### 1.2.5 Automorphism groups

As described in Example 1.16, every group $G$ acts on itself by conjugation. The corresponding automorphisms $x \mapsto x^g$, $g \in G$, are called *inner* and the group of all inner automorphisms is denoted by $\mathrm{Inn}(G)$. Further, we denote the *outer automorphisms* $\mathrm{Aut}(G)/\mathrm{Inn}(G)$ by $\mathrm{Out}(G)$. By the first isomorphism theorem,

$$G/Z(G) \cong \mathrm{Inn}(G).$$

With a slight abuse of language we also refer to any automorphism that is not inner as an outer automorphism.

Suppose we have a group $G \cong \mathrm{Sym}_n$. Does it make sense to speak of a transposition in $G$ or does the notion of a transposition in $G$ depend on the isomorphism chosen between $G$ and $\mathrm{Sym}_n$? Since transpositions in $\mathrm{Sym}_n$ form a single conjugacy class they are stable under inner isomorphisms of $\mathrm{Sym}_n$. The question is therefore equivalent to asking if there are outer automorphisms of $\mathrm{Sym}_n$ that don't preserve transpositions. This turns out to be a nontrivial question and the answer provided by Theorem 1.27 shows that it does make sense to speak of transpositions in $G \cong \mathrm{Sym}_n$ whenever $n \neq 6$. The fact that the symmetric group $\mathrm{Sym}_6$ has outer automorphisms was first proved by Otto L. Hölder in [Höl95].

**Definition 1.25.** *A group $G$ such that $Z(G) = 1$ and $\mathrm{Inn}(G) = \mathrm{Aut}(G)$ is called* complete.

**Lemma 1.26.** *An automorphism of the symmetric group $\mathrm{Sym}_n$ is inner if and only if it preserves transpositions.*

An elementary proof can be found for instance in [Rot95, 7.4]. Here, we demonstrate how this result follows from a standard theorem on Kneser graphs, see Definition 1.35, which is included in the next chapter.

**Proof.** Let $\psi$ be an automorphism of $\mathrm{Sym}_n$ that preserves transpositions. Identify the transpositions of $\mathrm{Sym}_n$ with 2-subsets of $[n]$. Two transpositions commute if and only if their images under $\psi$ commute. Since transpositions $(s_1, s_2)$, $(t_1, t_2)$ commute if and only if they are disjoint as sets, we see that $\psi$ induces an automorphism of the Kneser graph $K(n, 2)$. Since $\mathrm{Sym}_n$ is generated by transpositions, the action of the transposition preserving automorphisms of $\mathrm{Sym}_n$ on $K(n, 2)$ is faithful. By Corollary 1.59 the automorphism group of $K(n, 2)$ is isomorphic to $\mathrm{Sym}_n$ provided that $n > 4$. Since the center of $\mathrm{Sym}_n$ is trivial this implies that for $n > 4$ all transposition preserving automorphisms are inner. The cases $n \leq 4$ are easily checked. $\square$

**Theorem 1.27.** $\mathrm{Sym}_n$ *is complete if and only if* $n \neq 2, 6$.

**Proof.** Clearly, the center of $\mathrm{Sym}_n$ is trivial whenever $n \neq 2$. Following [Rot95, 7.5], let $T_k$ denote the set of all permutations that are a product of $k$ disjoint transpositions. Observe that the sets $T_k$ partition the involutions of $\mathrm{Sym}_n$ into their conjugacy classes. Let $\psi$ be an automorphism of $\mathrm{Sym}_n$. Since $\psi$ permutes conjugacy classes we have $\psi(T_1) = T_k$ for some $k \leqslant n/2$. If we can show that $|T_1| \neq |T_k|$ for $k > 1$ then $\psi(T_1) = T_1$, that is $\psi$ preserves transpositions, and Lemma 1.26 implies that $\psi$ is inner.

By counting the number of ways to choose $k$ disjoint 2-subsets from $[n]$ and dividing by the number of their reorderings we find

$$|T_k| = \frac{1}{k!} \prod_{j=0}^{k-1} \binom{n-2j}{2} = \frac{1}{k! 2^k} n(n-1) \cdots (n-2k+1) = \frac{(n)_{2k}}{k! 2^k}$$

where $(x)_n$ denotes the *falling factorial* $x(x-1)\cdots(x-n+1)$. Consequently, $|T_1| = |T_k|$ if and only if

$$(n-2)_{2k-2} = k! 2^{k-1}.$$

Assume that $|T_1| = |T_k|$. Since $n \geqslant 2k$ we have $(n-2)_{2k-2} \geqslant (2k-2)!$ and by induction $(2k-2)! > k! 2^{k-1}$ whenever $k \geqslant 4$. Consequently, $k = 2$ or $k = 3$. Since $(n-2)_2 \neq 4$ not depending on the value of $n$, we are left with $k = 3$. If $n > 6$ then $(n-2)_4 \geqslant 5! = 120 > 24$ so that $n \geqslant 2k$ implies $n = 6$.

$\mathrm{Sym}_2$ has a nontrivial center, while $\mathrm{Sym}_6$ indeed admits outer automorphisms, a fact first proved in [Höl95]. Actually, $\mathrm{Aut}(\mathrm{Sym}_6) \cong \mathrm{Sym}_6 \rtimes \mathbb{Z}/2$ so that there essentially is one outer automorphism of $\mathrm{Sym}_6$. More details can be found in [Rot95]. $\square$

**Remark 1.28.** See also Remark 1.36 for a hint and motivation as to why $\mathrm{Sym}_6$ admits an outer automorphism.

## 1.3 Graph theory

In this section we shall briefly introduce the basic notions of graph theory. We refer to [Har94] or [GR01] for more details and thorough introductions. A *graph* is a set $\Gamma$ together with an irreflexive and symmetric binary relation $\perp$ on $\Gamma$. The elements of $\Gamma$ are called *vertices*, and $\perp$ is called the *adjacency relation* of the graph. Two vertices $x, y \in \Gamma$ are said to be *adjacent* if $x \perp y$. If $x$ and $y$ are adjacent then we refer to $\{x, y\}$ as an *edge*, and we denote the set of all edges by $E(\Gamma)$. We usually regard the adjacency relation $\perp$ as being implicitly given and identify the graph with its set of vertices $\Gamma$. In cases where this can lead to confusion we will write "$x \perp y$ in $\Gamma$" to emphasize that $\perp$ refers to the adjacency relation of the graph $\Gamma$. On the other hand, we write $V(\Gamma)$ when we wish to consider the vertices of the graph $\Gamma$ just as a set.

Let $\Gamma$ be a graph. The open *neighborhood* of $x \in \Gamma$ is the set

$$x^\perp \triangleq \{y \in \Gamma : y \perp x\}$$

of all vertices adjacent to $x$, and the vertices $x^\perp$ are called the neighbors of $x$. As $x \notin x^\perp$ we refer to $\{x\} \cup x^\perp$ as the *closed neighborhood* of $x$. Let $X \subseteq \Gamma$ be a set of vertices. We write

$$X^\perp \triangleq \bigcap_{x \in X} x^\perp$$

to denote the *common neighborhood* of $X$. A *path* in $\Gamma$ is a sequence $x_1, x_2, ..., x_n$ of vertices such that $x_i \perp x_{i+1}$ for $i \in [n-1]$ and all the $x_i$ are distinct. The path $x_1$, $x_2, ..., x_n$ is said to be of length $n$ connecting $x_1$ and $x_n$. The *distance* of two vertices $x$, $y$ is the minimal length of a path connecting $x$ and $y$. We write $d_\Gamma(x, y)$ for the distance of $x, y \in \Gamma$. Note that $d_\Gamma(x, y) = \infty$ if there is no path connecting $x$ and $y$. A graph $\Gamma$ is said to be *connected* if any two vertices can be connected by a path. If $X$, $Y \subseteq \Gamma$ are sets of vertices then the distance of $X$ and $Y$ is

$$d_\Gamma(X, Y) \triangleq \inf_{x \in X, y \in Y} d_\Gamma(x, y).$$

A *vertex-labeled* (respectively *edge-labeled*) graph is a graph $\Gamma$ together with a map $\mathrm{Lab}_\Gamma \colon \Gamma \to L$ (respectively $\mathrm{Lab}_\Gamma \colon E(\Gamma) \to L$) for some set of labels $L$. Such labelings can be visualized as colorings of the vertices (respectively edges) with $|L|$ colors. We will be particularly interested in the case of vertex-labeled graphs with $|L| = 2$ and we will call such a graph *bichromatic*. For reasons that will become clear later we will usually distinguish the vertices of a bichromatic graph as *short* versus *long* instead of being differently colored. When depicting a bichromatic graph we will draw long vertices as filled dots and short vertices as unfilled dots. If $\Gamma$ is a graph then we denote with $\Gamma^s$ (respectively $\Gamma^\ell$) the bichromatic graph obtained from $\Gamma$ by considering all vertices short (respectively long).

**Remark 1.29.** Recall that the adjacency relation $\perp$ of a graph $\Gamma$ was defined to be irreflexive which excludes the possibility of an edge connecting a vertex to itself. Further, $\perp$ was assumed to be symmetric which means that edges are undirected.

Let $\Gamma$ and $\Lambda$ be graphs. $\Lambda$ is said to be a *subgraph* of $\Gamma$ if $V(\Lambda) \subseteq V(\Gamma)$ and $E(\Lambda) \subseteq E(\Gamma)$. A subgraph $\Lambda$ of $\Gamma$ is said to *induced* if vertices $x$, $y$ of $\Lambda$ are adjacent in $\Lambda$ if and only if they are adjacent in $\Gamma$. We will often identify an induced subgraph with the set of its vertices. If $\Gamma$ is bichromatic then we will refer to the induced subgraph on the long (respectively short) vertices of $\Gamma$ as the long (respectively short) induced subgraph of $\Gamma$.

A map $\psi \colon \Gamma \to \Lambda$ between graphs $\Gamma$, $\Lambda$ is said to be a *homomorphism* if it preserves the adjacency relation, that is if

$$x \perp y \implies \psi(x) \perp \psi(y)$$

for all $x, y \in \Gamma$. A *homomorphism of labeled graphs* $\Gamma$, $\Lambda$ is required to also preserve the labeling meaning that $\mathrm{Lab}_\Gamma(x) = \mathrm{Lab}_\Gamma(y)$ implies $\mathrm{Lab}_\Lambda(\psi(x)) = \mathrm{Lab}_\Lambda(\psi(y))$ for every $x, y \in \Gamma$ (respectively $x, y \in E(\Gamma)$). In particular, a homomorphism of bichromatic graphs has to map short vertices to short vertices and long vertices to long vertices.

**Remark 1.30.** In other words, a graph homomorphism $\psi\colon \Gamma \to \Lambda$ is a mapping between vertices such that an edge is mapped to an edge, $\psi(E(\Gamma)) \subseteq E(\Lambda)$. Consequently, two adjacent vertices of $\Gamma$ can not be identified by $\psi$. In particular, graph homomorphisms between two given graphs don't need to exist. See for example the problem of colorability which is described in Remark 1.53.

As usual, a homomorphism $\psi\colon \Gamma \to \Lambda$ is called an isomorphism if there exists a homomorphism $\psi'\colon \Lambda \to \Gamma$ such that $\psi \circ \psi'$ and $\psi' \circ \psi$ are the identity on $\Lambda$ respectively $\Gamma$. Two graphs $\Gamma, \Lambda$ are said to be isomorphic, denoted by $\Gamma \cong \Lambda$, if an isomorphism between them exists. An isomorphism $\Gamma \to \Gamma$ is called an automorphism of $\Gamma$. We denote the group of all the graph automorphisms of $\Gamma$ by $\mathrm{Aut}(\Gamma)$.

We shall also consider maps between graphs which either map adjacent vertices to adjacent vertices or identify them and hence are only "almost" graph homomorphisms. Such a map will be called a partial homomorphism.

**Definition 1.31.** *A* partial graph homomorphism *is a map* $\psi\colon \Gamma \to \Lambda$ *such that*

$$x \perp y \quad \implies \quad \psi(x) \perp \psi(y) \ \lor \ \psi(x) = \psi(y)$$

*for all* $x, y \in \Gamma$.

An important example of partial graph homomorphisms are contractions which will be introduced in Definition 1.47.

The *complement* of a graph $X$ is the graph $\bar X$ on the vertices of $X$ with two vertices adjacent in $\bar X$ if they are not adjacent in $X$. It's easy to see that for any graph $\Gamma$ either $\Gamma$ or its complement $\bar\Gamma$ is connected.

### 1.3.1 Examples of graphs

**Example 1.32.** $C_n$ denotes the *circuit* of length $n$. Similarly, we denote by $C_\infty$ the infinite cycle, that is the graph on the integers with consecutive numbers adjacent.

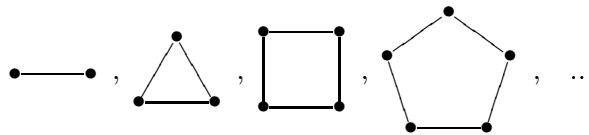

**Example 1.33.** $P_n$ denotes a *path* of length $n$.

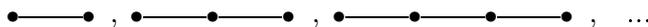

Note that a path in the previously defined sense, that is a sequence $x_1 \perp x_2 \perp \ldots \perp x_n$ of adjacent vertices, is equivalent to a (not necessarily induced) subgraph isomorphic to $P_n$.

**Example 1.34.** $K_n$ denotes the *complete graph* on $n$ vertices.

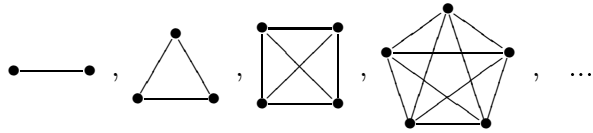

An (induced) subgraph that is a complete graph on $n$ vertices is also called an *n-clique*. We will usually refer to 2-cliques as edges and to 3-cliques as *triangles*. Likewise, an induced subgraph that is isomorphic to $\overline{K_n}$ is called an *n-coclique*. A coclique is also referred to as an *independent set* and the size $\alpha(\Gamma)$ of the largest independent set of a graph $\Gamma$ is called the *independence number* of $\Gamma$. Likewise, the size $\omega(\Gamma)$ of the largest clique of $\Gamma$ is called its *clique number*.

**Example 1.35.** The *Kneser graph* $K(n, k)$ is the graph on the $k$-element subsets of $[n]$ with two such subsets adjacent whenever they are disjoint.

Accordingly, $K(n, 1) \cong K_n$. Moreover, for even $n$ the Kneser graph $K(n, n/2)$ is a disjoint union of edges. If $k > n/2$ then $K(n, k)$ is a graph without edges. Consequently, we usually assume that $n \geqslant 2k + 1$ to exclude these two extremal cases. The graph $K(5, 2)$ is called the *Petersen graph* and is depicted below.

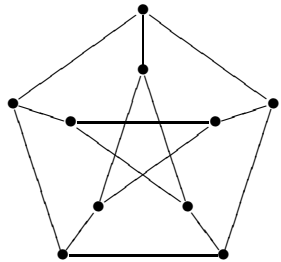

**Remark 1.36.** $K(6, 2)$ is a graph on 15 vertices each of which represents a 2-subset of $[6]$. Call a 3-clique of $K(6, 2)$ a line. We count that there are 15 lines and that each vertex is contained in exactly 3 lines. The following is an attempt to draw the vertices and lines of $K(6, 2)$.

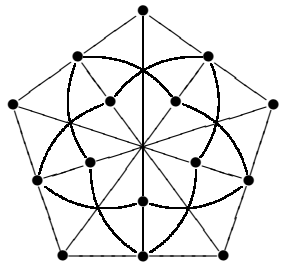

The vertices are the 15 dots, and the lines in this picture are the ten straight ones through three dots (five given by the sides of the outer pentagon and another five given by the lines through the barycenter) along with the five curved lines connecting two inner dots via an outer dot. By construction, the actual graph $K(6, 2)$ can be obtained from this picture by adding the invisible edge between the two "boundary points" of each line.

We will prove in Corollary 1.59 that the automorphism group of $K(6, 2)$ is $\text{Sym}_6$ acting on the underlying set $[6]$. Recall from Theorem 1.27 that $\text{Sym}_6$ is the only symmetric group which admits an outer automorphism, and from Lemma 1.26 that an outer automorphism necessarily doesn't preserve transpositions. Since a line represents a triple of pairwise disjoint 2-subsets of $[6]$ we count that $K(6, 2)$ contains 15 lines as well as 15 vertices. In fact, there exists a duality between points and lines in the associated geometry which gives rise to the outer automorphism of $\text{Sym}_6$. In particular, an outer automorphism of $\text{Sym}_6$ maps transpositions to products of three disjoint transpositions. See for instance [Rot95, 7.12] for an explicit construction of the outer automorphism of $\text{Sym}_6$. A nice account with further references is the article [JR82].

There are numerous ways to define graphs from a given group. In this text we will be especially interested in the following construction.

**Definition 1.37.** *Let $G$ be a group and $X \subseteq G$. The* commuting graph *of $G$ on $X$ is the graph with vertices $X$ in which two vertices $g, h \in X$ are adjacent whenever $g$ and $h$ commute.*

Note that $G$ acts on its commuting graph on $X$ by conjugation if and only if $X^g \subseteq X$ for all $g \in G$, that is if $X$ is a union of conjugacy classes of $G$. If $X$ consists of two conjugacy classes then we can regard the commuting graph on $X$ as a bichromatic graph by declaring vertices of one conjugacy class to be short and vertices of the other conjugacy class to be long.

**Example 1.38.** The complete graphs $K_n$ are the commuting graphs of the cyclic group $\mathbb{Z}/n$ on itself.

**Example 1.39.** The Kneser graph $K(n, 2)$ is the commuting graph of the symmetric group $\text{Sym}_n$ on its transpositions. This is because two distinct transpositions commute if and only if their support is disjoint.

To define another class of interesting graphs we have to introduce some terminology from linear algebra first. Let $\mathbb{F}$ be a field, and $V$ a vector space over $\mathbb{F}$. A *quadratic form* on $V$ is a mapping $Q \colon V \to \mathbb{F}$ such that $Q(\lambda v) = \lambda^2 \, Q(v)$ for all $\lambda \in \mathbb{F}$, $v \in V$, and such that the map $B \colon V \times V \to \mathbb{F}$ defined by

$$B(v, w) = Q(v + w) - Q(v) - Q(w)$$

is bilinear. In this case, $B$ is said to be the bilinear form associated to $Q$.

**Remark 1.40.** If the characteristic of $\mathbb{F}$ is different from 2 we can recover $Q$ from its associated bilinear form $B$ by

$$Q(v) = \frac{1}{2} \, B(v, v).$$

However, if the characteristic of $\mathbb{F}$ equals 2 then there may be distinct quadratic forms with the same associated bilinear form. For instance, let $V = \mathbb{F}_2^2$. The mappings $Q^+ \colon V \to \mathbb{F}_2$ and $Q^- \colon V \to \mathbb{F}_2$ defined by

$$Q^+(x) = x_1 x_2, \quad Q^-(x) = x_1^2 + x_2^2 + x_1 x_2,$$

are quadratic forms. For both of them the associated bilinear form is the standard symplectic form $B \colon V \times V \to \mathbb{F}_2$ given by

$$B(x,y) = x_1 y_2 + x_2 y_1.$$

In fact, up to isomorphism $Q^+$, $Q^-$ are the only quadratic forms that $B$ is associated to.

Let $V = \mathbb{F}_2^{2n}$, and consider the standard symplectic form $B \colon V \times V \to \mathbb{F}_2$ defined in Example 1.8. Let $Q^+ \colon V \to \mathbb{F}_2$ and $Q^- \colon V \to \mathbb{F}_2$ be the quadratic forms determined by

$$\begin{aligned} Q^+(x) &= x_1 x_2 + x_3 x_4 + \ldots + x_{n-1} x_n, \\ Q^-(x) &= x_1^2 + x_2^2 + x_1 x_2 + x_3 x_4 + \ldots + x_{n-1} x_n. \end{aligned}$$

Up to isomorphism $Q^+$ and $Q^-$ are the only quadratic forms that $B$ is associated to.

**Definition 1.41.** *Let $V = \mathbb{F}_2^{2n}$, and $B$ the standard symplectic form on V. The symplectic graph $\mathcal{S}p_2(2n)$ is the graph on the nonzero vectors $V \setminus \{0\}$ with two vectors $x, y$ adjacent whenever $B(x,y) = 0$.*

*Let $\varepsilon \in \{\,+\,,\,-\,\}$, and let $Q^\varepsilon$ be the quadratic form defined above. $\mathcal{NS}p^\varepsilon(2n)$ denotes the induced subgraph of $\mathcal{S}p_2(2n)$ on the vectors that are nonsingular under $Q^\varepsilon$ (that is those vectors $x$ such that $Q^\varepsilon(x) = 1$).*

A *polar subspace* of $V$ with respect to a bilinear form $B$ is a subset of the form

$$X^\pi \triangleq \{v \in V \colon (\forall x \in X) \ B(v,x) = 0\}$$

for some $X \subset V$. A polar subspace of $\mathcal{S}p_2(2n)$ is the induced subgraph on the vertices contained in a polar subspace of $\mathbb{F}_2^{2n}$ with respect to the standard symplectic form.

For more details about symplectic forms over $\mathbb{F}_2^{2n}$ and about the graphs just defined we refer the interested reader to [HS85, Section 2]. An implementation of the graphs $\mathcal{S}p_2(2n)$ and $\mathcal{NS}p^\varepsilon(2n)$ in SAGE that we will use for later computations can be found in the appendix.

### 1.3.2 Graph constructions

**Definition 1.42.** *(Sums of graphs) Let $\Gamma, \Lambda$ be graphs.*

- *The* disjoint union *$\Gamma \sqcup \Lambda$ is the graph with vertex set $V(\Gamma) \sqcup V(\Lambda)$ such that two vertices are adjacent if they are both contained and adjacent in one of $\Gamma$ or $\Lambda$. We write $n \cdot \Gamma$ for the disjoint union of $n$ copies of $\Gamma$.*

- *The* join $\Gamma + \Lambda$ *is the graph with vertex set* $V(\Gamma) \sqcup V(\Lambda)$ *such that two vertices are adjacent if they are adjacent in* $\Gamma \sqcup \Lambda$ *or if one of them is contained in* $\Gamma$ *and the other in* $\Lambda$.

Both the disjoint union $\Gamma \sqcup \Lambda$ and the join $\Gamma + \Lambda$ are made from a copy of $\Gamma$ and a copy of $\Lambda$. While in $\Gamma \sqcup \Lambda$ no elements from $\Gamma$ are adjacent to elements of $\Lambda$, in $\Gamma + \Lambda$ each element from $\Gamma$ is adjacent to each element of $\Lambda$.

Let $\Gamma$ be a graph. A maximal connected subgraph of $\Gamma$ is said to be a *component* of $\Gamma$. By construction, every graph is the disjoint union of its components. If $\Gamma$ is bichromatic then a short (respectively long) component of $\Gamma$ is a maximal connected subgraph of $\Gamma$ with only short (respectively only long) vertices. Equivalently, a short (respectively long) component of $\Gamma$ is a component of the short (respectively long) induced subgraph of $\Gamma$.

**Example 1.43.** Let $k \geqslant 2$ for notational reasons. The *complete multipartite graph* $K_{n_1,n_2,\ldots,n_k}$ is the graph on the disjoint union of $n_j$-element sets $V_{n_j}$ such that two vertices are adjacent whenever they lie in different $V_{n_j}$. In other words,

$$K_{n_1,\ldots,n_k} \triangleq \overline{K_{n_1}} + \ldots + \overline{K_{n_k}}.$$

Note that the complete graph $K_n$ is isomorphic to the complete multipartite graph $K_{1,1,\ldots,1}$ on $n$ vertices.

**Definition 1.44.** *(Products of graphs)* *Let* $\Gamma, \Lambda$ *be graphs.*

- *The* Cartesian product $\Gamma \square \Lambda$ *is the graph with vertex set* $V(\Gamma) \times V(\Lambda)$ *and*

  $$(x_1, y_1) \perp (x_2, y_2) \quad :\!\!\iff \quad (x_1 \perp x_2 \wedge y_1 = y_2) \vee (y_1 \perp y_2 \wedge x_1 = x_2).$$

- *The* composition $\Gamma[\Lambda]$ *is the graph with vertex set* $V(\Gamma) \times V(\Lambda)$ *and*

  $$(x_1, y_1) \perp (x_2, y_2) \quad :\!\!\iff \quad x_1 \perp x_2 \vee (y_1 \perp y_2 \wedge x_1 = x_2).$$

**Remark 1.45.** The composition $\Gamma[\Lambda]$ of two graphs $\Gamma, \Lambda$ is sometimes referred to as the wreath product $\Gamma \wr \Lambda$ of $\Gamma$ and $\Lambda$. A reason for this terminology is that

$$\mathrm{Aut}(\Gamma) \wr \mathrm{Aut}(\Lambda) \leqslant \mathrm{Aut}(\Gamma \wr \Lambda).$$

Moreover, equality holds in many cases.

**Lemma 1.46.** *If* $\Gamma, \Lambda$ *are connected graphs then their Cartesian product* $\Gamma \square \Lambda$ *is connected as well.*

**Proof.** Let $(x, y)$ and $(x', y')$ be two vertices of $\Gamma \square \Lambda$. By assumption, there exist paths $x \perp x_1 \perp x_2 \perp \ldots \perp x_n \perp x'$ and $y \perp y_1 \perp y_2 \perp \ldots \perp y_m \perp y'$ in $\Gamma$ respectively $\Lambda$. By definition of the Cartesian product $(x, y) \perp (x_1, y) \perp (x_2, y) \perp \ldots \perp (x_n, y) \perp (x', y)$ is a path connecting $(x, y)$ and $(x', y)$. Likewise, $(x', y) \perp (x', y_1) \perp (x', y_2) \perp \ldots \perp (x', y_n) \perp (x', y')$ is a path connecting $(x', y)$ and $(x', y')$. $\square$

**Definition 1.47.** *(Contractions) Let $\Gamma$ be a graph and $\Pi$ a partition of $\Gamma$. The contraction $\Gamma/\Pi$ is the graph on $\Pi$ such that two sets $A, B \in \Pi$ are adjacent whenever there is $a \in A$ and $b \in B$ such that $a \perp b$ in $\Gamma$.*

*Likewise, if $\Gamma$ is a bichromatic graph and $\Pi$ a partition of $\Gamma$ into sets of either short or long vertices then $\Gamma/\Pi$ is the bichromatic graph on $\Pi$ with adjacency as just defined and $A \in \Pi$ a short (respectively long) vertex of $\Gamma/\Pi$ if every $a \in A$ is a short (respectively long) vertex of $\Gamma$.*

**Remark 1.48.** The contraction of two vertices can be nicely visualized. Just imagine a graph drawn on a planar surface and move the two vertices closer and closer to each other. The contraction is what you get when these vertices come to rest on top of each other and are regarded as a single vertex. Analogously, one may visualize arbitrary contractions.

Note that the quotient map $\pi \colon \Gamma \to \Gamma/\Pi$ is a surjective partial homomorphism. It is a graph homomorphism if and only if any $W \in \Pi$ is a coclique of $\Gamma$. For $x \in \Gamma$ we denote with $[x]_\Pi$ the unique element of $\Pi$ containing $x$. We record the following observation.

**Proposition 1.49.** *Let $\Gamma$ be a graph and $\Pi$ a partition of $\Gamma$. A partial graph homomorphism $\psi \colon \Gamma \to \Gamma$ induces a partial graph homomorphism $\psi_\Pi$ on the contraction $\Gamma/\Pi$ given by $[x]_\Pi \mapsto [\psi(x)]_\Pi$ if and only if for any $X \in \Pi$ there is a $Y \in \Pi$ such that $\psi(X) \subseteq Y$. In this case, if $\Pi$ is finite and if $\psi$ is a graph automorphism then $\psi_\Pi$ is a graph automorphism.*

**Proof.** The condition that for every $X \in \Pi$ there is a (necessarily unique) $Y \in \Pi$ such that $\psi(X) \subseteq Y$ is equivalent to $\psi([x]_\Pi) \subseteq [\psi(x)]_\Pi$ for every $x \in \Gamma$. This in turn is exactly the well-definedness of the map

$$\psi_\Pi \colon \Pi \to \Pi, \quad [x]_\Pi \mapsto [\psi(x)]_\Pi.$$

To see that this map defines a partial graph homomorphism, let $X \perp Y \in \Gamma/\Pi$. By construction of $\Gamma/\Pi$ this means that we find $x \in X$ and $y \in Y$ such that $x \perp y$ in $\Gamma$. Since $\psi$ is a partial graph homomorphism we have $\psi(x) = \psi(y)$ or $\psi(x) \perp \psi(y)$ which implies that $[\psi(x)]_\Pi = [\psi(y)]_\Pi$ or $[\psi(x)]_\Pi \perp [\psi(y)]_\Pi$.

Additionally, assume that $\psi$ is an automorphism. Thus $\psi_\Pi$ is surjective which by finiteness of $\Pi$ implies that $\psi_\Pi$ is in fact a bijection. In particular, for each $X \in \Pi$ there is $Y \in \Pi$ such that $\psi(X) = Y$. Therefore $\psi$ and $\psi^{-1}$ induce partial graph homomorphisms $\psi_\Pi$ and $(\psi^{-1})_\Pi$. By construction, the compositions $\psi_\Pi \circ (\psi^{-1})_\Pi$ and $(\psi^{-1})_\Pi \circ \psi_\Pi$ are the identity on $\Gamma/\Pi$. We conclude that $\psi_\Pi$ is a graph automorphism. $\square$

**Example 1.50.** A particularly important and natural instance of a graph contraction is the following. Let $\Gamma$ be a graph, and let $\Pi$ be the partition of the vertices of $\Gamma$ into sets of vertices that have the same closed neighborhood. We will denote the graph $\Gamma/\Pi$ by $\Gamma^*$ and refer to it as the *reduced graph* of $\Gamma$. A graph $\Gamma$ is said to be reduced if $\Gamma = \Gamma^*$. Suppose that $\Gamma$ is finite. Then the conditions of Proposition 1.49 are met and we conclude that every graph automorphism of $\Gamma$ gives rise to a graph automorphism of $\Gamma^*$. In particular, whenever a group $G$ acts on $\Gamma$ then it also acts on $\Gamma^*$.

**Example 1.51.** Composing a graph with a complete graph can be seen as a partial converse to reducing a graph. Namely, if $\Gamma$ is a reduced graph then $(\Gamma[K_n])^* = \Gamma$.

### 1.3.3 Graph colorings

**Definition 1.52.** *Let $\Gamma$ be a graph. An $n$-coloring of $\Gamma$ is a mapping $\Gamma \to [n]$ assigning colors to the vertices of $\Gamma$. It is called* proper *if adjacent vertices are assigned distinct colors. $\Gamma$ is called $n$-colorable if there exists a proper $n$-coloring. The least value $n$ such that $\Gamma$ is $n$-colorable is called the* chromatic number $\chi(\Gamma)$ *of $\Gamma$.*

$\Gamma$ is called *bipartite* if its vertices are a disjoint union $V(\Gamma) = V_1 \sqcup V_2$ such that every edge connects a vertex from $V_1$ with a vertex from $V_2$ or, equivalently, if $\Gamma$ is a subgraph of a complete bipartite graph $K_{m,n}$. Notice that a graph is bipartite if and only if it is 2-colorable.

**Remark 1.53.** An $n$-coloring of $\Gamma$ can also be thought of as a graph homomorphism $\psi \colon \Gamma \to K_n$. To see this, recall that adjacent vertices of $\Gamma$ have to be mapped to adjacent vertices under the homomorphism $\psi$. Since the range of $\psi$ is the complete graph $K_n$ this condition is equivalent to adjacent vertices being mapped to distinct vertices in $K_n$ which resembles the definition of an $n$-coloring. For this reason, general graph homomorphisms $\psi \colon \Gamma \to \Lambda$ are also called $\Lambda$-colorings of $\Gamma$.

**Lemma 1.54.** *The chromatic number of the Cartesian product $\Gamma_1 \square \Gamma_2$ is*

$$\chi(\Gamma_1 \square \Gamma_2) = \chi(\Gamma_1) \vee \chi(\Gamma_2).$$

**Proof.** Let $f_i$ be proper colorings of $\Gamma_i$ using $\chi(\Gamma_i)$ colors, and set $\chi = \chi(\Gamma_1) \vee \chi(\Gamma_2)$. Define the $\chi$-coloring

$$f \colon \Gamma_1 \square \Gamma_2 \to \{0, 1, \ldots, \chi\}, \quad (x_1, x_2) \mapsto f_1(x_1) + f_2(x_2) \bmod \chi.$$

We claim that $f$ is a proper coloring of $\Gamma_1 \square \Gamma_2$. Let $(x_1, x_2)$ be a vertex of $\Gamma_1 \square \Gamma_2$. Neighbors of $(x_1, x_2)$ of the form $(\ldots, x_2)$ are contained alongside with $(x_1, x_2)$ in the induced subgraph $\Gamma_1 \square \{x_2\} \cong \Gamma_1$. This subgraph which we identify with $\Gamma_1$ can be properly colored by $f_1$ and accordingly by $f_1 + f_2(x_2) \bmod \chi$ as well. The same argument applies to neighbors of $(x_1, x_2)$ of the form $(x_1, \ldots)$. $\square$

### 1.3.4 Graph automorphisms

Let $\Gamma$ be a graph. Recall that we denote the group of all graph automorphisms $\Gamma \to \Gamma$ by $\mathrm{Aut}(\Gamma)$. Notice that a graph and its complement have the same automorphism group, that is $\mathrm{Aut}(\Gamma) = \mathrm{Aut}(\bar{\Gamma})$ for any graph $\Gamma$. We will only be interested in the automorphism groups of connected graphs $\Gamma$ which is justified by the following elementary result. This observation and a lot more information on automorphism groups of graphs constructed from simpler graphs can found for instance in [Har94, Chapter 14].

**Lemma 1.55.** *Let $\Gamma_1, ..., \Gamma_n$ be nonisomorphic connected graphs. Let $\Gamma$ be the disjoint union of $m_i$ copies of $\Gamma_i$. Then*

$$\mathrm{Aut}(\Gamma) \cong \mathrm{Aut}(\Gamma_1) \wr S_{m_1} \times \mathrm{Aut}(\Gamma_2) \wr S_{m_2} \times \cdots \times \mathrm{Aut}(\Gamma_n) \wr S_{m_n}. \qquad \square$$

By definition, $\mathrm{Aut}(\Gamma)$ acts on $\Gamma$. We often attribute properties of this action to the graph $\Gamma$ as well. For instance, we say that $\Gamma$ is *(vertex-)transitive* if $\mathrm{Aut}(\Gamma)$ acts transitively on the vertices of $\Gamma$.

Let $G$ be a group acting on $\Gamma$, and let $\Lambda$ be another graph. Note that $G$ also acts on the induced subgraphs of $\Gamma$ which are isomorphic to $\Lambda$. The next two examples present terminology associated to two special cases that will be of particular importance later on.

**Example 1.56.** If $G$ acts on $\Gamma$ then $G$ also acts on the triangles of $\Gamma$, that is on the induced subgraphs of $\Gamma$ that are isomorphic to $K_3$. If this induced action is transitive then we say that *$G$ acts transitively on the triangles* of $\Gamma$. As remarked before, we usually identify an induced subgraph with its set of vertices. A triangle of $\Gamma$ then corresponds to a set of three vertices $\{x_1, x_2, x_3\} \subset \Gamma$ which are pairwise adjacent. Let $\{y_1, y_2, y_3\}$ be another triangle of $\Gamma$. If $G$ acts transitive then there is $g \in G$ such that

$$g \cdot \{x_1, x_2, x_3\} = \{y_1, y_2, y_3\}.$$

Recall that $g \cdot \{x_1, x_2, x_3\} = \{g \cdot x_1, g \cdot x_2, g \cdot x_3\}$.

If $\{x_1, x_2, x_3\}$ is a triangle then we refer to the tuple $(x_1, x_2, x_3)$ as an *oriented triangle*. Clearly, $G$ acts on the oriented triangles of $\Gamma$ as well.

Suppose that $\Gamma$ is a bichromatic graph. Then we say that two oriented triangles $(x_1, x_2, x_3)$ and $(y_1, y_2, y_3)$ are of the same type if $x_i$ and $y_i$ are of the same type for all $i \in [3]$. An oriented triangle is said to be short (respectively long) if all its vertices are short (respectively long). Note that $G$ acts on the oriented triangles of $\Gamma$ of the same type.

**Example 1.57.** The terminology introduced in Example 1.56 for triangles analogously applies to paths. Recall that an $n$-path in $\Gamma$ is an induced subgraph that is isomorphic to $P_n$. If $\{x_1, x_2, ..., x_n\}$ is an $n$-path in $\Gamma$ then $(x_1, x_2, ..., x_n)$ is said to be an oriented $n$-path if $x_1 \perp x_2 \perp ... \perp x_n$. $G$ acts on (oriented) $n$-paths of $\Gamma$.

Again, if $\Gamma$ is bichromatic then we say that two oriented $n$-paths $(x_1, x_2, ..., x_n)$ and $(y_1, y_2, ..., y_n)$ are of the same type if $x_i$ and $y_i$ are of the same type for all $i \in [n]$. Clearly, $G$ acts on oriented $n$-paths of $\Gamma$ of the same type.

In the sequel, we wish to find the automorphism group $\mathrm{Aut}(K(n, k))$ of the Kneser graphs $K(n, k)$. As pointed out before, we may assume $n > 2k$ to exclude the extremal (and uninterestingly easy) cases. First, we observe that $\mathrm{Sym}_n \leqslant \mathrm{Aut}(K(n, k))$ induced by the action of $\mathrm{Sym}_n$ on the underlying set $[n]$. As it turns out we actually have equality, that is every automorphism of $K(n, k)$ is induced by a permutation of its underlying set. To prove this, we wish to make use of the simple fact that any automorphism of a graph $\Gamma$ has to permute its maximum independent sets. We therefore first classify the maximum independent sets of $K(n, k)$.

**Theorem 1.58. ([EKR61])** *Let $n \geqslant 2k$. Then*

$$\alpha(K(n,k)) = \binom{n-1}{k-1}.$$

*Further, if $n > 2k$ then for every maximum independent set $X$ of $K(n,k)$ there is $x \in [n]$ such that $X$ are those $k$-subsets of $[n]$ containing $x$.*

The beautiful proof we give for the first statement is due to Gyula O. H. Katona [Kat72]. A proof for the second claim that is close to ours can be found in [GR01]. Note that the case $n = 2k$ has to be excluded for the second statement since for instance $K(4,2)$ contains the maximum independent set $\{\{1,2\}, \{2,3\}, \{3,1\}\}$ which is not of the postulated form. A cyclic ordering of $[n]$ is a permutation $\psi \in \text{Sym}_{[n]}$ of order $n$ which can hence be interpreted as a successor function. The natural ordering corresponds to $(1, 2, ..., n)$. An interval with respect to a cyclic ordering $\psi$ is a set of successive elements, that is a set $\{a_1, ..., a_k\}$ such that $\psi(a_i) = a_{i+1}$ for $i \in \{1, 2, ..., k-1\}$. The right (respectively left) cyclic shift of a set $A \subseteq [n]$ with respect to $\psi$ is the set $\psi(A)$ (respectively $\psi^{-1}(A)$).

**Proof.** Let $n \geqslant 2k$ and let $X$ be a maximum independent set of $K(n,k)$. Consider the natural cyclic ordering of $[n]$. Without loss we may assume that $X$ contains the set $\{1, 2, ..., k\}$. For any $i \in \{1, 2, ..., k-1\}$ there are two intervals that contain exactly one of $i$ and $i+1$. Since $n \geqslant 2k$ these two intervals are disjoint, whence at most one of them is contained in $X$. Thus $X$ contains at most $k$ intervals. If we additionally assume that $X$ contains exactly $k$ intervals and that $n > 2k$ then these intervals necessarily share one element. Of course, the same arguments work for arbitrary cyclic orderings of $[n]$.

We now count in two ways the average number of intervals in $X$ with respect to some cyclic ordering. There are $(n-1)!$ cyclic orderings of $[n]$ and the previous argument showed that at most $k$ intervals per ordering can be contained in $X$. On the other hand, every subset of $X$ is an interval in exactly $k!(n-k)!$ many orderings. Hence

$$k \geqslant \frac{|X|\, k!\, (n-k)!}{(n-1)!}$$

and we conclude that

$$|X| \leqslant \binom{n-1}{k-1}.$$

Since those $k$-subsets of $[n]$ that contain a fixed $x \in [n]$ form independent sets of $K(n,k)$ we actually have equality in the preceding bound.

For proving the second statement suppose that $n > 2k$ and that $X$ is a maximum independent set. The above arguments imply that $X$ contains exactly $k$ intervals with respect to any cyclic ordering. Further, we observed that for a fixed ordering these $k$ intervals share an element. Considering the natural ordering we may thus assume without loss of generality that $X$ contains the sets

$$\{1, 2, ..., k\},\ \{2, 3, ..., k+1\},\ ...,\ \{k, k+1, ..., 2k-1\} \qquad (1.1)$$

which are those intervals sharing $k$. We are now going to show that $X$ actually contains any $k$-subset of $[n]$ which contains $k$.

To this end, consider a cyclic ordering of the form

$$(..., c, 1, 2, ..., 2k-1, ...)$$

for some $c \in \{2k, ..., n\}$. Since the sets from (1.1) are intervals with respect to this ordering as well, $X$ contains no further intervals. In particular, $X$ contains no sets of the form $\{c, 1, 2, ..., k-1\}$ for any $c \in \{2k, ..., n\}$. Now, let $c \in \{k+1, ..., 2k-1\}$ and consider the natural ordering with $c$ and $2k$ replaced. Recall that $2k < n$. Then $\{1, 2, ..., k\}$ is an interval but the left cyclic shift $\{n, 1, 2, ..., k-1\}$ is not contained in $X$. Thus $X$ contains the $k$ right cyclic shifts of $\{1, 2, ..., k\}$ analogous to (1.1). The same argument as before then shows that $X$ does not contain $\{c, 1, 2, ..., k-1\}$. Summarizing, we found that $X$ contains no sets of the form $\{c, 1, 2, ..., k-1\}$ for any $c \in \{k+1, ..., n\}$.

Let $A$ a $k$-subset of $[n]$ such that $k \in A$. Since $n > 2k$ we find an element $c \in \{k+1, ..., n\}$ that is not contained in $A$. We write $A = \{a_1, ..., a_{r-1}, k, a_{r+1}, ..., a_k\}$ where $a_1, ..., a_{r-1}$ are the elements of $A$ contained in $\{1, 2, ..., k-1\}$. Let $b_1, ..., b_{k-r}$ be the remaining elements of $\{1, 2, ..., k-1\}$, and consider a cyclic ordering of the form

$$(..., c, b_1, b_2, ..., b_{k-r}, a_1, ..., a_{r-1}, k, a_{r+1}, ..., a_k, ...)$$

Observe that $\{c, b_1, b_2, ..., b_{k-r}, a_1, ..., a_{r-1}\} = \{c, 1, 2, ..., k-1\}$ is not contained in $X$ while $\{b_1, b_2, ..., b_{k-r}, a_1, ..., a_{r-1}, k\} = \{1, 2, ..., k\}$ is. Therefore the $k$ right cyclic shifts of $\{1, 2, ..., k\}$ with respect to this ordering are contained in $X$. In particular, $A \in X$. □

**Corollary 1.59.** *Let $n > 2k$. Then $\mathrm{Aut}(K(n, k)) = \mathrm{Sym}_n$.*

**Proof.** For every $x \in [n]$ denote with $M(x)$ the set of $k$-subsets of $[n]$ containing $x$. According to Theorem 1.58 the sets $M(x)$ are exactly the maximum independent sets of $K(n, k)$. The automorphisms $\mathrm{Aut}(K(n, k))$ act on these $n$ maximum independent sets $M(x)$, and since this action is easily seen to be faithful we find that $\mathrm{Aut}(K(n, k)) \leqslant \mathrm{Sym}_n$. □

## 1.4 Coxeter groups

### 1.4.1 Basic properties

A Coxeter group is a group generated by *involutions* (that is elements of order 2) only subject to certain particularly simple relations. Starting with this abstract definition we will relate Coxeter groups to groups generated by reflections in Euclidean space. While we attempt to provide all necessary definitions and results (without proof), basic knowledge of Coxeter groups is strongly advised. Our presentation is heavily based on [Hum92] which we warmly recommend as an introduction. Other good references include [Coh08] and [BB05] as well as the classical [Bou02].

**Definition 1.60.** *A* Coxeter group *is a group $W$ with presentation*

$$W = \langle S | (st)^{m_{s,t}} = 1 \colon s, t \in S \rangle$$

*where $S$ is a finite set, $m_{s,t} \in \{1, 2, ..., \infty\}$, $m_{s,t} = m_{t,s}$, and $m_{s,t} = 1$ if and only if $s = t$. The pair $(W, S)$ is called a* Coxeter system. *Further, we refer to $S$ as the* Coxeter generators *and to $|S|$ as the* rank *of $(W, S)$.*

**Remark 1.61.** Our assumption that $S$ is finite is a common one although "a good part of the theory goes through for arbitrary $S$", see [Hum92]. However, we will mainly be interested in Coxeter systems $(W, S)$ for which in fact $W$ is finite.

Note that $m_{s,t} = 2$ if and only if $s$ and $t$ commute. The condition that $m_{s,s} = 1$ requires the element $s \in S$ to satisfy $s^2 = 1$. A priori we only know that the order of $st$ divides $m_{s,t}$ but according to Lemma 1.62 $m_{s,t}$ indeed is the order of $st$. In particular, the elements of $S$ are distinct involutions in $W$.

**Lemma 1.62.** *Let $(W, S)$ be a Coxeter system as in Definition 1.60. Then for any $s, t \in S$ the order of the product $st$ is $m_{s,t}$.*

**Proof.** See for instance [Hum92, Proposition 5.3] or [BB05, Proposition 1.1.1 and Corollary 1.4.8]. □

The values $m_{s,t}$ are therefore determined by the Coxeter system itself. There is a nice way to encode this data in a diagram which will prove particularly useful in the classification of the finite Coxeter groups.

**Definition 1.63.** *For a Coxeter system $(W, S)$ as in Definition 1.60 we define its associated* Coxeter diagram *as the edge-labeled graph on $S$ where $s$ and $t$ are connected whenever $m_{s,t} > 2$, in which case the label assigned is $m_{s,t}$.*

When drawing a Coxeter diagram we usually omit depicting labels of value 3. A Coxeter diagram is said to be *simply laced* if each edge has label 3. The *rank* of a Coxeter diagram is the number of its vertices which by construction coincides with the rank of the associated Coxeter system.

**Example 1.64.** Let $m \geqslant 3$. The dihedral group $\mathrm{Dih}_m$ can be represented as

$$\mathrm{Dih}_m = \langle s, t | s^2 = t^2 = (st)^m = 1 \rangle,$$

and hence is a Coxeter group. Note however that for odd $m$ we also have

$$\mathrm{Dih}_{2m} \cong \langle a, b, c | a^2 = b^2 = c^2 = (ac)^2 = (bc)^2 = (ab)^m = 1 \rangle.$$

This demonstrates that a Coxeter group may admit entirely different Coxeter systems. This example can be found for instance in [Bah05, 1.2.3].

We also allow $m = \infty$ and call the corresponding group $\text{Dih}_\infty = \langle s, t | s^2 = t^2 = 1 \rangle$ the infinite dihedral group.

**Remark 1.65.** A Coxeter group $W$ is said to be *rigid* if for every two Coxeter systems $(W, S)$ and $(W, S')$ there exists an automorphism $\psi \colon W \to W$ such that $\psi(S) = S'$. Equivalently, a Coxeter group $W$ is rigid if and only if all Coxeter systems $(W, S)$ have isomorphic associated Coxeter diagrams. Example 1.64 shows that the dihedral groups $\text{Dih}_m$ are not rigid for odd $m$. A good reference for rigidity of Coxeter groups is [Bah05].

In view of Lemma 1.62, we find that there is a one-to-one correspondence between Coxeter systems and Coxeter diagrams (up to isomorphism, of course). We therefore say that the Coxeter system $(W, S)$ is of *type $M$* where $M$ is the Coxeter diagram associated to $(W, S)$. With a slight abuse of language we also say that a Coxeter group $W$ is of type $M$ meaning that there exists a Coxeter system $(W, S)$ of type $M$.

**Lemma 1.66.** *Let $(W, S)$ be a Coxeter system, $T \subseteq S$, and set $W_T = \langle T \rangle$. Then $(W_T, T)$ is a Coxeter system. $W_T$ and its conjugates are said to be* parabolic subgroups *of $W$.*

**Proof.** See for instance [Hum92, Theorem 5.5]. $\square$

**Definition 1.67.** *A Coxeter system $(W, S)$ is said to be* reducible *if its Coxeter diagram is connected. Otherwise, it is called* irreducible.

**Remark 1.68.** Let $(W, S)$ be a Coxeter system, and $S = S_1 \sqcup S_2 \sqcup ... \sqcup S_n$ where the sets $S_i$ correspond to the connected components of the Coxeter diagram. Then the Coxeter systems $(W_{S_i}, S_i)$ are irreducible, and $W$ is the direct product of the parabolic subgroups $W_{S_i}$. For more details see for instance [Hum92, 6.1].

**Example 1.69.** Let's return to Example 1.64 and consider $\text{Dih}_m = \langle s, t | s^2 = t^2 = (st)^m = 1 \rangle$ where $m \geqslant 3$. The Coxeter system $(\text{Dih}_m, \{s, t\})$ is irreducible. However, for odd $m$ we also have the previously mentioned

$$\text{Dih}_{2m} \cong \langle a, b, c | a^2 = b^2 = c^2 = (ac)^2 = (bc)^2 = (ab)^m = 1 \rangle.$$

Consequently, the Coxeter system $(\text{Dih}_{2m}, \{a, b, c\})$ is reducible. This again illustrates the fact that different choices of Coxeter generators for the same Coxeter group give rise to essentially different Coxeter systems. Notice that the above example shows that in particular $\text{Dih}_{2m} \cong \text{Dih}_m \times \mathbb{Z}/2$ for odd $m$.

Let $n = |S|$ be the rank of the Coxeter system $(W, S)$. The *Coxeter matrix* $A = (a_{s,t}) \in \mathbb{R}^{S \times S}$ associated to $(W, S)$ is defined by

$$a_{s,t} \triangleq -\cos \frac{\pi}{m_{s,t}}$$

with the understanding that $1/\infty = 0$. By construction, $A$ is symmetric and hence induces a symmetric bilinear form $B$ on $V = \mathbb{R}^n$. Namely, after fixing a basis $\{\alpha_s\}_{s \in S}$ of $V$ we define the *Coxeter form* $B$ by setting

$$B(\alpha_s, \alpha_t) = a_{s,t}$$

and extending linearly. Note that $B(\alpha_s, \alpha_s) = 1$ while $B(\alpha_s, \alpha_t) \leqslant 0$ whenever $s \neq t$. We construct an action $\sigma$ of $W$ on $V$ by mapping $s \mapsto r_s$ where $r_s$ is the *generalized reflection* through $\alpha_s$ given by

$$r_s \colon v \mapsto v - 2B(\alpha_s, v)\alpha_s.$$

Notice that a generalized reflection indeed generalizes the usual notion of a reflection in Euclidean space in that it is a linear map which fixes a hyperplane pointwise and sends some nonzero vector to its negative (which in the Euclidean case is orthogonal to the hyperplane).

**Theorem 1.70.** *Let $(W, S)$ be a Coxeter system of rank $n$. Then the action $\sigma \colon W \to \mathrm{GL}(n, \mathbb{R})$ defined above is faithful.*

**Proof.** See for instance [Hum92, Corollary 5.4]. $\square$

We can therefore identify $W$ with the subgroup $\sigma(W)$ of $\mathrm{GL}(n, \mathbb{R})$. By construction, $W$ preserves the Coxeter form $B$ on $V$. Note that in general $B$ is not positive definite or even nondegenerate and hence does not induce a Euclidean geometry on $V$.

**Example 1.71.** Consider again the dihedral groups $\mathrm{Dih}_m = \langle s, t \,|\, s^2 = t^2 = (st)^m = 1 \rangle$ which have the Coxeter graph

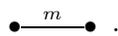

The associated Coxeter matrix is given by

$$\begin{pmatrix} 1 & -\cos\frac{\pi}{m} \\ -\cos\frac{\pi}{m} & 1 \end{pmatrix}$$

and is positive semidefinite. It is nondegenerate if and only if $m < \infty$.

Suppose that we have a finite Coxeter system $(W, S)$ with positive definite Coxeter matrix. Then $V$ endowed with the Coxeter form $B$ is isomorphic to the $n$-dimensional Euclidean space $\mathbb{R}^n$ with the standard inner product. Accordingly, the action $W \to \mathrm{GL}(n, \mathbb{R})$ identifies $W$ not only with a subgroup of $\mathrm{GL}(n, \mathbb{R})$ but with a subgroup of the orthogonal group $O(n, \mathbb{R})$ generated by reflections. Such a group is called a *reflection group*. The following result therefore shows that every finite Coxeter group is a finite reflection group. In fact, the converse is also true. A finite reflection can always be realized as a Coxeter group.

**Theorem 1.72.** *Let $(W, S)$ be a Coxeter system. $W$ is finite if and only if its Coxeter form is positive definite.*

**Proof.** See for instance [Hum92, Theorem 6.4]. □

**Theorem 1.73.** *Let $W$ be a finite reflection group. Then there exists a set of reflections $S \subset W$ such that $(W, S)$ is a Coxeter system.*

**Proof.** See for instance [Hum92, Theorem 1.9] or [Coh08, Theorem 5.1.4]. □

### 1.4.2 Finite Coxeter groups

In the sequel we will only be interested in finite Coxeter groups. According to Theorem 1.72 these are characterized by having a positive definite associated Coxeter form, and can hence be regarded as finite reflection groups. Since the Coxeter diagram is determined by the Coxeter form (and vice versa) we also call the Coxeter diagram positive definite whenever the Coxeter form is.

**Theorem 1.74.** *A Coxeter diagram is positive definite if and only if its connected components are isomorphic to one of the diagrams $A_n$ for $n \geqslant 1$, $B_n$ for $n \geqslant 2$, $D_n$ for $n \geqslant 4$, $E_6$, $E_7$, $E_8$, $F_4$, $G_2$, $H_3$, $H_4$ or $I_2(m)$ for $m \geqslant 3$ depicted below (the subindex corresponding to the number of vertices).*

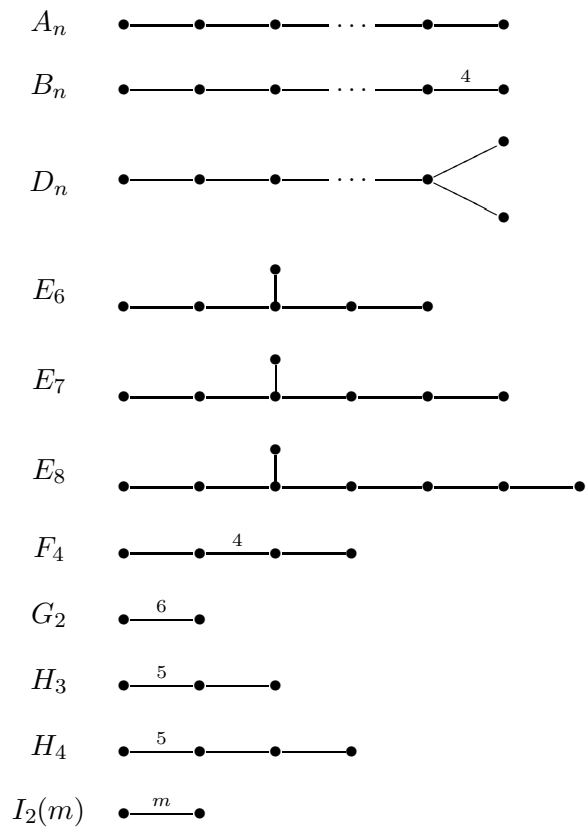

**Proof.** See for instance [Hum92, Theorem 2.7] or [Coh08, Theorem 5.3.3]. □

**Remark 1.75.** Coxeter diagrams corresponding to finite Coxeter groups are also called *spherical*. This terminology is motivated by the geometric realization provided by Theorem 1.70, and similarly there are Coxeter diagrams dubbed *affine* and *hyperbolic*, see for instance [Hum92, 6.5] and [Hum92, 6.8], which provide the two most important classes of infinite Coxeter groups. In this terminology, Theorem 1.74 classifies all spherical Coxeter diagrams.

Recall that we say that a Coxeter system is of type $M$ where $M$ is its associated Coxeter diagram. We have already seen that the Coxeter systems of type $I_2(m)$ are the dihedral groups generated by two reflections. In the next few examples we explicitly describe the other three infinite families of finite Coxeter systems, namely those of type $A_n$, $B_n$ and $D_n$.

**Example 1.76.** Consider the Coxeter diagram $A_n$

$$s_1 \text{------} s_2 \text{------} s_3 \text{------} \cdots \text{------} s_{n-1} \text{------} s_n \;.$$

The associated Coxeter system generated by $s_1, s_2, ..., s_n$ is defined by the presentation

$$\langle s_1, ..., s_n | (s_i s_j)^{m_{ij}} \colon i, j \in [n]\rangle$$

where (with $i,j$ considered modulo $n$)

$$m_{ij} = \begin{cases} 1 & i = j \\ 2 & |i-j| > 1 \\ 3 & |i-j| = 1 \end{cases}.$$

It is easy to verify that the symmetric group $\mathrm{Sym}_{n+1}$ with $s_i$ chosen to be the generating transpositions $(i\ i+1)$ satisfies this presentation. In fact, $(\mathrm{Sym}_{n+1}, \{s_1, s_2, ..., s_n\})$ is a Coxeter system of type $A_n$ a proof of which can be found for instance in [BB05, Proposition 1.5.4].

**Example 1.77.** Let $n \geqslant 2$. Denote with $\mathrm{Sym}_n^B$ the group of *signed permutations* of $[n]$, that is the group of permutations $\sigma$ of $-[n] \cup [n]$ such that

$$\sigma(-n) = -\sigma(n).$$

Let $\sigma_i$ be the permutation $(-i\ i)$ just changing the sign of $\pm i$, and let $s_i$ be the permutation $(-i\ -i-1)(i\ i+1)$ which corresponds to the transposition $(i\ i+1)$ of $[n]$. The $\sigma_i$ generate a group isomorphic to $(\mathbb{Z}/2)^n$ and the $s_i$ generate a group that we identify with $\mathrm{Sym}_n$. Since the sign changes $(\mathbb{Z}/2)^n$ are normal in $\mathrm{Sym}_n^B$ we obtain using Lemma 1.24 that $\mathrm{Sym}_n^B$ is the semidirect product

$$\mathrm{Sym}_n^B = (\mathbb{Z}/2)^n \rtimes \mathrm{Sym}_n = \mathbb{Z}/2 \wr \mathrm{Sym}_n.$$

$\mathrm{Sym}_n^B$ is generated by the signed transpositions $s_1, s_2, ..., s_{n-1}$ along with the single sign change $\sigma_n$. Again, one checks that the relations encoded in the Coxeter diagram

$$s_1 \text{—} s_2 \text{—} s_3 \text{—} \cdots \text{—} s_{n-1} \overset{4}{\text{—}} \sigma_n$$

of type $B_n$ are satisfied in $\mathrm{Sym}_n^B$. Some more work shows that indeed $(\mathrm{Sym}_n^B, \{s_1, s_2, ..., s_{n-1}, \sigma_n\})$ is a Coxeter system of type $B_n$ as is proved for instance in [BB05, Proposition 8.1.3].

**Example 1.78.** We adopt the notations from Example 1.77. Define $\mathrm{Sym}_n^D$ to be the subgroup of $\mathrm{Sym}_n^B = (\mathbb{Z}/2)^n \rtimes \mathrm{Sym}_n$ consisting of those permutations that involve only an even number of sign changes, that is the subgroup generated by the signed transpositions $s_1, s_2, ..., s_{n-1}$ along with the double sign change $\sigma_{n-1}\sigma_n$. By the same argument as for $\mathrm{Sym}_n^B$,

$$\mathrm{Sym}_n^D = (\mathbb{Z}/2)^{n-1} \rtimes \mathrm{Sym}_n.$$

Set $t = s_{n-1}\sigma_{n-1}\sigma_n$. It is straightforward to verify that $\mathrm{Sym}_n^D$ satisfies the relations encoded in the Coxeter diagram

$$s_1 \text{—} s_2 \text{—} s_3 \text{—} \cdots \text{—} s_{n-2} \diagup_{\displaystyle t}^{\displaystyle s_{n-1}}$$

of type $D_n$. Again, more work is needed to prove that $(\mathrm{Sym}_n^D, \{s_1, s_2, ..., s_{n-1}, t\})$ indeed is a Coxeter system of type $D_n$. A proof can be found for instance in [BB05, Proposition 8.2.3].

### 1.4.3 Crystallographic Coxeter groups

A *lattice* in $\mathbb{R}^n$ is a set $L$ consisting of the integer linear combinations of some vectors $v_1, v_2, ..., v_m$ which generate $\mathbb{R}^n$. In this case, we denote the lattice by $L = \mathbb{Z}\,v_1 + \mathbb{Z}\,v_2 + ... + \mathbb{Z}\,v_m$. Every lattice can in fact be written as $\mathbb{Z}\,w_1 + \mathbb{Z}\,w_2 + ... + \mathbb{Z}\,w_n$ for a basis $w_1, w_2, ..., w_n$ of $\mathbb{R}^n$ which we also refer to as a basis for the lattice, see [Gar97, 18.4]. The automorphisms of a lattice are the maps $g \in \mathrm{GL}(n, \mathbb{R})$ such that $gL \subseteq L$. In particular, a group $G \leqslant \mathrm{GL}(n, \mathbb{R})$ acts on $L$ if and only if $gL \subseteq L$ for every $g \in G$.

In this section we will be interested in Coxeter systems that correspond to subgroups of $O(n, \mathbb{R})$ acting on a lattice. Finite Coxeter systems of this kind are said to be *crystallographic*. The following together with the classification of finite Coxeter systems given in Theorem 1.74 allows us to determine the finite crystallographic Coxeter systems.

**Proposition 1.79.** *A finite Coxeter system $(W, S)$ as in Definition 1.60 is crystallographic only if $m_{s,t} \in \{1, 2, 3, 4, 6\}$ for all $s, t \in S$.*

**Proof.** See for instance [Hum92, Proposition 2.8]. □

Accordingly, Coxeter systems of type $H_3$, $H_4$, and $I_2(m)$ for $m \neq 3, 4, 6$ are not crystallographic. Since $I_2(3) = A_2$, $I_2(4) = B_2$ and $I_2(6) = G_2$ we are left with types $A_n$, $B_n$, $D_n$, $E_6$, $E_7$, $E_8$, $F_4$ and $G_2$ (note that this justifies the redundant introduction of the type $G_2$). For a Coxeter system of each of these types one can explicitly construct a lattice that it acts on. This is done for instance in [Hum92, 2.10]. We will present these lattices later on in terms of (not yet introduced) crystallographic root systems. Consequently, we have the following classification of crystallographic Coxeter systems.

**Corollary 1.80.** *An irreducible Coxeter system is crystallographic if and only if it is of type $A_n$, $B_n$, $D_n$, $E_6$, $E_7$, $E_8$, $F_4$ or $G_2$.*

In the sequel, we briefly sketch the concept of root systems and how they relate to the lattices that crystallographic Coxeter systems act on. More details can be found for instance in [Hum92, Chapter I]. Let $(W, S)$ be a finite Coxeter system. By Theorem 1.72 and Theorem 1.73 we can realize $W$ as a subgroup of $O(n, \mathbb{R})$ such that the involutions $S$ are orthogonal reflections. Recall that an orthogonal reflection $s$ is always of the form

$$s_\alpha \colon v \mapsto v - \frac{(v, \alpha)}{(\alpha, \alpha)} \alpha$$

where $\alpha$ is any vector on the line that $s$ reflects through. Let $R \subset W$ be the set of all reflections in $W$. If we choose vectors $\alpha_r$ for $r \in R$ such that $s_{\alpha_r} = r$ and such that the set

$$\Phi \triangleq \{\pm \alpha_r \colon r \in R\}$$

is left invariant under $W$ then $\Phi$ is called a *root system* and the vectors $\pm \alpha_r$ are called *roots* (note that we can always construct a root system by letting all the $\alpha_r$ have the same length). The roots $\alpha_s$ for $s \in S$ are called *simple roots* and we denote the set of all simple roots by $\Delta$. The simple roots $\Delta$ form a basis of $V$, and a root $\alpha$ is said to be *positive* (respectively *negative*) if it has nonnegative (respectively nonpositive) coordinates with respect to the basis $\Delta$. It turns out that every root is either positive or negative. $\Phi$ is said to be *irreducible* if $(W, S)$ is irreducible. Equivalently, $\Phi$ is irreducible if and only if there is no decomposition $\Phi = \Phi_1 \sqcup \Phi_2$ such that each root in $\Phi_1$ is orthogonal to each root in $\Phi_2$. Clearly, the Coxeter system $(W, S)$ is determined by the tuple $(\Phi, \Delta)$. Indeed, $W$ is the group generated by the reflections $s_\alpha$ for $\alpha \in \Phi$, and $S$ is the set of reflections $s_\alpha$ for $\alpha \in \Delta$. We say that $W$ is the *Weyl group* generated by the root system $\Phi$.

**Example 1.81.** The set $\Phi$ of vectors

$$\left\{ \begin{pmatrix} \cos \pi k/m \\ \sin \pi k/m \end{pmatrix} \colon k \in \{0, 1, ..., 2m-1\} \right\} \subset \mathbb{R}^2$$

forms a root system. $\Phi$ is irreducible provided that $m \neq 2$. The Weyl group $W$ generated by $\Phi$ is isomorphic to the dihedral group $\text{Dih}_m$ of order $2m$. Let $\alpha_1, \alpha_2$ be two roots of $\Phi$ which are at an angle of $\pi/m$, and let $s_1, s_2$ be the reflections through $\alpha_1$, $\alpha_2$. Then $(s_1 \, s_2)^m = 1$, and $(W, \{s_1, s_2\})$ is a Coxeter system of type $I_2(m)$. Consequently, the roots $\{\alpha_1, \alpha_2\}$ can be chosen as simple roots of $\Phi$.

Notice that if $m$ is even then

$$\left\{ \begin{pmatrix} \cos \pi k/m \\ \sin \pi k/m \end{pmatrix} : k \in \{0, 2, ..., 2m-2\} \right\} \cup \left\{ 2 \begin{pmatrix} \cos \pi k/m \\ \sin \pi k/m \end{pmatrix} : k \in \{1, 3, ..., 2m-1\} \right\}$$

is a different root system in $\mathbb{R}^2$ which generates the same Weyl group.

It turns out that if $(W, S)$ is a crystallographic Coxeter system than a lattice that $W$ acts on can be obtained from a well-chosen root system $\Phi$. Namely, there is a root system $\Phi$ such that $W$ acts on the lattice generated by the roots of $\Phi$. Such a root system is said to be *crystallographic*. Accordingly, the set $\Delta$ of simple roots is a basis for the lattice generated by $\Phi$. Let $\Phi$ be a crystallographic root system with associated Coxeter system $(W, S)$. If $(W, S)$ is of type $A_n$, $D_n$, $E_6$, $E_7$ or $E_8$ then the action of $W$ on the roots is transitive, and the roots therefore necessarily have the same length. However, if $(W, S)$ is of type $B_n$, $F_4$ or $G_2$ then there are two orbits of roots and the roots in one orbit are shorter than the roots in the other orbit. We accordingly refer to the roots of $\Phi$ as *short* versus *long* ones. Recall that the reflections of the Weyl group generated by $\Phi$ arise as reflections through the roots of $\Phi$. We refer to the reflections arising from short (respectively long) roots as *short (respectively long) root reflections*. An irreducible crystallographic root system $\Phi$ is (up to isomorphism) determined by the type of its associated Coxeter system together with the information which of the roots of $\Phi$ are short and which are long (in the case where there are only roots of one length it is common to call them long). We can incorporate this additional piece of information into the Coxeter diagram by coloring the vertices according to their length. Such a bichromatic edge-labeled graph is called a *Dynkin diagram*. We will say that a Dynkin diagram is *crystallographic* if it arises from a crystallographic root system.

The classification of crystallographic root systems is very close to the classification of finite Coxeter groups in Theorem 1.74. The only additional feature turns out to be that for a crystallographic root system there is a dual root system which reverses the role of short and long vertices. As mentioned above, crystallographic root systems of type $A_n$, $D_n$, $E_6$, $E_7$ or $E_8$ have roots of one length only whence the dual root system adds nothing new. We therefore denote the associated Dynkin diagrams by $A_n$, $D_n$, $E_6$, $E_7$ or $E_8$ as well. On the other hand, a crystallographic root system of type $B_n$, $F_4$ or $G_2$ contains short as well as long roots. For types $B_2$, $F_4$ and $G_2$ the dual root systems are again isomorphic to the original root systems, but for type $B_n$ we get a new family of crystallographic root systems whenever $n \geqslant 3$. We denote the associated Dynkin diagrams by $B_n$ and $C_n$. Summarizing, we have the following classification of connected crystallographic Dynkin diagrams.

**Theorem 1.82.** *A connected Dynkin diagram is crystallographic if and only if it is isomorphic to one of the Dynkin diagrams $A_n$ for $n \geqslant 1$, $B_n$ for $n \geqslant 2$, $C_n$ for $n \geqslant 2$, $D_n$ for $n \geqslant 4$, $E_6$, $E_7$, $E_8$, $F_4$, $G_2$ which are depicted below (the subindex corresponding to the number of vertices).*

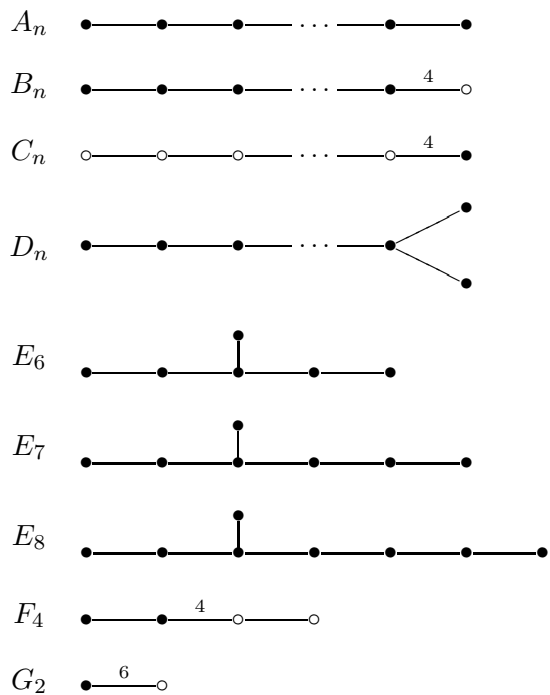

□

By forgetting the color of the vertices a Dynkin diagram can always be interpreted as a Coxeter diagram. No confusion will therefore result from treating Dynkin diagrams as Coxeter diagrams when appropriate.

Let $M$ be a connected crystallographic Dynkin diagram as classified in Theorem 1.82, and let $\Phi$ be a crystallographic root system of type $M$. Recall that the Weyl group $W$ generated by $\Phi$ can be turned into a Coxeter system of type $M$ once we choose a set of simple roots $\Delta \subset \Phi$. Further, recall that each reflection in $W$ is either a short or a long root reflection (with respect to $\Phi$). We denote the Weyl group $W$ together with the notion of a short (respectively long) root reflection by $W(M)$. Of course, this is only well-defined up to isomorphism.

**Remark 1.83.** Let $(W, S)$ be a finite irreducible crystallographic Coxeter system of rank at least 2 with corresponding crystallographic root system $\Phi$ and simple roots $\Delta$. We describe a somewhat surprising way to find the isomorphism type of the centralizer subgroup of a long root reflection $s_\alpha$, $\alpha \in \Phi$, which is described for instance in [Hum92, 2.11]. First, note that

$$C_W(s_\alpha) = C_W(\alpha).$$

Since $W$ acts transitively on long roots we may pick a particular one. To this end, we choose the unique highest root $\beta_0$ with respect to the partial order $\preccurlyeq$ defined by $\beta_1 \preccurlyeq \beta_2$ whenever all the coordinates of $\beta_2 - \beta_1$ with respect to $\Delta$ are nonnegative. Then $\beta_0$ is a long root, and $\beta_0 \notin \Delta$. Let $s_0$ be the reflection through $\beta_0$. As it turns out, the stabilizer of $s_0$ not just contains the simple root reflections it commutes with but is actually generated by them. In other words,

$$C_W(s_0) = \langle s_0 \rangle \times \langle s \in S \colon [s, s_0] = 1 \rangle.$$

According to Lemma 1.66 the centralizer $C_W(s_0)$ is a Coxeter group with Coxeter generators $s_0$ and the elements of $S$ commuting with $s_0$. Suppose we extend the Coxeter diagram to contain $s_0$. The simple reflections commuting with $s_0$ are exactly those not connected to it. Consequently, we obtain the Coxeter diagram for $C_W(s_0)$ by deleting the neighbors of $s_0$. We therefore have a way to read off the isomorphism type of $C_W(s_0)$ once we know how this extended Coxeter diagram looks like. Now, here is the sneaky part. This extended Coxeter diagram is just the corresponding (see for instance [Hum92, 2.5]) affine Coxeter diagram.

A general way to find the centralizer of any $s \in S$ in an arbitrary Coxeter system ($W$, $S$) based on [Bri96] can be found in [Bah05, 2.2.2].

We are now going to describe constructions for the irreducible crystallographic root systems classified in Theorem 1.82. These descriptions are based on [Hum92, 2.10] and [Bou02]. Throughout, let $\varepsilon_1, \varepsilon_2, ..., \varepsilon_n$ denote the standard basis of $\mathbb{R}^n$, and

$$L_n \triangleq \mathbb{Z}\varepsilon_1 + \mathbb{Z}\varepsilon_2 + ... + \mathbb{Z}\varepsilon_n$$

the *standard lattice*. We endow $\mathbb{R}^n$ with the standard inner product.

**Example 1.84.** Let $\Phi(A_n)$ be the vectors of squared length 2 in the standard lattice $L_{n+1}$ which are orthogonal to $\varepsilon_1 + \varepsilon_2 + ... + \varepsilon_{n+1}$. $\Phi(A_n)$ then consists of the $n(n+1)$ vectors

$$\varepsilon_i - \varepsilon_j, \quad 1 \leqslant i \neq j \leqslant n+1.$$

The Weyl group $W(A_n)$ generated by the root system $\Phi(A_n)$ is the symmetric group $\text{Sym}_{n+1}$ which faithfully acts on $\mathbb{R}^{n+1}$ by permuting the basis vectors $\varepsilon_1, \varepsilon_2, ..., \varepsilon_{n+1}$. In particular, the reflection through $\varepsilon_i - \varepsilon_j$ corresponds to the transposition $(i\ j)$ in $\text{Sym}_{n+1}$. Recall from Example 1.76 that $\text{Sym}_{n+1}$ is indeed a Coxeter group of type $A_n$.

**Example 1.85.** Let $\Phi(B_n)$ be the vectors of squared length 1 or 2 in the standard lattice $L_n$. Accordingly, $\Phi(B_n)$ consists of the $2n$ short roots

$$\pm\varepsilon_i, \quad 1 \leqslant i \leqslant n,$$

and the $2n(n-1)$ long roots

$$\pm\varepsilon_i \pm \varepsilon_j, \quad 1 \leqslant i \neq j \leqslant n.$$

The Weyl group $W(B_n)$ generated by $\Phi(B_n)$ is the signed permutation group $\text{Sym}_n^B = (\mathbb{Z}/2)^n \rtimes \text{Sym}_n$ which acts faithfully on $\mathbb{R}^n$ by changing the signs of and permuting the basis vectors $\varepsilon_1, \varepsilon_2, ..., \varepsilon_n$. Using the notation of Example 1.77, the short root reflection through $\varepsilon_i$ corresponds to the sign change $\sigma_i$ in $\text{Sym}_n^B$, and the long root reflection through $\varepsilon_i - \varepsilon_{i+1}$ corresponds to the transposition $s_i$. Recall that $\text{Sym}_n^B$ is indeed a Coxeter group of type $B_n$.

**Example 1.86.** $\Phi(C_n)$ is the dual root system of $\Phi(B_n)$ which is obtained by switching the role of short and long roots. $\Phi(C_n)$ consists of the $2n(n-1)$ short roots

$$\pm\varepsilon_i \pm \varepsilon_j, \quad 1 \leqslant i \neq j \leqslant n,$$

and the $2n$ long roots

$$\pm 2\varepsilon_i, \quad 1 \leqslant i \leqslant n.$$

By construction, the Weyl groups generated by the root systems $\Phi(C_n)$ and $\Phi(B_n)$ are identical. However, a long root reflection of the Weyl group $W(C_n)$ is a short root reflection of $W(B_n)$ while a short root reflection of $W(C_n)$ is a long root reflection of $W(B_n)$.

**Example 1.87.** Let $\Phi(D_n)$ be the vectors of squared length 2 in the standard lattice $L_n$. Consequently, $\Phi(B_n)$ consists of the $2n(n-1)$ roots

$$\pm \varepsilon_i \pm \varepsilon_j, \quad 1 \leqslant i \neq j \leqslant n.$$

Observe that $\Phi(D_n) \subset \Phi(B_n)$. Analogous to Example 1.85 one finds that the Weyl group $W(D_n)$ generated by $\Phi(D_n)$ is isomorphic to $\mathrm{Sym}_n^D = (\mathbb{Z}/2)^{n-1} \rtimes \mathrm{Sym}_n$ introduced in Example 1.78 as a Coxeter group of type $D_n$.

**Example 1.88.** Let $L'$ be the sublattice of $L_8$ given by vectors $c_1 \varepsilon_1 + c_2 \varepsilon_2 + \ldots + c_8 \varepsilon_8$ such that the sum of the $c_i$ is even. Consider the lattice

$$L \triangleq L' + \mathbb{Z}\left(\frac{1}{2}(\varepsilon_1 + \varepsilon_2 + \ldots + \varepsilon_8)\right).$$

Let $\Phi(E_8)$ be the vectors of squared length 2 in $L$. $\Phi(E_8)$ then consists of the $112 + 128 = 240$ roots

$$\begin{aligned}\pm \varepsilon_i \pm \varepsilon_j, \quad & 1 \leqslant i \neq j \leqslant 8,\\ \tfrac{1}{2}\left(\sum c_i \varepsilon_i\right), \quad & c_i = \pm 1, \prod c_i = 1.\end{aligned}$$

Let $\Phi(E_7)$ be the roots of $\Phi(E_8)$ orthogonal to $\varepsilon_7 + \varepsilon_8$. $\Phi(E_7)$ consists of the $60 + 2 + 64 = 126$ roots

$$\begin{aligned}\pm \varepsilon_i \pm \varepsilon_j, \quad & 1 \leqslant i \neq j \leqslant 6,\\ \pm(\varepsilon_7 - \varepsilon_8), & \\ \tfrac{1}{2}\left(\sum c_i \varepsilon_i\right), \quad & c_i = \pm 1, \prod c_i = 1, c_7 = -c_8.\end{aligned}$$

Let $\Phi(E_6)$ be the roots of $\Phi(E_7)$ orthogonal to $\varepsilon_6 + \varepsilon_8$. Consequently, $\Phi(E_6)$ consists of the $40 + 32 = 72$ roots

$$\begin{aligned}\pm \varepsilon_i \pm \varepsilon_j, \quad & 1 \leqslant i \neq j \leqslant 5,\\ \tfrac{1}{2}\left(\sum c_i \varepsilon_i\right), \quad & c_i = \pm 1, \prod c_i = 1, c_6 = c_7 = -c_8.\end{aligned}$$

**Example 1.89.** Consider the lattice

$$L \triangleq L_n + \mathbb{Z}\left(\frac{1}{2}(\varepsilon_1 + \varepsilon_2 + \varepsilon_3 + \varepsilon_4)\right).$$

Let $\Phi(F_4)$ be the vectors of squared length 1 or 2 in $L$. $\Phi(F_4)$ consists of the $8+16=24$ short roots

$$\pm \varepsilon_i, \quad 1 \leqslant i \leqslant 4, \qquad \frac{1}{2}(\pm \varepsilon_1 \pm \varepsilon_2 \pm \varepsilon_3 \pm \varepsilon_4),$$

and the 24 long roots

$$\pm \varepsilon_i \pm \varepsilon_j, \quad 1 \leqslant i \neq j \leqslant 4.$$

**Example 1.90.** Let $\Phi(G_2)$ be the vectors of squared length 2 or 6 in the standard lattice $L_3$ which are orthogonal to $\varepsilon_1 + \varepsilon_2 + \varepsilon_3$. Consequently, $\Phi(G_2)$ consists of the 6 short roots

$$\varepsilon_i - \varepsilon_j, \quad 1 \leqslant i \neq j \leqslant 3,$$

and the 6 long roots

$$\pm(2\varepsilon_i - \varepsilon_j - \varepsilon_k), \quad \{i,j,k\} = [3].$$

In the sequel we will mostly work with the crystallographic root systems $\Phi(A_n)$, $\Phi(B_n)$, $\Phi(C_n)$, $\Phi(D_n)$, $\Phi(F_4)$. These are summarized below.

| | | |
|---|---|---|
| $A_n$ | $\varepsilon_i - \varepsilon_j$ | $1 \leqslant i \neq j \leqslant n+1$ |
| $B_n$ | $\pm \varepsilon_i \pm \varepsilon_j$ | $1 \leqslant i \neq j \leqslant n$ |
| | $\pm \varepsilon_i$ | $1 \leqslant i \leqslant n$ |
| $C_n$ | $\pm \varepsilon_i \pm \varepsilon_j$ | $1 \leqslant i \neq j \leqslant n$ |
| | $\pm 2\varepsilon_i$ | $1 \leqslant i \leqslant n$ |
| $D_n$ | $\pm \varepsilon_i \pm \varepsilon_j$ | $1 \leqslant i \neq j \leqslant n$ |
| $F_4$ | $\pm \varepsilon_i \pm \varepsilon_j$ | $1 \leqslant i \neq j \leqslant 4$ |
| | $\pm \varepsilon_i$ | $1 \leqslant i \leqslant 4$ |
| | $(\pm \varepsilon_1 \pm \varepsilon_2 \pm \varepsilon_3 \pm \varepsilon_4)$ | |

### 1.4.4 Relations between Coxeter groups

By construction, as given in Example 1.78, the Weyl group $W(D_n)$ is embedded in $W(B_n)$ as a subgroup of index 2. In fact, recall from Example 1.77 that $W(B_n) = (\mathbb{Z}/2)^n \rtimes \mathrm{Sym}_n$ acts faithfully on $\mathbb{R}^n$ by permuting and changing the signs of the standard basis vectors $\varepsilon_i$. As in Example 1.77 let $s_i$ permute $\varepsilon_i$ and $\varepsilon_{i+1}$, and let $\sigma_i$ be the sign change on $\varepsilon_i$. Recall that

$$W(B_n) = \langle s_1, s_2, \ldots, s_{n-1}, \sigma_n \rangle,$$

and

$$W(D_n) = \langle s_1, s_2, \ldots, s_{n-1}, s_{n-1}\sigma_{n-1}\sigma_n \rangle.$$

Consequently, $W(B_n) = W(D_n) \rtimes \langle \sigma_n \rangle$ according to Lemma 1.24.

**Lemma 1.91.** $W(B_n) \cong W(D_n) \rtimes \mathrm{Sym}_2$. □

We now wish to illustrate that this relation can also be seen by looking at the Dynkin diagram of $W(D_n)$.

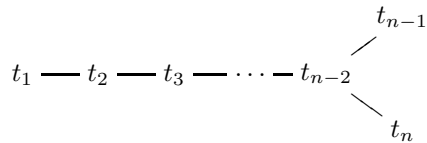

Observe that this diagram has precisely one nontrivial automorphism for $n \geqslant 5$. This automorphism has order 2 and permutes the two vertices $t_{n-1}$ and $t_n$. Since $W(D_n)$ is the Coxeter group generated by $t_1, t_2, \ldots, t_n$ this graph automorphism lifts to a unique automorphism $\psi$ of $W(D_n)$. Consider the group $W(D_n) \rtimes \langle \psi \rangle$ which is generated by $t_1, t_2, \ldots, t_{n-1}, \psi$. By definition, $\psi$ centralizes $t_1, \ldots, t_{n-2}$, and $\psi t_{n-1} = t_n \psi$. Consequently, we have $(t_{n-1} \psi)^2 = t_{n-1} t_n$ and thus $(t_{n-1} \psi)^4 = 1$. This shows that $t_1$, $t_2, \ldots, t_{n-1}, \psi$ satisfy the defining Coxeter relations of type $B_n$. Taking into account the number of elements we deduce that $W(D_n) \rtimes \langle \psi \rangle \cong W(B_n)$.

This way, we also get an understanding why $W(B_n)$ contains two conjugacy classes of reflections. One class are the conjugates of the reflections $t_i$ coming from $W(D_n)$ which by construction are left invariant under the action of $\psi$, while the second class of conjugates are the reflections conjugate to $\psi$.

While for $n \geqslant 5$ there is just the one observed automorphism, the Coxeter diagram of type $D_4$ is more symmetric with $\mathrm{Sym}_3$ as its automorphism group.

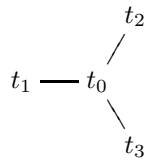

This phenomenon is sometimes referred to as *triality*. Let $\psi_1$ be the automorphism permuting $t_1$ and $t_2$, and $\psi_2$ the automorphism permuting $t_2$ and $t_3$. Then $\langle \psi_1, \psi_2 \rangle \cong \mathrm{Sym}_3$ and $W(D_4) \rtimes \langle \psi_1, \psi_2 \rangle$ is generated by $t_0, t_1, \psi_1, \psi_2$. By definition, $\psi_1$ centralizes $t_1$, and both $\psi_1$ and $\psi_2$ centralize $t_0$. Further, we have $\psi_1 t_1 = t_2 \psi_1$ whence, as above, $(t_1 \psi_1)^2 = t_1 t_2$ and $(t_1 \psi_1)^4 = 1$. In other words, $t_0, t_1, \psi_1, \psi_2$ satisfy the relations encoded into the Coxeter diagram

$$t_0 \text{———} t_1 \overset{4}{\text{———}} \psi_1 \text{———} \psi_2 \ .$$

This diagram is of type $F_4$. Accordingly, $W(D_4) \rtimes \langle \psi_1, \psi_2 \rangle$ is a quotient of $W(F_4)$. Since $|W(D_4) \rtimes \mathrm{Sym}_3| = 2^{4-1} \, 4! \, 3! = 1152 = |W(F_4)|$ the number of elements of these groups coincide. Consequently, we have proved the following.

**Lemma 1.92.** $W(F_4) \cong W(D_4) \rtimes \mathrm{Sym}_3$. □

Again, this construction reveals that the Coxeter generators of $W(F_4)$ belong to two distinct conjugacy classes, namely the one coming from $W(D_4)$ made up by conjugates of the $t_i$ and the second consisting of the conjugates of the $\psi_i$.

# 2  Local recognition of graphs

Let $\Gamma$ be a graph. Recall that if $x \in \Gamma$ then $x^\perp$ denotes the set of vertices adjacent to $x$. The induced subgraph on $x^\perp$ is called the *local graph* at $x$. A graph $\Gamma$ is said to be *locally homogeneous* if all its local graphs are isomorphic, that is if there exists a graph $\Delta$ such that the induced subgraph on $x^\perp$ is isomorphic to $\Delta$ for all $x \in \Gamma$. In this case, $\Gamma$ is said to be locally $\Delta$, and $\Delta$ is referred to as the local graph of $\Gamma$. If $\Gamma$ is locally homogeneous then we denote its local graph by $\Delta(\Gamma)$. Of course, this is well-defined only up to isomorphism. If $\Lambda$ is another locally homogeneous graph such that $\Delta(\Lambda) \cong \Delta(\Gamma)$ then we say that $\Lambda$ is *locally like* $\Gamma$.

**Remark 2.1.** In literature, some authors refer to the local graph at a vertex $x$ as the *link* of $x$. Accordingly, locally homogeneous graphs are also referred to as graphs with constant link. This terminology stems from the fact that, when considering a graph as a simplicial complex, the notion of the local graph at a vertex coincides with the already existing notion of the link of a vertex.

Notice that a graph $\Gamma$ is locally $\Delta$ if and only if all its connected components are locally $\Delta$. Accordingly, we usually restrict the discussion to connected graphs.

**Remark 2.2.** Which graphs occur as local graphs? This question was posed by Alexander A. Zykov in [Zyk64] and is commonly referred to as the *Trahtenbrot-Zykov problem*. An introduction to this problem can be found for instance in the first part of [BC75]. An excellent resource for techniques for proving whether a graph is a local graph is [Hal85]. Moreover, [Hal85] describes the graphs of order up to 6 that are local graphs.

Is there an algorithm to decide which graphs are local graphs? It is proved in [Bul73] that this question is algorithmically unsolvable when infinite graphs are allowed. When restricting to finite graphs this question is still open as is noted in [Bug90]. To illustrate that it makes a difference to permit infinite graphs we refer to [Hal85, 4.16] where it is proved that there is no finite graph $\Gamma$ that is locally $K_{1,3} \sqcup \overline{K_2}$ while infinite graphs are constructed which are locally $K_{1,3} \sqcup \overline{K_2}$.

Here, we will be interested in the problem of characterizing a locally homogeneous graph in terms of its local graph. We say that a connected locally homogeneous graph $\Gamma$ is *locally recognizable* if up to isomorphism $\Gamma$ is the only connected graph that is locally $\Delta(\Gamma)$. Lemma 2.3 shows that the complete graphs $K_n$ are locally recognizable. On the other hand, Example 2.4 demonstrates that the circuit graphs $C_n$ are not locally recognizable.

**Lemma 2.3.** *A connected graph $\Gamma$ is locally $K_{n-1}$ if and only if $\Gamma \cong K_n$.*

**Proof.** By assumption, every vertex of $\Gamma$ has $n-1$ neighbors each two of which are adjacent. Let $x$ and $y$ be adjacent vertices. We find that the neighbors of $y$ are given by $x$ together with the $n-2$ neighbors of $x$ besides $y$. Since $\Gamma$ is connected this implies that the closed neighborhood of $x$ exhausts the vertices of $\Gamma$. □

**Example 2.4.** The circuit $C_n$ is locally homogeneous. If $n \geqslant 4$ its local graph is $\overline{K_2}$, the graph with two nonadjacent vertices.

Thus, in contrast to the complete graphs $K_n$, the circuit graphs $C_n$ are not locally recognizable for $n \geqslant 4$ (for $n < 4$ we get the complete graphs again). However, the local structure still conveys information.

**Lemma 2.5.** *A connected graph $\Gamma$ is locally $\overline{K_2}$ if and only if $\Gamma \cong C_n$ for some $n \in \{4, 5, ..., \infty\}$.*

**Proof.** Starting with a vertex $x_1 \in \Gamma$, we have two neighbors $x_0$, $x_2$ which are not adjacent. $x_2$ therefore has another neighbor $x_3$. Either $x_3$ is adjacent to $x_0$ (whence $\Gamma \cong C_4$) or it is adjacent to another vertex $x_4$. Again, $x_4$ is either adjacent to $x_0$ (whence $\Gamma \cong C_5$) or we find yet another vertex $x_5$, and so on.

$$\cdots\cdots x_0 \text{—} x_1 \text{—} x_2 \text{—} x_3 \text{—} x_4 \cdots\cdots$$

If this process doesn't terminate then $\Gamma \cong C_\infty$. □

In other words, all the circuits locally look the same but faced with a graph $\Gamma$ that locally looks like a circuit we are still able to recognize it as a circuit. We just can't tell how big a circuit $\Gamma$ is.

**Lemma 2.6.** *Let $\Gamma_1$ be locally $\Delta_1$, and $\Gamma_2$ be locally $\Delta_2$. Then the Cartesian product $\Gamma_1 \square \Gamma_2$ is locally $\Delta_1 \sqcup \Delta_2$.*

**Proof.** Let $(x_1, x_2) \in \Gamma_1 \square \Gamma_2$. By definition, a vertex $(y_1, y_2)$ of the Cartesian product is adjacent to $(x_1, x_2)$ if and only if

$$(x_1 \sim y_1 \wedge x_2 = y_2) \vee (x_2 \sim y_2 \wedge x_1 = y_1).$$

Hence the neighborhood of $(x_1, x_2)$ in $\Gamma_1 \square \Gamma_2$ is isomorphic to the disjoint union of the neighbors of $x_1$ in $\Gamma_1$ and the neighbors of $x_2$ in $\Gamma_2$. □

Consequently, the class of all local graphs is closed under disjoint unions. In other words, a graph is a local graph if all of its connected components are local graphs.

**Remark 2.7.** If a graph $\Gamma$ is transitive then $\Gamma$ is locally homogeneous. Not surprisingly, the converse is not true. The smallest graphs which are locally homogeneous but not transitive have 10 vertices, and up to isomorphism there are exactly three such graphs, depicted below.

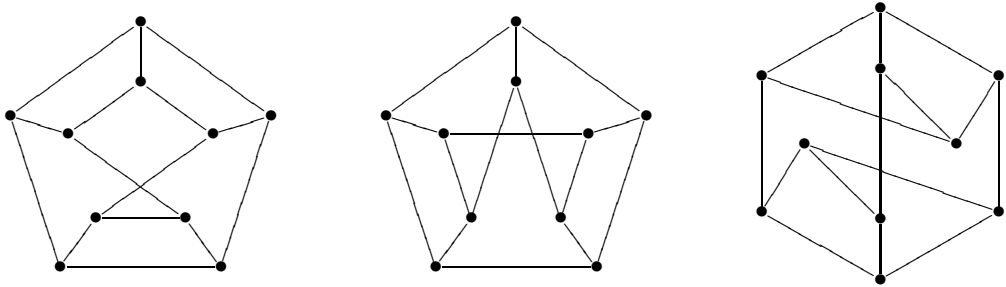

These graphs and further information can be found in [Hal85].

## 2.1 Some local recognition results

Notice that the Kneser graph $K(n,k)$ is locally homogeneous with local graph $K(n-k,k)$. Jonathan I. Hall proves in [Hal87] that for $n$ sufficiently large compared to $k$ the Kneser graphs are locally recognizable.

**Theorem 2.8.** ([Hal87]) *Let $k \geqslant 1$, $n \geqslant 3k+1$. A connected graph $\Gamma$ is locally $K(n,k)$ if and only if $\Gamma \cong K(n+k,k)$.* □

Define $n(k)$ to be the minimal $n$ such that a locally $K(n,k')$ graph can be recognized as $K(n+k',k')$ for any $k' \geqslant k$. According to Theorem 2.8 $n(k) \leqslant 3k+1$. On the other hand, we clearly have $n(k) \geqslant 2k+1$ since the graphs $K(2k,k)$ are degenerate. There are 3 graphs that are locally $K(6,2)$, see Theorem 2.10, whence $n(2) = 7$. However, it might still be possible that $n(k) = 2k+1$ for all $k \geqslant 3$ as is noted in [Hal87]. In particular, it is not even yet known if there are further connected graphs that are locally $K(7,3)$.

**Theorem 2.9.** ([Hal80]) *A connected graph $\Gamma$ is locally $K(5,2)$ if and only if $\Gamma$ is isomorphic to one of the graphs*

- $K(7,2)$,
- $3 \cdot K(7,2)$,
- $\Sigma L_{2,25}$.

*In particular, $|\Gamma| \in \{21, 63, 65\}$.* □

Here, the graph $3 \cdot K(7,2)$ is the 3-fold cover of $K(7,2)$, and $\Sigma L_{2,25}$ is the graph on the conjugates of the unique nontrivial field automorphism of $\mathbb{F}_{25}$ in the special semi-linear group $\Sigma L(2, 25)$ with two elements adjacent whenever they commute, see Example 1.21. For more details about these graphs we refer to [Hal80].

**Theorem 2.10.** ([HS85]) *A connected graph $\Gamma$ is locally $K(6,2)$ if and only if $\Gamma$ is isomorphic to one of the graphs*

- $K(8,2)$,

- $\mathcal{S}p_2(6)$ minus $\{x\} \cup x^\perp$ for some $x$,
- $\mathcal{N}\mathcal{S}p^-(6)$.

In particular, $|\Gamma| \in \{28, 32, 36\}$. □

See Definition 1.41 for the graphs $\mathcal{S}p_2(6)$ and $\mathcal{N}\mathcal{S}p^-(6)$.

Jonathan I. Hall and Ernest E. Shult actually prove a lot more in [HS85]. They characterize the graphs that are locally cotriangular in the following sense. A graph is said to be *cotriangular* if every pair $x, y$ of nonadjacent vertices is contained in a cotriangle $\{x, y, z\}$ (that is a 3-cocique) such that every other vertex is adjacent to either all or exactly one of the vertices $x, y, z$. Observe that a join $\Gamma + \Lambda$ is cotriangular if and only if both $\Gamma$ and $\Lambda$ are. As in Example 1.50 denote with $\Gamma^*$ the reduced graph of $\Gamma$, that is the graph $\Gamma/\Pi$ where $\Pi$ partitions the vertices of $\Gamma$ into sets of vertices with the same closed neighborhood. Then $\Gamma$ is cotriangular if and only if $\Gamma^*$ is. A graph $\Gamma$ is called completely reduced in this context whenever $\Gamma^* = \Gamma$ and $\Gamma$ can't be decomposed into $\Gamma_1 + \Gamma_2$ with nonempty $\Gamma_1, \Gamma_2$. A classification of all cotriangular graphs is given by the following theorem.

**Theorem 2.11.** ([HS85, Cotriangle Theorem]) *A finite completely reduced graph is cotriangular if and only if it is isomorphic to one of the graphs*

$$K(n, 2), \ n \geqslant 2; \quad \mathcal{S}p_2(2n), \ n \geqslant 2; \quad \mathcal{N}\mathcal{S}p^\varepsilon(2n), \ \varepsilon = \pm 1, \ n \geqslant 3.$$ □

The graphs $\mathcal{S}p_2(2n)$ and $\mathcal{N}\mathcal{S}p^\varepsilon(2n)$ have been introduced in Definition 1.41, and describe the orthogonality relation with respect to a nondegenerate symplectic form on an even-dimensional vector space over $\mathbb{F}_2$. The graphs $K(2,2) \cong K_1$ and $K(3,2) \cong \overline{K_3}$ are considered degenerate. Let $\mathcal{D}$ denote the set of graphs $\Gamma$ such that $\Gamma^*$ is a finite completely reduced cotriangular graph.

**Theorem 2.12.** ([HS85, Main Theorem]) *Let $\Gamma$ be connected and locally $\mathcal{D}$. Then either $\Gamma$ is locally $\{K_1, \overline{K_3}\}$ or $\Gamma$ is isomorphic to one of the following graphs*

- $K(n, 2)$ *where* $n \geqslant 7$,
- $\mathcal{S}p_2(2n)$ *possibly with a polar subspace deleted*,
- $\mathcal{H}_{2n}^\varepsilon(T)$, $\mathcal{G}_{2n}^\varepsilon$,
- $3 \cdot K(7, 2)$, $\Sigma L_{2,25}$, *or* $\mathcal{N}_6^+(3)$. □

The graphs $\mathcal{H}_{2n}^\varepsilon(T)$, $\mathcal{G}_{2n}^\varepsilon$ are derived from the graph $\mathcal{S}p_2(2n)$. For a precise definition and description of these graphs as well as $\mathcal{N}_6^+(3)$ we refer to [HS85]. As Hall and Shult remark, the polar subspaces of $\mathcal{S}p_2(2n)$ as well as the graphs $\mathcal{H}_{2n}^\varepsilon(T)$ can be cataloged so that Theorem 2.12 indeed classifies all graphs that are locally $*$-finite nondegenerate indecomposable cotriangular graphs. It is further remarked that a similar classification can still be done when the indecomposability condition is dropped. Note that the case $k = 2$ of Theorem 2.8 as well as Theorem 2.9 and Theorem 2.10 are special cases of the classification in Theorem 2.12. Likewise, the following result which can also be found in [BH77] can be obtained from Theorem 2.12.

**Theorem 2.13.** ([HS85, Theorem 5]) *Let $\Gamma$ be connected and locally $\mathcal{S}p_2(2n)$ for some $n \geqslant 2$. Then $\Gamma$ is isomorphic to one of the following graphs*

- $\mathcal{NS}p^+(2n+2)$,
- $\mathcal{NS}p^-(2n+2)$,
- $\mathcal{S}p_2(2n+2)$ *minus* $\{x\} \cup x^\perp$ *for some $x$.*                                     $\square$

We finally state a local recognition result due to Ralf Gramlich. Following [Gra02] we define.

**Definition 2.14.** *The* space-complement graph $\mathrm{SC}_{n,k}(q)$ *is the graph with vertices the pairs $(x, X)$ where $x$ is a subspace of the projective space $\mathbb{P}_n(\mathbb{F}_q)$ of dimension $k$ and $X$ is a subspace of codimension $k$ such that $x$ and $X$ don't intersect. Further, two pairs $(x, X)$ and $(y, Y)$ are defined to be adjacent whenever $x \subseteq Y$ and $y \subseteq X$.*

**Remark 2.15.** The graphs $\mathrm{SC}_{n,k}(q)$ can be seen as $q$-analogs (see Remark 1.6) of the Kneser graphs. To illustrate this, we attempt to give a meaning to $\mathrm{SC}_{n,k}(1)$. As indicated in Remark 1.6, the projective space $\mathbb{P}_n(\mathbb{F}_q)$ is a $q$-analog of the $n+1$-element set $[n+1]$ with $k$-dimensional subspaces corresponding to $k+1$-element subsets of $[n+1]$. Adopting this language, let $x$ be a $k$-dimensional subspace of $[n+1]$ and $X$ a subspace of codimension $k$ that doesn't intersect with $x$. Accordingly, $X$ has to be the complement of $x$. Therefore the vertices of $\mathrm{SC}_{n,k}(1)$ are the $k$-dimensional subspaces of $[n+1]$. Two such subspaces $x$, $y$ are adjacent in $\mathrm{SC}_{n,k}(1)$ whenever $x$ is contained in the complement of $y$ and $y$ is contained in the complement of $x$. Since this is equivalent to $x$ and $y$ being disjoint, the graph $\mathrm{SC}_{n,k}(1)$ is just the Kneser graph $K(n+1, k+1)$. Theorem 2.8 then states that a connected graph $\Gamma$ is locally $\mathrm{SC}_{n,k}(1)$ if and only if $\Gamma \cong \mathrm{SC}_{n+k+1,k}(1)$ provided that $k \geqslant 0$, $n \geqslant 3(k+1)$.

**Theorem 2.16.** ([Gra02, 2.5.4]) *Let $k \geqslant 0$, $n \geqslant 4(k+1) - 1$. A connected graph $\Gamma$ is locally $\mathrm{SC}_{n,k}(q)$ if and only if $\Gamma \cong \mathrm{SC}_{n+k+1,k}(q)$.*                    $\square$

The graph $\mathrm{SC}_{n,0}(q)$ is referred to as the *point-hyperplane graph* and the graph $\mathrm{SC}_{n,1}(q)$ accordingly as the *line-hyperline graph* of the projective space $\mathbb{P}_n(\mathbb{F}_q)$. A proof for the local recognition of line-hyperline graphs can be found in [Gra04] where $n = 6$ is permitted provided that $q \neq 2$. See [CCG05] for the local recognition of point-hyperplane graphs.

## 2.2 Local recognition of simply laced Weyl graphs

Throughout this section, let $M$ be a connected crystallographic Dynkin diagram. Consequently, $M$ is one of the diagrams classified in Theorem 1.82. We defined $W(M)$ to be the Weyl group generated by a crystallographic root system of type $M$ along with the notion of a reflection of $W(M)$ being a short or a long root reflection. Further, recall from Definition 1.37 that the commuting graph of a group $G$ on elements $X \subset G$ is the graph with vertices $X$ and two elements adjacent whenever they commute.

**Definition 2.17.** *Let $M$ be a connected crystallographic Dynkin diagram. The* Weyl graph $\mathbb{W}(M)$ *is the commuting graph of $W(M)$ on its reflections. If $M$ is not simply laced then we regard $\mathbb{W}(M)$ as a bichromatic graph with short (respectively long) vertices given by the short (respectively long) root reflections of $W(M)$.*

If $M$ is simply laced then all reflections in $W(M)$ are conjugate. The Weyl graph $\mathbb{W}(M)$ is therefore locally homogeneous. On the other hand, if $M$ is not simply laced then there are two conjugacy classes of reflections in $W(M)$, namely short and long root reflections. Accordingly, the bichromatic graph $\mathbb{W}(M)$ is locally homogeneous as well. Note that the (long) local graph of $\mathbb{W}(M)$ can be obtained using Remark 1.83.

**Remark 2.18.** We have the following generalization of Weyl graphs to arbitrary Coxeter systems. Let $(W,S)$ be a Coxeter system of rank $n$. In Theorem 1.70 we constructed an embedding of $W$ into $\mathrm{GL}(n,\mathbb{R})$ such that the Coxeter generators $s \in S$ are generalized reflections. For this reason one refers to all conjugates of some $s \in S$ as reflections of $(W,S)$. The *reflection graph* of $(W,S)$ is the commuting graph of $W$ on the reflections of $(W,S)$.

Let $\Phi$ be a root system corresponding to $W(M)$. Recall that the reflections of $W(M)$ are given as reflections through the roots of $\Phi$. Notice that the reflections through roots $\alpha$ and $\beta$ commute if and only if the roots $\alpha$ and $\beta$ are orthogonal. Accordingly, the Weyl graph $\mathbb{W}(M)$ can be constructed as the graph on the roots of $\Phi$, with roots $\alpha$ and $-\alpha$ identified for each $\alpha \in \Phi$, such that two roots are adjacent whenever they are orthogonal. This observation allows us to construct the graphs $\mathbb{W}(M)$ from the previous descriptions of the root systems $\Phi(M)$.

**Example 2.19.** Recall from Example 1.84 that $\Phi(A_n) \subset \mathbb{R}^{n+1}$ consists of the roots $\varepsilon_i - \varepsilon_j$ for $i \neq j \in [n+1]$. For $i < j$ denote with $y_{i,j}$ the reflection through $\pm(\varepsilon_i - \varepsilon_j)$. $\mathbb{W}(A_n)$ is the graph with vertices $y_{i,j}$, $i < j \in [n+1]$ and adjacency given by

$$y_{i,j} \perp y_{k,l} \iff \{i,j\} \cap \{k,l\} = \emptyset.$$

Consequently, the Weyl graph $\mathbb{W}(A_n)$ is isomorphic to the Kneser graph $K(n+1,2)$. In particular, $\mathbb{W}(A_n)$ is connected if and only if $n \geqslant 4$, in which case Theorem 1.59 implies that $\mathrm{Aut}(\mathbb{W}(A_n)) \cong W(A_n)$.

**Example 2.20.** In Example 1.87 we constructed $\Phi(D_n) \subset \mathbb{R}^n$ as containing the roots $\pm \varepsilon_i \pm \varepsilon_j$ for $i \neq j \in [n]$. For $i < j$ denote with $y_{i,j}$ the reflection through $\pm(\varepsilon_i - \varepsilon_j)$, and with $y_{j,i}$ the reflection through $\pm(\varepsilon_i + \varepsilon_j)$. Accordingly, $\mathbb{W}(D_n)$ is the graph with vertices $y_{i,j}$, $i \neq j \in [n]$ and adjacency given by

$$y_{i,j} \perp y_{k,l} \iff \{i,j\} \cap \{k,l\} = \emptyset \ \lor \ (k,l) = (j,i).$$

We observe that $\mathbb{W}(D_n)$ is isomorphic to the composition graph $K(n,2)[K_2]$. Accordingly, $\mathbb{W}(D_n)^* \cong K(n,2)$ when $n \geqslant 5$. In particular, $\mathbb{W}(D_n)$ is connected if and only if $n \geqslant 5$. For $n \geqslant 5$ we also infer from Theorem 1.59 that

$$\mathrm{Aut}(\mathbb{W}(D_n)) \cong (\mathbb{Z}/2\mathbb{Z})^{\binom{n}{2}} \rtimes \mathrm{Sym}_n$$

where $\text{Sym}_n$ acts in the natural way on the set of all 2-subsets of $[n]$ indexing the copies of $\mathbb{Z}/2\mathbb{Z}$.

**Example 2.21.** Let $n \geqslant 2$. Following the construction in Example 1.85 the root system $\Phi(B_n)$ consists of the roots $\pm \varepsilon_i$ and $\pm \varepsilon_i \pm \varepsilon_j$ for $i \neq j \in [n]$. Denote with $y_{i,i}$ the reflection through $\pm \varepsilon_i$. As in the previous example, for $i < j$ denote with $y_{i,j}$ the reflection through $\pm(\varepsilon_i - \varepsilon_j)$ and with $y_{j,i}$ the reflection through $\pm(\varepsilon_i + \varepsilon_j)$. Accordingly, we see that $\mathbb{W}(B_n)$ is the bichromatic graph with vertices $y_{i,j}$, $i, j \in [n]$, where the $y_{i,i}$ are short and the $y_{i,j}$ with $i \neq j$ are long vertices, and adjacency is given by

$$y_{i,j} \perp y_{k,l} \iff \{i,j\} \cap \{k,l\} = \emptyset \ \vee \ (k,l) = (j,i).$$

In particular, $\mathbb{W}(B_n)$ is connected if and only if $n \geqslant 3$. Let $n \geqslant 3$. For any $i \in [n]$ the vertex $y_{i,i}$ is characterized as the unique short neighbor of the long vertices $y_{k,l}$ for which $i \notin \{k, l\}$. Henceforth, an automorphism of $\mathbb{W}(B_n)$ is determined by how it acts on the long induced subgraph of $\mathbb{W}(B_n)$ which is $\mathbb{W}(D_n)$. Consequently, $\text{Aut}(\mathbb{W}(B_n))$ is isomorphic to $\text{Aut}(\mathbb{W}(D_n))$.

**Proposition 2.22.** *We have the following isomorphisms.*

- $\mathbb{W}(E_6) \cong \mathcal{N}\mathcal{S}p^-(6)$,
- $\mathbb{W}(E_7) \cong \mathcal{S}p_2(6)$,
- $\mathbb{W}(E_8) \cong \mathcal{N}\mathcal{S}p^+(8)$.

**Proof.** Verified in Proposition A.1. □

We now provide local recognition results for the Weyl graphs $\mathbb{W}(M)$ for the case that $M$ is a simply laced Dynkin diagram. These results are basically restatements of local recognition results for locally cotriangular graphs as studied in [HS85] and presented in the previous section.

**Theorem 2.23.** *Let $n \geqslant 6$. A connected graph $\Gamma$ is locally $\mathbb{W}(A_n)$ if and only if $\Gamma \cong \mathbb{W}(A_{n+2})$.*

**Proof.** By Example 2.19 $\mathbb{W}(A_n)$ is isomorphic to the Kneser graph $K(n+1, 2)$. Accordingly, Theorem 2.8 applies. □

**Theorem 2.24.** *Let $n \geqslant 7$. A connected graph $\Gamma$ is locally $\mathbb{W}(A_1) + \mathbb{W}(D_n)$ if and only if $\Gamma \cong \mathbb{W}(D_{n+2})$.*

**Proof.** Exploiting the assumed local structure we see that the vertices of $\Gamma$ come in pairs. Namely, for each vertex $x$ there exists a unique vertex $x'$ such that the closed neighborhoods of $x$ and $x'$ agree. Therefore $\Gamma \cong \Gamma^*[K_2]$. With Example 2.20 in mind, $\Gamma^*$ is seen to be connected and locally $\mathbb{W}(A_{n-1})$. According to Theorem 2.23 this shows that $\Gamma^* \cong \mathbb{W}(A_{n+1})$. The claim now follows from Example 2.20. □

**Theorem 2.25.** *Let $\Gamma$ be a connected graph.*

- *If $\Gamma$ is locally $\mathbb{W}(A_5)$ and $|\Gamma| = 36$ then $\Gamma \cong \mathbb{W}(E_6)$.*

- *If $\Gamma$ is locally $\mathbb{W}(D_6)$ then $\Gamma \cong \mathbb{W}(E_7)$.*
- *If $\Gamma$ is locally $\mathbb{W}(E_7)$ and $|\Gamma| = 120$ then $\Gamma \cong \mathbb{W}(E_8)$.*

**Proof.** According to Theorem 2.10 there are three connected graphs that are locally $\mathbb{W}(A_5)$. Among these, $\mathcal{NS}p^-(6)$ is the only graph with 36 vertices. The claim therefore follows from the isomorphism $\mathbb{W}(E_6) \cong \mathcal{NS}p^-(6)$ stated in Proposition 2.22.

For the second claim, recall from Example 2.20 that $\mathbb{W}(D_6)^* \cong K(6, 2)$. Hence Theorem 2.12 applies, and we can use $\mathbb{W}(E_7) \cong \mathcal{S}p_2(6)$ as provided by Proposition 2.22.

By Proposition 2.22 we also have $\mathbb{W}(E_8) \cong \mathcal{NS}p^+(8)$ whence the third claim follows from Theorem 2.13.                                                                                      □

**Remark 2.26.** Recall that while $\mathbb{W}(E_6)$ and $\mathbb{W}(E_8)$ are not locally recognizable in the strictest sense there are only two further connected graphs each that are locally like $\mathbb{W}(E_6)$ and $\mathbb{W}(E_8)$ respectively.

**Remark 2.27.** In Remark 2.18 we introduced the notion of the reflection graph of a Coxeter system $(W, S)$ as the graph on the conjugates of $S$ such that two elements are adjacent whenever they commute. In the case that $(W, S)$ is irreducible and crystallographic the reflection graph is just the Weyl graph (with the types of vertices forgotten). According to the classification in Theorem 1.74 the only finite irreducible Coxeter systems that are not crystallographic are those of types $H_3$, $H_4$ and $I_2(m)$. Let $M$ be one of the corresponding Coxeter diagrams, let $(W, S)$ be of type $M$, and denote with $\mathbb{W}(M)$ the reflection graph on $(W, S)$.

In the case that $M$ is $I_2(m)$, one easily checks, using for instance the associated root system described in Example 1.81 and observing that orthogonal roots only appear in pairs and only when $m$ is even, that

$$\mathbb{W}(I_2(m)) \cong \begin{cases} m \cdot K_1 & \text{if } 2 \nmid m, \\ m/2 \cdot K_2 & \text{if } 2 \mid m. \end{cases}$$

Now, let $M$ be one of the Coxeter diagrams $H_3$ and $H_4$. In Proposition A.2 it is verified that $\mathbb{W}(H_3)$ is isomorphic to $5 \cdot K_3$, and that the reflection graph $\mathbb{W}(H_4)$ is a connected graph on 60 vertices that is locally $\mathbb{W}(H_3)$. Notice that by Proposition 1.46 and Proposition 2.6 the Cartesian product of five copies of $K_4$ is another example of a connected graph that is locally $\mathbb{W}(H_3)$

## 2.3 Local recognition of non-simply laced Weyl graphs

Let $M$ be a connected crystallographic Dynkin diagram, as classified in Theorem 1.82, which is not simply laced. Recall that this means that the root system $\Phi(M)$ contains roots of two lengths. Consequently, the Weyl graph $\mathbb{W}(M)$ is bichromatic. Recall that for a graph $\Gamma$ we denote with $\Gamma^s$ (respectively $\Gamma^\ell$) the bichromatic graph obtained from $\Gamma$ by considering all vertices as short (respectively long).

**Example 2.28.** Notice that each root in the root system $\Phi(G_2)$ defined in Example 1.90 is orthogonal to exactly one other root, and that two such orthogonal roots are of different type. The Weyl graph $\mathbb{W}(G_2)$ therefore is isomorphic to $3 \cdot (K_1^s + K_1^\ell)$. In other words, $\mathbb{W}(G_2)$ consists of three disjoint edges of mixed type.

In the sequel we are interested in local recognition results for the Weyl graphs $\mathbb{W}(B_n)$, $\mathbb{W}(C_n)$ and $\mathbb{W}(F_4)$. Extending our previous terminology we say that a bichromatic graph is *locally homogeneous* if the local graphs at short vertices are all isomorphic to some bichromatic graph $\Delta_s$ and if the local graphs at long vertices are all isomorphic to some bichromatic graph $\Delta_\ell$. In this case we say that $\Delta_s$ is the *short local graph* of $\Gamma$ and that $\Delta_\ell$ is the *long local graph* of $\Gamma$. If $\Gamma$ is a bichromatic locally homogeneous graph then we denote its short local graph by $\Delta_s(\Gamma)$ and its long local graph by $\Delta_\ell(\Gamma)$. If $\Lambda$ is another bichromatic locally homogeneous graph such that $\Delta_s(\Lambda) \cong \Delta_s(\Gamma)$ as well as $\Delta_\ell(\Lambda) \cong \Delta_\ell(\Gamma)$ then we say that $\Lambda$ is *locally like* $\Gamma$.

**Theorem 2.29. (joint with Gramlich, Hall)** *Let $n \geqslant 4$, and let $\Gamma$ be a connected bichromatic locally homogeneous graph with $\Delta_s(\Gamma) \cong \mathbb{W}(B_{n+1})$ and $\Delta_\ell(\Gamma) \cong \mathbb{W}(A_1)^\ell + \mathbb{W}(B_n)$. Then $\Gamma \cong \mathbb{W}(B_{n+2})$.*

**Proof.** Let $X$ be a short component of $\Gamma$ and $x \in X$ a short vertex. The short induced subgraph of $x^\perp$ is a clique on $n+1$ elements which by Lemma 2.3 implies that $X$ is a clique on $n+2$ elements. By assumption, the long neighbors of $x$ induce a subgraph isomorphic to the long induced subgraph of $\mathbb{W}(B_{n+1})$. This subgraph is isomorphic to $\mathbb{W}(D_{n+1})$. In particular, it is connected for $n \geqslant 4$, see Example 2.20. This implies that all long neighbors of $x$ are contained in a single long component $Y$ of $\Gamma$. Consider a short vertex $x_1 \in X$ adjacent to $x$. Again, all long neighbors of $x_1$ lie in one long component of $\Gamma$. But looking at $\{x, x_1\}^\perp \subset x^\perp$ we see that $x$ and $x_1$ share long neighbors whence this component has to be $Y$ as well. Since $X$ is connected this shows that all long vertices adjacent to some vertex of $X$ are contained in $Y$. Likewise, let $y \in Y$. The short induced subgraph of $y^\perp$ is a clique on $n$ vertices and thus in particular connected. Again, we see that for a long vertex $y_1$ adjacent to $y$ the common neighbors $\{y, y_1\}^\perp$ contain a short vertex. Therefore the same argument as before shows that all short vertices adjacent to some vertex of $Y$ are contained in $X$. Since $\Gamma$ is connected this proves that $X$ and $Y$ are the only short respectively long components of $\Gamma$.

We count the number of long vertices by counting the long neighbors of the $n+2$ short vertices of $\Gamma$. By assumption, a short vertex has $(n+1)n$ long neighbors. Further, two short vertices have $n(n-1)$ long neighbors in common, three short vertices have $(n-1)(n-2)$ long neighbors in common, and so on. Thus there are

$$\binom{n+2}{1}(n+1)n - \binom{n+2}{2}n(n-1) + \ldots + (-1)^{n+1}\binom{n+2}{n}2 = (n+2)(n+1)$$

long vertices in $\Gamma$. Note that for the above equation we exploited that the alternating sum of the binomial coefficients equals zero, that is

$$\sum_{k=0}^{n}(-1)^k \binom{n}{k} = 0.$$

Let $x_1, x_2, \ldots, x_{n+2}$ be the short vertices of $\Gamma$. $\Gamma$ is locally $\mathbb{W}(B_{n+1})$ at short vertices which implies that for $i \neq j \in [n+2]$ the common neighborhood

$$\{x_r \colon r \notin \{i,j\}\}^\perp$$

contains exactly two long vertices which we denote by $y_{i,j}$ and $y_{j,i}$. Since a long vertex is adjacent to exactly $n$ short vertices the $y_{i,j}$ thus defined are all distinct. By construction, $y_{i,j} \perp y_{j,i}$. Further, the $y_{i,j}$ exhaust $Y$ because $\Gamma$ contains exactly $(n+2)(n+1)$ long vertices. Given two vertices $y_{i,j}$ and $y_{k,l}$ we find $m \in [n+2] \setminus \{i,j,k,l\}$ whence

$$y_{i,j}, y_{k,l} \in x_m^\perp \cong \mathbb{W}(B_{n+1}).$$

$y_{i,j}$ is characterized in $x_m^\perp$ as one of the two long vertices contained in $\{x_r \colon r \notin \{i,j,m\}\}^\perp$. Likewise, $y_{k,l}$ is characterized in $x_m^\perp$ as one of the two long vertices contained in $\{x_r \colon r \notin \{k,l,m\}\}^\perp$. Consequently,

$$y_{i,j} \perp y_{k,l} \quad \Longleftrightarrow \quad \{i,j\} \cap \{k,l\} = \emptyset$$

for $\{i,j\} \neq \{k,l\}$. According to Example 2.21, $\Gamma \cong \mathbb{W}(B_{n+2})$. □

**Corollary 2.30.** *Let $n \geqslant 4$, and let $\Gamma$ be a connected bichromatic locally homogeneous graph with $\Delta_s(\Gamma) \cong \mathbb{W}(A_1)^s + \mathbb{W}(C_n)$ and $\Delta_\ell(\Gamma) \cong \mathbb{W}(C_{n+1})$. Then $\Gamma \cong \mathbb{W}(C_{n+2})$.* □

## 2.4   Graphs locally like $\mathbb{W}(F_4)$

Consider the Weyl graph $\mathbb{W}(F_4)$. Exploiting the description of the crystallographic root system $\Phi(F_4)$ given in Example 1.89 we find that $\mathbb{W}(F_4)$ is a connected bichromatic locally homogeneous graph on 24 vertices with short local graph $\mathbb{W}(B_3)$ and long local graph $\mathbb{W}(C_3)$. Recall that we agreed to depict long vertices as filled dots and short vertices as unfilled dots. We use the description given in Example 2.21 to draw the local graphs of $\mathbb{W}(F_4)$ as

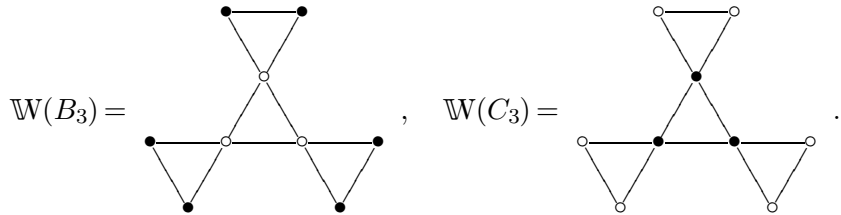

We will shortly see that $\mathbb{W}(F_4)$ is not locally recognizable. Before we turn to investigating additional constraints under which we seek to recognize $\mathbb{W}(F_4)$ nonetheless, we study connected bichromatic graphs $\Gamma$ which are locally like $\mathbb{W}(F_4)$. That is, we study locally homogeneous bichromatic graphs such that $\Delta_s(\Gamma) \cong \mathbb{W}(B_3)$ and $\Delta_\ell(\Gamma) \cong \mathbb{W}(C_3)$. The results we obtain then guide our way in determining appropriate conditions under which we will be (almost) able to recognize $\mathbb{W}(F_4)$. An easy but crucial observation to start with is the following.

**Proposition 2.31.** *Let $\Gamma$ be locally like $\mathbb{W}(F_4)$. The short (respectively long) induced subgraph of $\Gamma$ is isomorphic to a disjoint union of 4-cliques.* □

Let $\Gamma$ be a bichromatic graph that is locally like $\mathbb{W}(F_4)$. Observe that the graph obtained from $\Gamma$ by exchanging the roles of short and long vertices is locally like $\mathbb{W}(F_4)$ as well. Results that we obtain for short vertices of graphs locally like $\mathbb{W}(F_4)$ are therefore also true for long vertices.

Paraphrasing Proposition 2.31, the vertices of $\Gamma$ come in 4-cliques of the same type. In order to simplify things it is natural to investigate the collapsed graph $\Gamma/\Pi$ where $\Pi$ is the partition of $\Gamma$ into short and long 4-cliques. To this end, we analyze how these 4-cliques relate to each other.

**Proposition 2.32.** *Let $\Gamma$ be locally like $\mathbb{W}(F_4)$, and $x_1, x_2, x_3, x_4$ a short 4-clique in $\Gamma$. Let $i \neq j$ and $k \neq l$.*

- *$\{x_i, x_j\}^\perp$ is locally $K_2^s \sqcup K_2^\ell$. In particular, for any pair $x_i, x_j$ there exist unique long vertices $y_{i,j}, y_{j,i}$ contained in $\{x_i, x_j\}^\perp$.*
- *$\{x_i, x_j, x_k\}^\perp$ contains no long vertex if $i, j, k$ are distinct. In particular, the vertices $y_{i,j}$ are all distinct.*
- *There are exactly 12 long vertices adjacent to at least one of the $x_i$, namely the above vertices $y_{i,j}$.*
- *$y_{i,j} \perp y_{k,l}$ implies that $\{k,l\} = \{i,j\}$ or $\{k,l\} \cap \{i,j\} = \emptyset$.*

**Proof.** Exploiting the local structure at $x_i$ we see that every short adjacent pair $x_i$, $x_j$ has exactly two long neighbors in common which we will (arbitrarily) denote by $y_{i,j}$ and $y_{j,i}$. Accordingly, $y_{i,j} \perp y_{j,i}$. Looking at the neighbors of a vertex $y_{i,j}$ reveals that $x_i$ and $x_j$ are the only short vertices among $x_1, x_2, x_3, x_4$ which are adjacent to $y_{i,j}$. Consequently, the $y_{i,j}$ are 12 distinct vertices. Since three adjacent short vertices share no long neighbors we count that exactly

$$\binom{4}{1} 6 - \binom{4}{2} 2 = 12$$

long vertices are neighbored to at least one of the vertices $x_1, x_2, x_3, x_4$. Consequently, the long neighbors of the $x_i$ are precisely the vertices $y_{i,j}$.

For the last claim, assume that $y_{i,j} \perp y_{k,l}$ and $\{k, l\} \cap \{i, j\} = \{i_0\}$. A look at the neighbors of $x_{i_0}$ shows that this is a contradiction. □

**Proposition 2.33.** *Let $\Gamma$ be locally like $\mathbb{W}(F_4)$, and let $\Pi$ be the partition of $\Gamma$ into short and long 4-cliques.*

- *The contraction $\Gamma/\Pi$ has short local graphs each isomorphic to one of $\{\overline{K_n}^\ell: n \in \{3, 4, 5, 6\}\}$ and long local graphs each isomorphic to one of $\{\overline{K_n}^s: n \in \{3, 4, 5, 6\}\}$.*
- *Every graph homomorphism $\psi \colon \Gamma \to \Gamma$ induces a graph homomorphism $\Gamma/\Pi \to \Gamma/\Pi$ given by $X \mapsto \psi(X)$. In particular, if $\Gamma$ is transitive then $\Gamma/\Pi$ is transitive as well.*

**Proof.** Let $X \in \Gamma/\Pi$ be a short vertex. By Proposition 2.31 the local graph at $X$ is isomorphic to $\overline{K_n}^\ell$ for some $n$. Since the long induced subgraph of the local graph at a short vertex in $\Gamma$ is isomorphic to $3 \cdot K_2^\ell$ we see that $n \geqslant 3$. On the other hand, $X = \{x_1, x_2, x_3, x_4\}$ is a 4-clique of $\Gamma$ and according to Proposition 2.32 there are 12 long vertices $y_{i,j}$ at distance 1 from $X$ in $\Gamma$. Since $y_{i,j}$ and $y_{j,i}$ are adjacent in $\Gamma$ they are identified in $\Gamma/\Pi$ which shows $n \leqslant 6$.

Let $\psi \colon \Gamma \to \Gamma$ be a graph homomorphism. With Proposition 2.31 in mind, it is clear that $\psi(X) \in \Pi$ for every $X \in \Pi$. By Proposition 1.49 this ensures that $\Gamma/\Pi \to \Gamma/\Pi$ given by $X \to \psi(X)$ is well-defined and a partial graph homomorphism (see Definition 1.31). Let $X \perp Y$ in $\Gamma/\Pi$. Consequently, either $\psi(X) = \psi(Y)$ or $\psi(X) \perp \psi(Y)$. But the local structure of $\Gamma/\Pi$ just proved implies that one of $X$ and $Y$ is short while the other is long. Therefore, $\psi(X)$ and $\psi(Y)$ are of different type as well. In particular, $\psi(X) \neq \psi(Y)$ which shows that $\psi$ is in fact a graph homomorphism. $\square$

We now do the reverse. Starting with a what-might-be-possible collapse of a graph $\Gamma$ that is locally like $\mathbb{W}(F_4)$ we try to reconstruct $\Gamma$. Since different $\Gamma$ can have isomorphic collapses this reconstruction necessarily is not unique.

**Lemma 2.34.** *For every connected bipartite graph $\Lambda$ that is locally $\overline{K_6}$ there is a connected bichromatic graph $\Gamma$ that is locally like $\mathbb{W}(F_4)$ such that $\Gamma/\Pi = \Lambda$ where $\Pi$ is the partition of $\Gamma$ into short and long 4-cliques.*

**Proof.** Let $\Lambda$ be a bipartite graph that is locally $\overline{K_6}$. Exploiting that $\Lambda$ is 2-colorable, we may identify $\Lambda$ with a bichromatic graph such that no two vertices of the same type are adjacent. Accordingly, $\Lambda$ is locally homogeneous with $\Delta_s(\Lambda) \cong \overline{K_6}^\ell$ and $\Delta_\ell(\Lambda) \cong \overline{K_6}^s$. For any vertex $x$ of $\Lambda$ choose a bijection

$$x^\perp \to \binom{[4]}{2}, \quad y \mapsto a(x, y)$$

between its six neighbors and the six 2-subsets of $[4]$. To every directed edge $(x, y)$ we thus assigned the 2-subset $a(x, y)$ of $[4]$. Construct the bichromatic graph $\Gamma$ from $\Lambda$ as follows. For every vertex $x \in \Lambda$ add a 4-clique $x_1, x_2, x_3, x_4$ of the same type as $x$ to $\Gamma$. Moreover, for $x, y \in \Lambda$ let $x_i$ and $y_j$ be adjacent in $\Gamma$ if and only if $x$ and $y$ are adjacent in $\Lambda$ and $(i, j) \in a(x, y) \times a(y, x)$. By construction, contracting the 4-cliques of $\Gamma$ produces $\Lambda$.

We claim that $\Gamma$ is locally like $\mathbb{W}(F_4)$. Let $x_i$ be a short vertex of $\Gamma$. By construction, $x_i$ is adjacent to exactly three other short vertices in $\Gamma$, namely those that it constitutes the short 4-clique $x_1, x_2, x_3, x_4$ with. Let $y_j$ be a long vertex adjacent to $x_i$ in $\Gamma$. Then $x$ and $y$ are adjacent in $\Lambda$ and $i \in a(x, y)$, $j \in a(y, x)$. In particular, we see that $x_i$ has six long neighbors in $\Gamma$ and that the long induced subgraph on $x_i^\perp$ is isomorphic to $3 \cdot K_2^\ell$. Let $k$ be the index such that $\{i, k\} = a(x, y)$. Observe that $x_k$ is the unique short vertex contained in $\{x_i, y_j\}^\perp$. Now, consider any short neighbor $x_k$ of $x_i$, and let $y$ be the unique neighbor of $x$ in $\Lambda$ such that $a(x, y) = \{i, k\}$. By construction, the common long neighbors of $x_i$ and $x_k$ are exactly the two adjacent vertices $y_j$, $j \in a(y, x)$. This shows that the short local graph of $\Gamma$ is isomorphic to $\mathbb{W}(B_3)$ and the same reasoning reasoning reveals that the long local graph of $\Gamma$ is isomorphic to $\mathbb{W}(C_3)$. $\square$

The construction in Lemma 2.34 is easily adapted to work for locally $\overline{K_3}$ graphs $\Lambda$ as well. In this case, we basically need to create twice as many edges between the vertices of the graph $\Gamma$ constructed from $\Lambda$.

**Lemma 2.35.** *For every connected bipartite graph $\Lambda$ that is locally $\overline{K_3}$ there is a connected bichromatic graph $\Gamma$ that is locally like $\mathbb{W}(F_4)$ such that $\Gamma/\Pi = \Lambda$ where $\Pi$ is the partition of $\Gamma$ into short and long 4-cliques.*

**Proof.** Let $\Lambda$ be a bipartite graph that is locally $\overline{K_3}$. Again, we may identify $\Lambda$ with a bichromatic graph such that $\Lambda$ is locally homogeneous with $\Delta_s(\Lambda) \cong \overline{K_3}^\ell$ and $\Delta_\ell(\Lambda) \cong \overline{K_3}^s$. For any vertex $x \in \Lambda$ choose a bijection

$$x^\perp \to \{\{1,2\},\{1,3\},\{1,4\}\}, \quad y \mapsto a(x,y)$$

between its three neighbors and the sets $\{1, 2\}$, $\{1, 3\}$, $\{1, 4\}$ (or any other three 2-subsets of $[4]$ which together with their complements exhaust the 2-subsets of $[4]$). This way we assign to every directed edge $(x, y)$ the 2-element set $a(x, y)$. Construct the graph $\Gamma$ from $\Lambda$ by adding for each vertex $x \in \Lambda$ a 4-clique $x_1, x_2, x_3, x_4$ of the same type as $x$ to $\Gamma$. Further, let $x_i \perp y_j$ in $\Gamma$ if and only if $x \perp y$ in $\Lambda$ and $(i,j) \in a(x,y) \times a(y,x)$ or $(i,j) \in \overline{a(x,y)} \times \overline{a(y,x)}$. By construction, $\Lambda$ is obtained when contracting the short and long 4-cliques of $\Gamma$.

As in the proof of Lemma 2.34 one checks that $\Gamma$ is indeed locally like $\mathbb{W}(F_4)$. □

**Corollary 2.36.** *There exist infinitely many finite connected bichromatic graphs that are locally like $\mathbb{W}(F_4)$.*

**Proof.** We claim that there are infinitely many finite connected bipartite graphs $\Lambda$ that are locally $\overline{K_6}$. For instance, we can use Lemma 1.46 and Lemma 2.6 to see that the graphs in the infinite family $C_k \square C_m \square C_n$ for $k, m, n \in \{4, 5, ..., \infty\}$ are connected and locally $\overline{K_6}$. Since cycles $C_n$ are 2-colorable whenever $n$ is even, Lemma 1.54 implies that $C_k \square C_m \square C_n$ is 2-colorable and hence bipartite whenever $k$, $m$, $n$ are all even.

Let $\Lambda$ be a connected bipartite graph that is locally $\overline{K_6}$. By Lemma 2.34, $\Lambda$ gives rise to at least one connected bichromatic graph $\Gamma_\Lambda$ that is locally like $\mathbb{W}(F_4)$ such that $\Lambda$ is reconstructable from $\Gamma_\Lambda$ by collapsing the short and long 4-cliques of $\Gamma_\Lambda$. Consequently, the graphs $\Gamma_\Lambda$ are nonisomorphic for nonisomorphic $\Lambda$. □

## 2.5 Recognition results for graphs locally like $\mathbb{W}(F_4)$

### 2.5.1 Properties of graphs locally like $\mathbb{W}(F_4)$

In the previous section we showed that the Weyl graph $\mathbb{W}(F_4)$ is not locally recognizable. In fact, in the course of the proof of Corollary 2.36 we explicitly constructed infinitely many connected graphs that are locally like $\mathbb{W}(F_4)$. In this section we study further properties of $\mathbb{W}(F_4)$ in order to characterize $\mathbb{W}(F_4)$ among the connected bichromatic graphs that are locally like $\mathbb{W}(F_4)$. To this end we start with some easy observations.

**Proposition 2.37.** *Let $\Gamma$ be a finite bichromatic graph that is locally like $\mathbb{W}(F_4)$. Then the numbers of short and long vertices in $\Gamma$ are the same.*

**Proof.** Let $n_s$ (respectively $n_\ell$) be the number of short (respectively long) vertices of $\Gamma$. We count the number of edges connecting a short and a long vertex in two ways. Looking at the neighborhood of a short vertex which by assumption is isomorphic to $\mathbb{W}(B_3)$ shows that every short vertex has six long neighbors. We conclude that there are exactly $6n_s$ edges connecting a short and a long vertex. Likewise, by looking at neighbors of long vertices we count $6n_\ell$ edges connecting a long and a short vertex. Thus $n_s = n_\ell$ as desired. $\square$

**Corollary 2.38.** *Let $\Gamma$ be a finite bichromatic graph that is locally like $\mathbb{W}(F_4)$. Then $|\Gamma|$ is divisible by $8$ and $|\Gamma| \geqslant 24$.*

**Proof.** By Proposition 2.31 the short (respectively long) vertices come in 4-cliques and by Proposition 2.37 there are as many short as long vertices. Hence $|\Gamma|$ is a multiple of 8. Looking at the neighborhood of a short vertex we see that there are at least three 4-cliques of long vertices in $\Gamma$. Consequently, there are at least 12 long vertices in $\Gamma$. It follows that $|\Gamma| \geqslant 24$. $\square$

Since $|\mathbb{W}(F_4)| = 24$ we see that, in a sense, $\mathbb{W}(F_4)$ is maximally tight among the graphs that are locally like $\mathbb{W}(F_4)$. There are several further properties besides the number of vertices that describe tightness of a graph, for instance its diameter. A third alternative to express tightness for a bichromatic graph is the notion of tight connectedness.

**Definition 2.39.** *A bichromatic graph is said to be* tightly connected *if every long vertex has a neighbor in every short component and vice versa.*

**Remark 2.40.** These three notions of tightness, however, are not local in nature where a local property is meant to be one which can be expressed in terms of the neighbors of each vertex. In this vague terminology, the property of a graph to be locally like $\mathbb{W}(F_4)$ is an instance of a local property because a graph has this property if and only if each short vertex has neighbors inducing a subgraph isomorphic to $\mathbb{W}(B_3)$ and each long vertex has neighbors inducing a subgraph isomorphic to $\mathbb{W}(C_3)$. On the other hand, for example the property of a graph to be tightly connected is not local because it requires each short vertex to have neighbors from every long component which involves the global knowledge about the entirety of all long vertices.

In order to find a more local notion to describe the tightness of $\mathbb{W}(F_4)$ we investigate the relation of vertices at distance 2.

**Proposition 2.41.** *Let $\Gamma$ be locally like $\mathbb{W}(F_4)$. Let $x, y \in \Gamma$ such that $d_\Gamma(x, y) = 2$.*

- *If $x, y$ are both short (respectively long) vertices then $\{x, y\}^\perp \cong \mu(x, y) \cdot K_1^\ell$ (respectively $K_1^s$) for some $\mu(x, y) \in \{1, 2, 3\}$.*

- If $x, y$ are of mixed type then $\{x, y\}^\perp \cong \mu_s(x, y) \cdot K_2^s \sqcup \mu_\ell(x, y) \cdot K_2^\ell$ for some $\mu_s(x, y), \mu_\ell(x, y) \in \{0, 1\}$.

**Proof.** First, let $x, y$ be two short vertices at distance 2. Since short vertices come in 4-cliques the common neighbors $\{x, y\}^\perp$ of $x, y$ can only contain long vertices. The structure of the local graph at such a long vertex $z$ implies that $\{x, y, z\}^\perp$ is empty. Accordingly, $\{x, y\}^\perp$ is a coclique on long vertices. Set $\mu(x, y) = |\{x, y\}^\perp|$. Analyzing the neighbors of $x$ we find that $\mu(x, y) \in \{1, 2, 3\}$.

Likewise, let $x$ be a short and $y$ a long vertex at distance 2. Let $z \in \{x, y\}^\perp$ be a, say short, vertex. By looking at the neighbors of $z$ we see that $\{x, y, z\}^\perp$ contains exactly one more short vertex. Hence, $\{x, y\}^\perp$ is a disjoint sum $\mu_s(x, y) \cdot K_2^s \sqcup \mu_\ell(x, y) \cdot K_2^\ell$ of pairs of short vertices and pairs of long vertices for some $\mu_s(x, y), \mu_\ell(x, y)$. Looking at $x^\perp$ and $y^\perp$ we further see that $\mu_s(x, y), \mu_\ell(x, y) \in \{0, 1\}$. □

We will be especially interested in the cases of graphs $\Gamma$ for which the parameters $\mu$, $\mu_s$, $\mu_\ell$ defined in Proposition 2.41 are constant. The following fact relates the values $\mu_s$ and $\mu_\ell$ to the local structure of the contractions $\Gamma/\Pi$ studied in the previous section, see Proposition 2.33.

**Proposition 2.42.** *Let $\Gamma$ be locally like $\mathbb{W}(F_4)$, and let $\Pi$ be the partition of $\Gamma$ into 4-cliques.*

- $\mu_s = \mu_\ell = 1$ *if and only if the contraction $\Gamma/\Pi$ is locally homogeneous with $\Delta_s(\Gamma/\Pi) \cong \overline{K_3}^\ell$ and $\Delta_\ell(\Gamma/\Pi) \cong \overline{K_3}^s$.*

- $\mu_s + \mu_\ell = 1$ *if and only if the contraction $\Gamma/\Pi$ is locally homogeneous with $\Delta_s(\Gamma/\Pi) \cong \overline{K_6}^\ell$ and $\Delta_\ell(\Gamma/\Pi) \cong \overline{K_6}^s$.*

**Proof.** Suppose that $\mu_s = \mu_\ell = 1$. Let $X = \{x_1, x_2, x_3, x_4\}$ be a short vertex of $\Gamma/\Pi$. Adopting the notation of Proposition 2.32 let $y_{i,j}$ and $y_{j,i}$ be the long vertices adjacent to both $x_i$ and $x_j$. Let $i \notin \{k, l\}$ so that $x_i$ and $y_{k,l}$ are at distance 2. Let $j$ be the index such that $\{i, j, k, l\} = [4]$. Since $\mu_\ell(x_i, y_{k,l}) = 1$ there are two long vertices adjacent to both $x_i$ and $y_{k,l}$. According to Proposition 2.32 the only possibilities are $y_{i,j}$ and $y_{j,i}$. Therefore $\{y_{i,j}, y_{j,i}, y_{k,l}, y_{l,k}\}$ is a long vertex of $\Gamma/\Pi$ and Proposition 2.32 implies that the local graph at $X$ is isomorphic to $\overline{K_3}^\ell$. Analogously for long local graphs.

On the other hand, assume that $\Gamma/\Pi$ is locally homogeneous with $\Delta_s(\Gamma/\Pi) \cong \overline{K_3}^\ell$ and $\Delta_\ell(\Gamma/\Pi) \cong \overline{K_3}^s$. Let $x$ be a short and $y$ be a long vertex of $\Gamma$ which are at distance 2. Observe that $x$ and $y$ are contained in adjacent equivalence classes of $\Gamma/\Pi$. Let $X = \{x_1, x_2, x_3, x_4\}$ be the short vertex of $\Gamma/\Pi$ where $x = x_i$ for some $i$. By assumption, $X$ is adjacent in $\Gamma/\Pi$ to exactly three long vertices which by Proposition 2.32 are of the form $\{y_{i,j}, y_{j,i}, y_{k,l}, y_{l,k}\}$ for $\{i, j, k, l\} = [4]$. $y$ therefore equals $y_{k,l}$ for some $\{k, l\} \not\ni i$. Choose $j$ as above. Then

$$\{x, y\}^\perp = \{x_i, y_{k,l}\}^\perp = \{x_k, x_l, y_{i,j}, y_{j,i}\},$$

and therefore $\mu_s(x, y) = \mu_\ell(x, y) = 1$.

The second equivalence is proved similarly. □

**Corollary 2.43.** *Let $\Gamma$ be locally like $\mathbb{W}(F_4)$, and $x_1, x_2, x_3, x_4$ a short 4-clique in $\Gamma$. As in Proposition 2.32 denote with $y_{i,j}$ and $y_{j,i}$ the long vertices adjacent to both $x_i$ and $x_j$. Let $\{i,j\} \cap \{k,l\} = \emptyset$.*

- *If $\mu_s = \mu_\ell = 1$ then $y_{i,j} \perp y_{k,l}$.*
- *If $\mu_s + \mu_\ell = 1$ then $y_{i,j} \not\perp y_{k,l}$.* □

**Remark 2.44.** For the Weyl graph $\mathbb{W}(F_4)$ the parameters $\mu$, $\mu_s$, $\mu_\ell$ are constant, and take the values $\mu = 3$ and $\mu_s = \mu_\ell = 1$ which is another instantiation of the tightness of $\mathbb{W}(F_4)$. Note that the property to have $\mu = 3$ and $\mu_s = \mu_\ell = 1$ is local in the sense of Remark 2.40 because it can be stated as follows: for each vertex $x$ and $y, z \in x^\perp$ it holds that $|\{y,z\}^\perp| = 3$ (respectively 4) if $y$ and $z$ are nonadjacent vertices of the same (respectively different) type.

### 2.5.2 Recognition results for $\mathbb{W}(F_4)$

The following theorem summarizes our recognition results for the Weyl graph $\mathbb{W}(F_4)$. Note that all of the provided conditions under which a graph $\Gamma$ is almost recognized as $\mathbb{W}(F_4)$ are statements which describe the tightness of $\Gamma$. We denote with $\Gamma_{24b}$ the graph on 24 vertices constructed in course of the proof of Proposition 2.46. An implementation of $\Gamma_{24b}$ in SAGE can be found in the appendix.

**Theorem 2.45.** *Let $\Gamma$ be a connected bichromatic graph that is locally like $\mathbb{W}(F_4)$. Assume that*

- *$|\Gamma| = 24$, or*
- *$\Gamma$ is tightly connected, or*
- *$\Gamma$ has diameter 2, or*
- *$\mu = 3$.*

*If one of these conditions holds then $\Gamma$ is isomorphic to $\mathbb{W}(F_4)$ or to $\Gamma_{24b}$. In particular,*

$$\mathrm{Aut}(\Gamma) \cong W(F_4)/Z.$$

We prove Theorem 2.45 by a series of propositions.

**Proposition 2.46. (joint with Gramlich, Hall)** *Let $\Gamma$ be a connected bichromatic graph that is locally like $\mathbb{W}(F_4)$. If $|\Gamma| = 24$ then $\Gamma \cong \mathbb{W}(F_4)$ or $\Gamma \cong \Gamma_{24b}$.*

**Proof.** As observed in Corollary 2.38, every graph that is locally like $\mathbb{W}(F_4)$ has at least 12 short and 12 long vertices. $\Gamma$ therefore consists of exactly 12 vertices of each type.

Let $x_1, x_2, x_3, x_4$ be a short 4-clique. Adopting the notation of Proposition 2.32, let $y_{i,j}$ and $y_{j,i}$ be the long vertices adjacent to both $x_i$ and $x_j$. The $y_{i,j}$ are 12 distinct vertices and therefore constitute the long vertices of $\Gamma$. It follows from Proposition 2.32 that the three long 4-cliques are given by $y_{i,j}, y_{j,i}, y_{k,l}, y_{l,k}$ for disjoint $\{i,j\}$ and $\{k,l\}$.

Each of the remaining eight short vertices has exactly two long neighbors in each of the three long 4-cliques. Let $x_5$ be one of remaining short vertices. The two neighbors of $x_5$ in a 4-clique $y_{i,j}, y_{j,i}, y_{k,l}, y_{l,k}$ are one of $y_{i,j}, y_{j,i}$ along with one of $y_{k,l}, y_{l,k}$. We ambiguously defined the vertices $y_{i,j}, y_{j,i}$ as the long vertices contained in $\{x_i, x_j\}^\perp$ so we may as well assume that $x_5$ is adjacent to $y_{i,j}$ and $y_{k,l}$ with $i < j$ and $k < l$. Let $x_6$ be the unique short vertex also adjacent to $y_{1,2}, y_{3,4}$. Likewise, let $x_7$ be the short vertex also adjacent to $y_{1,3}, y_{2,4}$, and $x_8$ the short vertex also adjacent to $y_{1,4}, y_{2,3}$. By construction, $x_5, x_6, x_7, x_8$ is a 4-clique. Notice that for instance $x_5, x_6 \in \{y_{1,2}, y_{3,4}\}^\perp$ implies that $x_7, x_8 \in \{y_{2,1}, y_{4,3}\}^\perp$. Altogether this determines the induced subgraph on $x_1, x_2, ..., x_8$ along with the vertices $y_{i,j}$.

Let $x_9, x_{10}, x_{11}, x_{12}$ be the remaining short 4-clique. We may assume that $x_9, x_{10}$ are the short vertices contained in $\{y_{1,2}, y_{4,3}\}^\perp$. Accordingly, $x_{11}, x_{12} \in \{y_{2,1}, y_{3,4}\}^\perp$. We may also assume that $x_9$ is contained in $\{y_{1,3}, y_{4,2}\}^\perp$ (because if both $x_9$ and $x_{10}$ were not contained in $\{y_{1,3}, y_{4,2}\}^\perp$ then both $x_{11}, x_{12} \in \{y_{1,3}, y_{4,2}\}^\perp$ which contradicts $x_{11}, x_{12} \in \{y_{2,1}, y_{3,4}\}^\perp$). Further, we may assume that $x_{11}$ is the second short vertex contained in $\{y_{1,3}, y_{4,2}\}^\perp$. Consider the two short vertices in $\{y_{1,4}, y_{3,2}\}^\perp$. These can be either $x_9, x_{12}$ or $x_{10}, x_{11}$, and either choice determines $\Gamma$. Denote with $\Gamma_{24a}$ the graph corresponding to the choice $x_9, x_{12} \in \{y_{1,4}, y_{3,2}\}^\perp$, and with $\Gamma_{24b}$ the graph corresponding to the choice $x_{10}, x_{11} \in \{y_{1,4}, y_{3,2}\}^\perp$. The following table summarizes adjacency involving the vertices $x_9, x_{10}, x_{11}, x_{12}$.

|  | by construction | $x_9, x_{10} \perp y_{1,2}, y_{4,3}$ | $x_{11}, x_{12} \perp y_{2,1}, y_{3,4}$ |
|---|---|---|---|
|  |  | $x_9, x_{11} \perp y_{1,3}, y_{4,2}$ | $x_{10}, x_{12} \perp y_{3,1}, y_{2,4}$ |
| $\Gamma_{24a}$ |  | $x_9, x_{12} \perp y_{1,4}, y_{3,2}$ | $x_{10}, x_{11} \perp y_{4,1}, y_{2,3}$ |
| $\Gamma_{24b}$ |  | $x_9, x_{12} \perp y_{4,1}, y_{2,3}$ | $x_{10}, x_{11} \perp y_{1,4}, y_{3,2}$ |

An implementation in SAGE of the graphs $\Gamma_{24a}$ and $\Gamma_{24b}$ can be found in the appendix. In particular, it is shown in Proposition A.3 using SAGE that $\Gamma_{24a}$ and $\Gamma_{24b}$ are nonisomorphic and that $\Gamma_{24a}$ is isomorphic to $\mathbb{W}(F_4)$. $\square$

**Proposition 2.47. (joint with Gramlich, Hall)** *Let $\Gamma$ be a connected bichromatic graph that is locally like $\mathbb{W}(F_4)$. If $\Gamma$ is tightly connected then $\Gamma \cong \mathbb{W}(F_4)$ or $\Gamma \cong \Gamma_{24b}$.*

**Proof.** Fix a short 4-clique $x_1, x_2, x_3, x_4$. Because of tightness every long vertex is adjacent to one of the $x_i$, and by Proposition 2.32 there are exactly 12 such long vertices. Thus $\Gamma$ consists of 12 long vertices. Likewise there are exactly 12 short vertices. Hence $|\Gamma| = 24$, and the claim follows from Proposition 2.46. $\square$

**Proposition 2.48.** *Let $\Gamma$ be a connected bichromatic graph that is locally like $\mathbb{W}(F_4)$. If $\Gamma$ has diameter 2 then $\Gamma \cong \mathbb{W}(F_4)$ or $\Gamma \cong \Gamma_{24b}$.*

**Proof.** Let $x_1, x_2, x_3, x_4$ be a short 4-clique. As in Proposition 2.32 let $y_{i,j}, y_{j,i}$ be the long vertices adjacent to both $x_i$ and $x_j$. Assume that there is a long vertex $v$ which is not among the 12 long vertices $y_{i,j}$. Because $v$ is not adjacent to any of the $x_i$ and since the diameter of $\Gamma$ is 2 we find a long vertex that connects $x_1$ and $v$. Without loss of generality let this long vertex be $y_{1,2}$. This prevents $y_{1,2}, y_{2,1}, y_{3,4}, y_{4,3}$ from forming a long 4-clique. By Proposition 2.32 there are thus long vertices $v_1, v_2$ not among the $y_{i,j}$ such that $y_{3,4}, y_{4,3}, v_1, v_2$ form a long 4-clique. Again, $v_1$ is not adjacent to any of the $x_i$ and hence is connected to $x_1$ by a long vertex. This is a contradiction since the long vertices adjacent to $x_1$ are the vertices $y_{1,j}, y_{j,1}$.

Consequently, $\Gamma$ contains no further long vertices besides the 12 vertices $y_{i,j}$. The same reasoning reveals that $\Gamma$ contains exactly 12 short vertices. We conclude that $|\Gamma| = 24$. According to Proposition 2.46 this proves that $\Gamma \cong \mathbb{W}(F_4)$ or $\Gamma \cong \Gamma_{24b}$. □

**Proposition 2.49.** *Let $\Gamma$ be a connected bichromatic graph that is locally like $\mathbb{W}(F_4)$. Let $x_1$ be a short vertex of $\Gamma$. If for any nonadjacent pair of long vertices $y$, $y' \in x_1^\perp$ we have $\mu(y, y') = 3$ and for any nonadjacent pair of a short and a long vertex $x, y \in x_1^\perp$ we have $\mu_\ell(x, y) = 1$ then $\Gamma \cong \mathbb{W}(F_4)$ or $\Gamma \cong \Gamma_{24b}$.*

**Proof.** Let $x_1, x_2, x_3, x_4$ be the short 4-clique containing $x_1$. As in Proposition 2.32 denote with $y_{i,j}$, $y_{j,i}$ the long vertices adjacent to both $x_i$ and $x_j$. The long vertex $y_{1,2}$ and the short vertex $x_3$ are nonadjacent neighbors of $x_1$. Since $\mu_\ell(x_3, y_{1,2}) = 1$ we find two long vertices in $\{x_3, y_{1,2}\}^\perp$. By Proposition 2.32 the only candidates are $y_{3,4}$ and $y_{4,3}$. Consequently, $y_{1,2}, y_{2,1}, y_{3,4}, y_{4,3}$ form a long 4-clique. Likewise, by considering the nonadjacent pairs $y_{1,3}, x_4$ respectively $y_{1,4}, x_2$ we find that $y_{1,3}, y_{3,1}, y_{2,4}, y_{4,2}$ respectively $y_{1,4}, y_{4,1}, y_{2,3}, y_{3,2}$ form a long 4-clique.

Consider the long 4-clique $y_{1,2}, y_{2,1}, y_{3,4}, y_{4,3}$. It follows from Proposition 2.32 that there are 8 short neighbors of this clique besides $x_1, x_2, x_3, x_4$. Denote these neighbors by $x_5, x_6..., x_{12}$. We may assume that

$$\{x_5, x_6, y_{1,2}, y_{3,4}\}, \quad \{x_7, x_8, y_{2,1}, y_{4,3}\}, \quad \{x_9, x_{10}, y_{1,2}, y_{4,3}\}, \quad \{x_{11}, x_{12}, y_{2,1}, y_{3,4}\}$$

form mixed 4-cliques. Looking at the neighborhood of $y_{1,2}$ we see that $x_5$ is not adjacent to $x_{10}$. Likewise, $x_5$ is not adjacent to $x_{12}$.

Now, consider the long vertex $y_{1,3}$. By assumption $\{y_{1,2}, y_{1,3}\}^\perp$ contains 2 short vertices besides $x_1$. Since these two short vertices are not adjacent we may assume them to be $x_5$ and $x_9$. Likewise, $y_{2,1}$ and $y_{1,3}$ share two short neighbors besides $x_1$ which we may assume to be $x_7$ and $x_{11}$. Summarizing, $y_{1,3}$ is adjacent to $x_5, x_7, x_9$ and $x_{11}$. Because $\{y_{1,3}, y_{3,1}\}^\perp$ contains no short vertex besides $x_1$ and $x_3$, we analogously find that $y_{3,1}$ is adjacent to $x_6, x_8, x_{10}$ and $x_{12}$. Because the short neighbors of $y_{1,3}$ and $y_{3,1}$ come in adjacent pairs we deduce that $x_5, x_6, x_7, x_8$ as well as $x_9, x_{10}, x_{11}, x_{12}$ form a short 4-clique.

In particular, we just showed that the six short neighbors of $y_{1,3}$ are found among $x_1$, $x_2, ..., x_{12}$. The same argument applies analogously to the other $y_{1,j}$ and $y_{i,1}$. Consider for instance the vertex $y_{2,4}$. Exploiting that $y_{2,4}$ is adjacent to both $y_{1,3}$ and $y_{3,1}$, and thus shares two short neighbors with each, we find that the six short neighbors of $y_{2,4}$ are among the $x_i$. Accordingly, we conclude that any short vertex adjacent to one of the 12 long vertices of the form $y_{i,j}$ is among the short vertices $x_1$, $x_2, ..., x_{12}$. By reciprocity as in Proposition 2.37 we find that all long neighbors of one of the 12 short vertices $x_1, x_2, ..., x_{12}$ are among the $y_{i,j}$. Since all of these vertices came in 4-cliques of the same type the connectedness of $\Gamma$ implies that $|\Gamma| = 24$. The claim follows by Proposition 2.46. □

**Corollary 2.50.** *Let $\Gamma$ be a connected bichromatic graph that is locally like $\mathbb{W}(F_4)$. If $\mu = 3$ and $\mu_s = \mu_\ell = 1$ then $\Gamma \cong \mathbb{W}(F_4)$ or $\Gamma \cong \Gamma_{24b}$.*

**Proof.** Just pick any short vertex $x_1 \in \Gamma$ and apply Proposition 2.49. □

**Proposition 2.51.** *Let $\Gamma$ be a bichromatic graph that is locally like $\mathbb{W}(F_4)$. If $\mu = 3$ then $\mu_s = \mu_\ell = 1$.*

**Proof.** Let $\Pi$ be the partition of $\Gamma$ into 4-cliques of the same type, and let $X = \{x_1, x_2, x_3, x_4\}$ be a short vertex of $\Gamma/\Pi$. As in Proposition 2.32 denote with $y_{i,j}$, $y_{j,i}$ the long vertices adjacent to both $x_i$ and $x_j$. Recall that $y_{i,j}$ and $y_{k,l}$ are at distance 2 if $|\{i, j\} \cap \{k, l\}| = 1$. Note that a long vertex has six short neighbors which come in adjacent pairs. Further, because $\mu = 3$ two long vertices at distance 2 share three short neighbors. Therefore, if $y$ is a long vertex neighbored to the adjacent pair of short vertices $x$, $x'$ and if $y'$ is a long vertex at distance 2 from $y$ then $y'$ is neighbored to exactly one $x$, $x'$.

Denote with $z_1$, $z_2$ a pair of adjacent short vertices besides $x_1$, $x_2$ neighbored to $y_{1,2}$. The vertices $y_{1,2}$ and $y_{1,3}$ are at distance 2. Without loss we may therefore assume that $y_{1,3}$ is adjacent to $z_1$. Let $z_3$ be the unique short vertex adjacent to $z_1$ and $y_{1,3}$, and let $z_4$ be the unique short vertex such that $z_1, z_2, z_3, z_4$ form a short 4-clique. Let $y$ be a long vertex at distance 2 from both $y_{1,2}$ and $y_{1,3}$. Then $y$ is either adjacent to both $z_1$ and $z_4$ or $y$ is adjacent to both $z_2$ and $z_3$. The vertices $y_{1,4}$ and $y_{4,1}$ are each long vertices at distance 2 from both $y_{1,2}$ and $y_{1,3}$. $y_{1,4}$ and $y_{4,1}$ are adjacent and thus share no short neighbors besides $x_1$ and $x_4$. Accordingly, we may assume that $y_{1,4}$ is adjacent to both $z_1$ and $z_4$, and that $y_{4,1}$ is adjacent to both $z_2$ and $z_3$. The following table summarizes the situation.

$$\begin{array}{ll} y_{1,2} \perp z_1, z_2 & y_{2,1} \perp z_3, z_4 \\ y_{1,3} \perp z_1, z_3 & y_{3,1} \perp z_2, z_4 \\ y_{1,4} \perp z_1, z_4 & y_{4,1} \perp z_2, z_3 \end{array}$$

Consider the vertex $y_{2,3}$. $y_{2,3}$ has distance 2 from both $y_{1,2}$ and $y_{1,3}$. Therefore either $y_{2,3} \perp z_1, z_4$ or $y_{2,3} \perp z_2, z_3$. In the former case, $y_{2,3}, y_{1,4} \in \{z_1, z_4\}^\perp$. Since $z_1$ and $z_4$ are adjacent the vertices $y_{2,3}$ and $y_{1,4}$ are adjacent as well. Likewise, the latter case implies that $y_{2,3}$ and $y_{4,1}$ are adjacent. In both cases, the long vertices $y_{2,3}, y_{3,2}, y_{1,4}, y_{4,1}$ form a 4-clique.

Analogously, one shows that the vertices $y_{i,j}, y_{j,i}, y_{k,l}, y_{l,k}$ form a 4-clique whenever the index sets $\{i, j\}$ and $\{k, l\}$ are disjoint. The local graph at $X$ in $\Gamma/\Pi$ therefore is isomorphic to $\overline{K_3}^\ell$. Since $X$ was an arbitrary short vertex, and since the same argument works for long vertices, we proved that $\Gamma/\Pi$ is locally homogeneous with $\Delta_s(\Gamma/\Pi) \cong \overline{K_3}^\ell$ and $\Delta_\ell(\Gamma/\Pi) \cong \overline{K_3}^s$. It follows from Proposition 2.42 that $\mu_s = \mu_\ell = 1$.  $\square$

**Corollary 2.52.** *Let $\Gamma$ be a connected bichromatic graph that is locally like $\mathbb{W}(F_4)$. If $\mu = 3$ then $\Gamma \cong \mathbb{W}(F_4)$ or $\Gamma \cong \Gamma_{24b}$.*

**Proof.** This is a consequence of Corollary 2.50 together with Proposition 2.51.  $\square$

**Proposition 2.53.** *Let $\Gamma$ be one of the graphs $\mathbb{W}(F_4)$ or $\Gamma_{24b}$. Then*

$$\mathrm{Aut}(\Gamma) \cong W(F_4)/Z.$$

*Furthermore, $\Gamma$ is transitive on vertices of the same type.*

**Proof.** This is proved in the appendix by Proposition A.3, Proposition A.4, Proposition A.5 and Proposition A.10. □

### 2.5.3 The parameters $\mu$, $\mu_s$, $\mu_\ell$

In Proposition 2.41 we associated the parameters $\mu$, $\mu_s$, $\mu_\ell$ to a graph that is locally like $\mathbb{W}(F_4)$. In light of Theorem 2.45 one is interested in conditions under which $\mu_s + \mu_\ell$ respectively $\mu$ are constant. Recall the terminology, introduced in Example 1.56 and Example 1.57, of a group acting transitively on triangles respectively paths of a graph.

**Proposition 2.54.** *Let $\Gamma$ be locally like $\mathbb{W}(F_4)$ and transitive on short respectively long oriented triangles. Then $\mu_s + \mu_\ell$ is constant.*

**Proof.** By construction, $\mu_s + \mu_\ell$ only takes the values 1 and 2. Suppose that $\mu_s + \mu_\ell$ is not constantly 1, and let $x$ and $y$ be nonadjacent short and long vertices such that $\mu_s(x, y) = \mu_\ell(x, y) = 1$. Let $x_1$, $x_2$ respectively $y_1$, $y_2$ be the short respectively long vertices contained in $\{x, y\}^\perp$. Denote with $y_0$ the unique long vertex contained in $\{x_1, x_2\}^\perp$ besides $y$. Accordingly, $y_0$ is adjacent to $y_1$ and $y_2$, and we have the following induced subgraph.

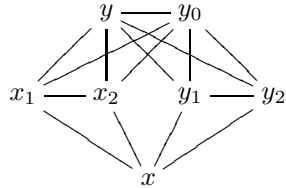

Let $x'$ and $y'$ be short and long vertices at distance 2. Without loss we may assume that $\mu_s(x', y') = 1$. Let $x'_1, x'_2$ be the short vertices contained in $\{x', y'\}^\perp$. Since $\Gamma$ is transitive on short oriented triangles there is an automorphism $\psi$ of $\Gamma$ such that

$$\psi(x, x_1, x_2) = (x', x'_1, x'_2).$$

Note that $y'$ equals either $\psi(y)$ or $\psi(y_0)$ because $y'$ is one of the two long vertices contained in $\{x'_1, x'_2\}^\perp = \{\psi(x_1), \psi(x_2)\}^\perp$. In both cases, we find that $\{x', y'\}^\perp$ contains the long vertices $\psi(y_1)$ and $\psi(y_2)$. Thus $\mu_\ell(x', y') = 1$. □

**Corollary 2.55.** *Let $\Gamma$ be locally like $\mathbb{W}(F_4)$ and transitive on short respectively long oriented triangles. Then $\mu_s = \mu_\ell = 1$ if and only if $\Gamma$ contains an induced subgraph isomorphic to*

.

**Proof.** Let $x$ be a short and $y$ a long vertex of this induced subgraph that are nonadjacent. Then $x$ and $y$ are at distance 2 in $\Gamma$ as well and $(\mu_s + \mu_\ell)(x, y) = 2$. By Proposition 2.54 this implies $\mu_s + \mu_\ell = 2$ or, equivalently, $\mu_s = \mu_\ell = 1$. □

**Proposition 2.56.** *Let $\Gamma$ be locally like $\mathbb{W}(F_4)$ and transitive on 3-paths of the same type. Then $\mu$ is constant on pairs of short vertices as well as on pairs of long vertices.*

**Proof.** Recall that $\mu$ takes values in $\{1, 2, 3\}$. Let $x_1$, $x_2$ be nonadjacent short vertices at distance 2 such that $\mu$ restricted to pairs of short vertices is maximal. Let $y \in \{x_1, x_2\}^\perp$. Consider another nonadjacent pair of short vertices $x'_1, x'_2$ at distance 2 and let $y' \in \{x'_1, x'_2\}^\perp$. $y$ and $y'$ are both long whence by assumption we find an automorphism $\psi \in \mathrm{Aut}(\Gamma)$ such that $\psi(x_1, y, x_2) = (x'_1, y', x'_2)$ or $\psi(x_1, y, x_2) = (x'_2, y', x'_1)$. Since $\{x'_1, x'_2\}^\perp = \psi(\{x_1, x_2\}^\perp)$ this proves that $\mu$ is constant on pairs of short vertices. Likewise for pairs of long vertices. □

To give a criterion under which $\mu$ is in fact constant we introduce the concept of maximal transitivity on neighbors. Let $\Gamma$ be a graph, and $x \in \Gamma$. Notice that the stabilizer $C_{\mathrm{Aut}(\Gamma)}(x)$ induces an action on the neighbors of $x$.

**Definition 2.57.** *A graph $\Gamma$ is said to be* maximally transitive on neighbors *if for each $x \in \Gamma$ the stabilizer $C_{\mathrm{Aut}(\Gamma)}(x)$ induces on $x^\perp$ all the automorphisms of the induced subgraph on $x^\perp$.*

**Remark 2.58.** Notice that for a graph to be maximally transitive on neighbors it is not sufficient that the stabilizer $C_{\mathrm{Aut}(\Gamma)}(x)$ is isomorphic to $\mathrm{Aut}(x^\perp)$ for each $x \in \Gamma$. In particular, we will see in Proposition A.9 that the Weyl graph $\mathbb{W}(F_4)$ is not maximally transitive on neighbors even though the automorphisms stabilizing a vertex are isomorphic to $W(B_3)$. See also Remark 3.11. We will, however, construct two graphs in the course of the proof of Theorem 2.62 which are locally like $\mathbb{W}(F_4)$ and which are maximally transitive on neighbors.

**Proposition 2.59.** *Let $\Gamma$ be locally like $\mathbb{W}(F_4)$. If $\Gamma$ is transitive and maximally transitive on neighbors then it is transitive on oriented triangles of the same type and on oriented 3-paths of the same type.*

**Proof.** Let $(x_1, x_2, x_3)$ and $(y_1, y_2, y_3)$ be two oriented triangles of the same type. We may assume that $x_1$ and $y_1$ are short vertices. Since $\Gamma$ is transitive we find an automorphism $\psi \in \mathrm{Aut}(\Gamma)$ such that $\psi(x_1) = y_1$. Let $y'_2 = \psi(x_2)$ and $y'_3 = \psi(x_3)$. Then $\psi(x_1, x_2, x_3) = (y_1, y'_2, y'_3)$ is an oriented triangle of the same type as $(x_1, x_2, x_3)$. By assumption, the induced subgraph on $x_1^\perp$ is isomorphic to $\mathbb{W}(B_3)$. Observe that the graph $\mathbb{W}(B_3)$ is transitive on oriented 2-paths of the same type. Maximal transitivity on neighbors applied to $x_1^\perp$ implies that we find an automorphism $\varphi \in \mathrm{Aut}(\Gamma)$ stabilizing $x_1$ such that $\varphi(y'_2, y'_3) = (y_2, y_3)$. Hence,

$$\varphi \circ \psi(x_1, x_2, x_3) = (y_1, y_2, y_3)$$

which proves transitivity on oriented triangles of the same type. Likewise for oriented 3-paths of the same type. □

**Proposition 2.60.** *Let $\Gamma$ be locally like $\mathbb{W}(F_4)$, transitive and maximally transitive on neighbors. Then $\mu$ is constant.*

**Proof.** By Proposition 2.56 and Proposition 2.59 $\mu$ is constant on pairs of short vertices as well as on pairs of long vertices.

Suppose that $\mu$ is not constantly 1. Then we find two nonadjacent vertices $x_1$, $x_2$ of the same type such that $\{x_1, x_2\}^\perp$ contains at least two distinct vertices $y_1$, $y_2$. Without loss we may assume $x_1, x_2$ to be short, and $y_1, y_2$ to be long. By Proposition 2.41 the vertices $y_1, y_2$ are nonadjacent. Clearly, $x_1, x_2 \in \{y_1, y_2\}^\perp$. Thus $\mu(x_1, x_2) \geqslant 2$ as well as $\mu(y_1, y_2) \geqslant 2$. Consequently, $\mu \geqslant 2$.

Finally, suppose that $\mu$ was not constant. Then $\mu \geqslant 2$, and $\mu$ is constant on pairs of short vertices and constant on pairs of long vertices. Without loss we assume that $\mu = 3$ on short vertices and $\mu = 2$ on long vertices. Let $x_1$, $x_2$ be nonadjacent short vertices and let $y_1$, $y_2$, $y_3$ be the three long vertices making up $\{x_1, x_2\}^\perp$. Let $y_0$ be the unique long vertex in $x_1^\perp$ adjacent to $y_1$. Since $y_0$ is not adjacent to $x_2$ we have the following induced subgraph.

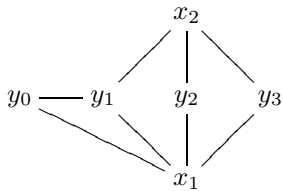

By maximal transitivity on neighbors of $x_1$ we find an automorphism $\psi$ such that

$$\psi(x_1, y_1, y_2, y_3) = (x_1, y_0, y_2, y_3).$$

Since $\mu = 2$ on long vertices, $x_2$ is the unique vertex in $\{y_2, y_3\}^\perp$ besides $x_1$. Therefore $\psi(x_2) = x_2$ which implies that $y_0 \perp x_2$ and hence $y_0 \in \{x_1, x_2\}^\perp$. This is a contradiction. □

### 2.5.4 Two further graphs locally like $\mathbb{W}(F_4)$

In Corollary 2.50 we recognized $\mathbb{W}(F_4)$ and $\Gamma_{24b}$ as the only connected bichromatic graphs that are locally like $\mathbb{W}(F_4)$ and for which $\mu = 3$ and $\mu_s = \mu_\ell = 1$ where $\mu$, $\mu_s$, $\mu_\ell$ are the parameters introduced in Proposition 2.41. This characterizes these two graphs as the tightest graphs that are locally like $\mathbb{W}(F_4)$. Furthermore, we showed in Proposition 2.51 that there is no graph locally like $\mathbb{W}(F_4)$ for which $\mu = 3$ and $\mu_s + \mu_\ell = 1$.

**Remark 2.61.** The next level of tightness in terms of the parameters $\mu$, $\mu_s$, $\mu_\ell$ is therefore achieved by graphs with $\mu = 2$ and $\mu_s = \mu_\ell = 1$. On the other hand, we characterized $\mathbb{W}(F_4)$ and $\Gamma_{24b}$ as the only bichromatic graphs on 24 vertices that are locally like $\mathbb{W}(F_4)$. In view of Corollary 2.38 the next level of tightness in terms of the number of vertices is achieved by graphs on 32 vertices.

Let $\Gamma$ be a bichromatic graph on 32 vertices, locally like $\mathbb{W}(F_4)$ with parameters $\mu_s = \mu_\ell = 1$, and let $\Pi$ be the partition of $\Gamma$ into 4-cliques. According to Proposition 2.42, $\Gamma/\Pi$ is locally homogeneous with short local graph $\overline{K_3}^\ell$ and long local graph $\overline{K_3}^s$. By construction, $|\Gamma| = 4 \cdot |\Gamma/\Pi|$, and hence $|\Gamma/\Pi| = 8$. The reader is invited to check that $\Gamma/\Pi$ is therefore, as an unlabeled graph, isomorphic to the 1-skeleton of a cube.

We close this section by characterizing two particularly symmetric graphs with parameters $\mu = 2$ and $\mu_s = \mu_\ell = 1$. Denote with $\Gamma_{32a}$ and $\Gamma_{32b}$ the two connected graphs on 32 vertices constructed in the course of the proof of Theorem 2.62.

**Theorem 2.62.** *Let $\Gamma$ be locally like $\mathbb{W}(F_4)$ and maximally transitive on neighbors. If $\mu = 2$ and $\mu_s = \mu_\ell = 1$, then $\Gamma \cong \Gamma_{32a}$ or $\Gamma \cong \Gamma_{32b}$.*

**Proof.** Let $x_1, x_2, x_3, x_4$ be a short 4-clique of $\Gamma$. Denote the long neighbors of $x_1$ by $y_1, y_2, ..., y_6$ and assume that $y_1 \perp y_2$, $y_3 \perp y_4$ and $y_5 \perp y_6$. Further, assume that $x_2 \perp y_1$, $x_3 \perp y_3$ and $x_4 \perp y_5$. Since $\mu = 2$, each nonadjacent pair $y_i, y_j$ has exactly one short neighbor in common besides $x_1$. Denote this short vertex by $z_{i,j}$. Since $\mu(x_1, z_{i,j}) = 2$, we have $\{x_1, z_{i,j}\}^\perp = \{y_i, y_j\}$. In particular, the $z_{i,j}$ with $i < j$ are 12 distinct short vertices.

We claim that $z_{i,j} \perp z_{i,k}$ if and only if $y_j \perp y_k$. Without loss we may assume that $i = 1$. The subgraph induced on the $z_{1,j}$ is a union of two 2-cliques. Suppose that $z_{1,3} \perp z_{1,5}$. By maximal transitivity on the neighbors of $x_1$ we find an automorphism of $\Gamma$ that switches $y_5$ and $y_6$ while fixing every other neighbor of $x_1$. Under this automorphism the vertices $z_{1,5}$ and $z_{1,6}$ are switched while the other $z_{1,j}$ are fixed. We therefore get the contradiction $z_{1,3} \perp z_{1,6}$. Hence $z_{1,3} \perp z_{1,4}$, and therefore $z_{1,5} \perp z_{1,6}$ as well.

This describes the subgraph of $\Gamma$ induced on the vertices $x_i$, $y_i$ and $z_{i,j}$. In particular, the vertices $z_{i,j}$ with $i \in \{i_0, i_0 + 1\}$ and $j \in \{j_0, j_0 + 1\}$ form 4-cliques for $i_0 \neq j_0 \in \{1, 3, 5\}$. Let $i \neq j \in \{1, 3, 5\}$. The vertices $y_i, z_{i,j}, z_{i,j+1}$ are pairwise adjacent. Denote with $y_{i,j}$ the unique (long) vertex adjacent to these three vertices.

By construction, $y_{1,3}$ is adjacent to $y_1$, $z_{1,3}$ and $z_{1,4}$. Likewise, $y_{1,5}$ is adjacent to $y_1$, $z_{1,5}$ and $z_{1,6}$. The long vertices $y_{1,3}$ and $y_1$ only have two common short neighbors which implies that $y_{1,3} \neq y_{1,5}$. Hence $y_1, y_2, y_{1,3}, y_{1,5}$ form a 4-clique. Since $z_{1,3}, z_{1,4}, z_{2,3}, z_{2,4}$ form a 4-clique we find that $\{y_{1,3}, y_2\}^\perp = \{z_{2,5}, z_{2,6}\}$. $y_{1,3}$ and $x_1$ are at distance 2 in $\Gamma$. Since $\mu_s(x_1, y_{1,3}) = 1$ and $x_2 \perp y_1$, the vertices $x_3$ and $x_4$ are adjacent to $y_{1,3}$. This determines $y_{1,3}^\perp$. Analogously, we find

$$
\begin{aligned}
y_{1,3}^\perp &= \{y_1, y_2, y_{1,5},\ x_3, x_4,\ z_{1,3}, z_{1,4},\ z_{2,5}, z_{2,6}\}, \\
y_{1,5}^\perp &= \{y_1, y_2, y_{1,3},\ x_3, x_4,\ z_{1,5}, z_{1,6},\ z_{2,3}, z_{2,4}\}, \\
y_{3,1}^\perp &= \{y_3, y_4, y_{3,5},\ x_2, x_4,\ z_{1,3}, z_{2,3},\ z_{4,5}, z_{4,6}\}, \\
y_{3,5}^\perp &= \{y_3, y_4, y_{3,1},\ x_2, x_4,\ z_{3,5}, z_{3,6},\ z_{1,4}, z_{2,4}\}, \\
y_{5,1}^\perp &= \{y_5, y_6, y_{5,3},\ x_2, x_3,\ z_{1,5}, z_{2,5},\ z_{3,6}, z_{4,6}\}, \\
y_{5,3}^\perp &= \{y_5, y_6, y_{5,1},\ x_2, x_3,\ z_{3,5}, z_{4,5},\ z_{1,6}, z_{2,6}\}.
\end{aligned}
$$

That determines the induced subgraph $\Gamma_2$ on the vertices $x_i, y_i, z_{i,j}, y_{i,j}$.

Consider the short 4-clique $z_{1,3}, z_{1,4}, z_{2,3}, z_{2,4}$. We already saw that $\{z_{1,3}, z_{1,4}\}^\perp$ contains the long vertices $y_1, y_{1,3}$ and that $\{z_{1,3}, z_{2,3}\}^\perp$ contains the long vertices $y_3, y_{3,1}$. Denote with $w_1, w_2$ the pair of adjacent long vertices contained in $\{z_{1,3}, z_{2,4}\}^\perp$, and with $w_3, w_4$ the pair of adjacent long vertices contained in $\{z_{1,4}, z_{2,3}\}^\perp$. Since $z_{1,3}$ and $z_{1,4}$ have no common long neighbors the vertices $w_i$ are distinct.

Now, consider the short 4-clique $z_{1,5}$, $z_{1,6}$, $z_{2,5}$, $z_{2,6}$. Note that every vertex of this short 4-clique is at distance 2 in $\Gamma_2$ from every vertex of the short 4-clique $z_{1,3}$, $z_{1,4}$, $z_{2,3}$, $z_{2,4}$ with exactly one vertex in $\Gamma_2$ connecting each pair. In particular, the vertices $z_{1,5}$ and $z_{1,3}$ are at distance 2 in $\Gamma$ and hence share exactly one long vertex besides $y_1$, say $w_1$ (the only other candidate is $w_2$). Hence $w_2$ is contained in $\{z_{1,6}, z_{1,3}\}^\perp$. Likewise, we may assume that $z_{1,5}$ and $z_{1,4}$ share $w_3$ while $z_{1,6}$ and $z_{1,4}$ are both adjacent to $w_4$. In particular, the $w_i$ form a long 4-clique. Further, $z_{2,6} \perp w_1, w_3$ and $z_{2,5} \perp w_2, w_4$.

Finally, consider the short 4-clique $z_{3,5}$, $z_{3,6}$, $z_{4,5}$, $z_{4,6}$. We have two choices for the two long vertices contained in $\{z_{3,5}, z_{4,6}\}^\perp$. Either $w_1, w_4 \in \{z_{3,5}, z_{4,6}\}^\perp$ or $w_2, w_3 \in \{z_{3,5}, z_{4,6}\}^\perp$. In both cases, the vertices of $\Gamma$ are precisely the 32 vertices $x_i$, $y_i$, $z_{i,j}$, $y_{i,j}$, $w_i$ and $\Gamma$ is determined by this final choice. Let $\Gamma_{32a}$ be the graph corresponding to the choice $w_1, w_4 \in \{z_{3,5}, z_{4,6}\}^\perp$, and $\Gamma_{32b}$ be the graph corresponding to the choice $w_2, w_3 \in \{z_{3,5}, z_{4,6}\}^\perp$. The following table summarizes adjacency involving the vertices $w_1, w_2, w_3, w_4$.

| | | |
|---|---|---|
| by construction | $w_1, w_2 \perp z_{1,3}, z_{2,4}$ | $w_3, w_4 \perp z_{1,4}, z_{2,3}$ |
| | $w_1, w_3 \perp z_{1,5}, z_{2,6}$ | $w_2, w_4 \perp z_{1,6}, z_{2,5}$ |
| $\Gamma_{32a}$ | $w_1, w_4 \perp z_{3,5}, z_{4,6}$ | $w_2, w_3 \perp z_{3,6}, z_{4,5}$ |
| $\Gamma_{32b}$ | $w_1, w_4 \perp z_{3,6}, z_{4,5}$ | $w_2, w_3 \perp z_{3,5}, z_{4,6}$ |

In the appendix we implement the graphs $\Gamma_{32a}$ and $\Gamma_{32b}$ in SAGE. In particular, we verify in Proposition A.6 using SAGE that $\Gamma_{32a}$ and $\Gamma_{32b}$ are nonisomorphic. In Proposition A.8 we prove that $\Gamma_{32a}$ and $\Gamma_{32b}$ are indeed maximally transitive on neighbors. □

# 3 Group theoretic applications

Let $G$ be a group acting on a graph $\Gamma$. In this section we explore how a local recognition result for the graph $\Gamma$ may be turned into a group theoretical statement about $G$. In fact, what we seek to do is to encode the recognition of the local structure of $\Gamma$ into a statement about the local structure of $G$. Here, the local structure of $G$ vaguely means a statement involving only few particular elements of $G$ and subgroups related to them but not the whole group $G$. Note that the graph theoretical recognition results we considered in this text were mostly concerned with graphs $\Gamma$ that are commuting graphs as introduced in Definition 1.37.

**Remark 3.1.** Let $G$ be a group, $x \in G$, and consider the commuting graph $\Gamma$ of $G$ on $x^G$. Recall from Definition 1.37 that this means that the vertices of $\Gamma$ are the $G$-conjugates of $x$ with two elements adjacent whenever they commute. The action of $G$ on $\Gamma$ is transitive which in particular means that $\Gamma$ is locally homogeneous. By construction, the neighborhood $x^\perp$ of $x$ in $\Gamma$ consists of the $G$-conjugates of $x$ which are contained in the centralizer $C_G(x)$, that is

$$\{x\} \cup x^\perp = x^G \cap C_G(x).$$

The conjugates of $x$, however, are considered a global information because it involves knowledge about the whole group $G$.

Let $y \in x^\perp$. We consider the following subgraph of $\Gamma$. Let $\Gamma'$ be the graph on the conjugates of $x$ with edge set $E(\Gamma')$ given by $\{x, y\}^G$. In other words, $\Gamma'$ has the same vertices as $\Gamma$ but two vertices $x_1, x_2$ are adjacent only if $\{x_1, x_2\} = \{x, y\}^g$ for some $g \in G$. Note that $\Gamma'$ is indeed a subgraph of the commuting graph $\Gamma$ because $x$ and $y$ commute. By construction, $G$ acts transitively on $\Gamma'$. Let us further assume that $(y, x) = (x, y)^u$ for some $u \in G$. This is the case for instance if $x$ and $y$ are conjugate by an involution. Then two vertices $x_1, x_2$ are adjacent if and only if $(x_1, x_2) = (x, y)^g$ for some $g \in G$. Consequently, the neighbors of $x$ in $\Gamma'$ are

$$x^\perp = \{y^g \colon g \in C_G(x)\}.$$

We consider this description local because it is in terms of the two elements $x$, $y$ and the centralizer $C_G(x)$ only.

In our applications the commuting graph $\Gamma$ is edge-transitive. Accordingly, the graphs $\Gamma'$ and $\Gamma$ agree.

**Example 3.2.** The commuting graph $\Gamma$ of the symmetric group $\mathrm{Sym}_n$ on its transpositions is the Kneser graph $K(n, 2)$, see Example 1.39. Let $n \geqslant 4$ and $x$, $y$ be two commuting transpositions. Define $\Gamma'$ to be the graph on the transpositions of $\mathrm{Sym}_n$ with edges $\{x, y\}^{\mathrm{Sym}_n}$. Exploiting that $\Gamma$ is edge-transitive, $\Gamma' = \Gamma$.

The following observations concerning graphs such as the graph $\Gamma'$ considered in Remark 3.1 are stated in a slightly more general setting so that they will be of use as well when we generalize to bichromatic graphs later on.

**Proposition 3.3.** *Let $G$ be a group, $x, y \in G$, and $\Gamma$ a graph with vertices the conjugates of $x$ and $y$ such that $\{x, y\}^G \subset E(\Gamma)$.*

- *If $G = \langle C_G(x), C_G(y) \rangle$ then $\Gamma$ is connected.*
- *If $G = \langle x^G, y^G \rangle$ and $G$ acts on $\Gamma$ by conjugation then the kernel of this action is $Z(G)$.*

**Proof.** For the first claim, suppose that $G = \langle C_G(x), C_G(y) \rangle$. Denote with $\Lambda$ the connected component of $\Gamma$ containing $x$ and $y$. Let $g \in G$ such that $x^g$ and $y^g$ are contained in $\Lambda$. Let $h \in C_G(x)$. Consequently, $x^{hg} = x^g \in \Lambda$. Moreover,

$$\{x, y\}^{hg} = \{x^g, y^{hg}\}$$

which implies that $y^{hg} \perp x^g$ in $\Gamma$. Hence $y^{hg}$ is contained in $\Lambda$ as well. The same argument applies to $h \in C_G(y)$. By assumption, any $g \in G$ can be written as $g = h_1 h_2 \cdots h_n$ with each $h_i$ contained either in $C_G(x)$ or $C_G(y)$. Connectedness of $\Gamma$ follows by induction.

Now, assume that $G = \langle x^G, y^G \rangle$ and that $G$ acts on $\Gamma$ by conjugation. Clearly, $Z(G)$ is contained in the kernel of the action. On the other hand, let $g \in G$ be an element acting trivially. This means that $g$ centralizes the conjugates of $x$ and $y$. But these generate $G$ and the claim follows. $\square$

## 3.1 Recognizing Sym$_n$

The following theorem provides a local characterization of the symmetric groups Sym$_n$. It is the paradigmatic example given in [GLS94, Theorem 27.1] to illustrate the strategy of recognizing a group from local information based on centralizers of involutions. The proof given here is based on the local recognition of Kneser graphs stated in Theorem 2.8. Actually, it is this characterization of the symmetric groups that was a motivation for Jonathan I. Hall to pursue the local recognition of Kneser graphs, see [Hal87].

**Theorem 3.4.** *Let $n \geqslant 7$, and let $G$ be a group with involutions $x, y$ such that*

- $C_G(x) = \langle x \rangle \times J$ *with* $J \cong \text{Sym}_n$,
- $C_G(y) = \langle y \rangle \times K$ *with* $K \cong \text{Sym}_n$,
- *$x$ is a transposition in $K$,*
- *$y$ is a transposition in $J$,*
- *$J \cap K$ contains an involution $z$ that is a transposition in both $J$ and $K$.*

*If $G = \langle J, K \rangle$, then $G \cong \text{Sym}_{n+2}$.*

**Proof.** $y$ and $z$ are both transpositions in $J$ and hence conjugate by an involution $u \in J$. Accordingly, $(x, y, z)^u = (x, z, y)$. Likewise, $x$ and $z$ are conjugate by an involution $v \in K$. We conclude that $x, y, z$ are all conjugate in $G$, and that $\langle u, v \rangle$ acts by conjugation on the set $\{x, y, z\}$ as $\mathrm{Sym}_3$. In particular, we find $w$ such that $(y, x) = (x, y)^w$. Let $\Gamma$ be the graph on the $G$-conjugates of $x$ with edges $\{x, y\}^G$. Note that two vertices $a, b$ in $\Gamma$ are adjacent if and only if $(a, b) = (x, y)^g$ for some $g \in G$. By Proposition 3.3 the graph $\Gamma$ is connected. The neighbors of $x$ are the $J$-conjugates of $y$ since

$$x^\perp = \{y^g \colon g \in C_G(x)\} = \{y^g \colon g \in J\}.$$

By construction, $\{x, y, z\}$ is a triangle of $\Gamma$. We claim that $G$ is transitive on the oriented triangles of $\Gamma$. Indeed, let $(a, b, c)$ be an oriented triangle of $\Gamma$. By edge-transitivity we find $g \in G$ such that $(a, b) = (x, y)^g$. Set $d = z^g$. Because $b, c, d$ are all neighbors of $a$ they are $J^g$-conjugates in $J^g$. Consequently, $b, c, d$ are transpositions in $J^g \cong \mathrm{Sym}_n$. Since $[b, c] = [b, d] = 1$ we find $h \in J^g$ such that $(b, c) = (b, d)^h$. Thus $(a, b, c) = (a, b, d)^h$ which shows that $(a, b, c)$ is indeed conjugate to $(x, y, z)$.

We observed that the neighbors of $x$ are the $J$-conjugates of $y$ which are exactly the transpositions of $J$. Two neighbors $a, b$ of $x$ are adjacent if and only if $(x, a, b)$ is a triangle in $\Gamma$. By the transitivity on triangles, this is the case if and only if we find $g \in G$ such that $(x, a, b) = (x, y, z)^g$ or, equivalently, if and only if we find $g \in J$ such that $(a, b) = (y, z)^g$. By assumption, $y, z$ are commuting transpositions of $J \cong \mathrm{Sym}_n$. But two transpositions $a, b$ in $J$ are conjugate to the two commuting transpositions $y, z$ if and only if they commute themselves. Since this is the case precisely when $a, b$ have disjoint support we see that $x^\perp$ is isomorphic to the Kneser graph $K(n, 2)$. By vertex-transitivity it follows that $\Gamma$ is locally $K(n, 2)$. Since $n \geqslant 7$ and $\Gamma$ is connected Theorem 2.8 implies that $\Gamma \cong K(n+2, 2)$.

The $J$-conjugates of $y$ generate $J$, and likewise the $K$-conjugates of $x$ generate $K$. The conjugates of $x$ thus generate $G = \langle J, K \rangle$ and Proposition 3.3 implies that $G/Z(G)$ acts faithfully on $\Gamma$. Note that $Z(G) \leqslant C_G(x)$ and hence $Z(G) \leqslant Z(C_G(x)) = \langle x \rangle$ by Proposition 1.9. Since $x \notin Z(G)$ we find that the center of $G$ is trivial. Consequently, $G$ identifies with a subgroup of $\mathrm{Aut}(\Gamma)$. According to Corollary 1.59 $\mathrm{Aut}(\Gamma)$ is isomorphic to $\mathrm{Sym}_{n+2}$. Since $G$ acts transitively on $\Gamma$ Lemma 1.15 finally implies that in fact $G \cong \mathrm{Sym}_{n+2}$. □

An alternative proof of Theorem 3.4 which does not require knowledge of the automorphism groups of the Kneser graphs is outlined in Remark 3.8.

We can prove a slightly more general version of Theorem 3.4 by making use of Theorem 1.27 which showed that the symmetric groups $\mathrm{Sym}_n$ have no outer automorphisms unless $n = 6$.

**Theorem 3.5.** *Let $n \geqslant 7$, and let $G$ be a group with involutions $x, y$ such that*

- $J \trianglelefteq C_G(x)$ *for some* $J \cong \mathrm{Sym}_n$,
- $K \trianglelefteq C_G(y)$ *for some* $K \cong \mathrm{Sym}_n$,
- $x$ *is a transposition in* $K$,

- *y is a transposition in J,*
- *J ∩ K contains an involution z that is a transposition in both J and K.*

If $G = \langle J, K \rangle$, then $G/Z(G) \cong \mathrm{Sym}_{n+2}$.

**Proof.** Basically, we may just copy the proof of Theorem 3.4. Let $\Gamma$ be the graph on $x^G$ with edges $\{x, y\}^G$. The crucial observation is that the neighbors of $x$ are again the $J$-conjugates of $y$. To see this, recall that by Theorem 1.27 all automorphisms of $J$ are inner. Since $J \trianglelefteq C_G(x)$, conjugation by some $g \in C_G(x)$ leaves $J$ invariant and hence induces an automorphism on $J$. We then find $g' \in J$ such that the conjugation action of $g$ and $g'$ agrees on $J$. Hence

$$x^\perp = \{y^g : g \in C_G(x)\} = \{y^{g'} : g' \in J\},$$

and the neighbors of $x$ are indeed the $J$-conjugates of $y$. The subsequent arguments in the proof of Theorem 3.4 only depend on this characterization of adjacency in $\Gamma$. We therefore obtain analogously that $G/Z(G)$ can be embedded in $\mathrm{Aut}(\Gamma) \cong \mathrm{Sym}_{n+2}$.

Note that $x$ is not contained in $J$ because $J$ has trivial center. Thus $C_G(x)$ contains the subgroup $\langle x \rangle \times J$ whose elements belong to distinct cosets of $G/Z(G)$. In particular, the centralizer of $x$ in $G/Z(G)$ has order at least $|\langle x \rangle \times J| = 2n!$ which equals the order of the stabilizer of $x$ in $\mathrm{Aut}(\Gamma)$. Again, Lemma 1.15 shows that $G/Z(G) \cong \mathrm{Sym}_{n+2}$. □

## 3.2 Recognizing Coxeter groups

### 3.2.1 Recognizing $W(A_n)$

Recall that the symmetric group $\mathrm{Sym}_n$ is a Coxeter group of type $A_{n-1}$. Theorem 3.4 can therefore be paraphrased as a recognition result for Coxeter groups of type $A_n$.

**Theorem 3.6.** *Let $n \geqslant 6$, and let $G$ be a group with involutions $x, y$ such that*

- $C_G(x) = \langle x \rangle \times J$ *with* $J \cong W(A_n)$,
- $C_G(y) = \langle y \rangle \times K$ *with* $K \cong W(A_n)$,
- *x is a reflection in K,*
- *y is a reflection in J,*
- *J ∩ K contains an involution z that is a reflection in both J and K.*

*If $G = \langle J, K \rangle$, then $G \cong W(A_{n+2})$.* □

Note that we can state Theorem 3.6 in the following equivalent way.

**Corollary 3.7.** *Let $n \geqslant 6$, and let $G$ be a group with involution $x$ such that*

- $C_G(x) = \langle x \rangle \times J$ *with* $J \cong W(A_n)$,

- $x^u$ is a reflection in $J$ for an involution $u \in G$,
- $J \cap J^u$ contains an involution $z$ that is a reflection in both $J$ and $J^u$.

If $G = \langle J, J^u \rangle$, then $G \cong W(A_{n+2})$.

**Proof.** Just set $y = x^u$ and $K = J^u$. Then $C_G(y) = \langle y \rangle \times K$. Exploiting that $u^2 = 1$, finally note that $x = y^u$ is a reflection in $K = J^u$. □

**Remark 3.8.** Recall that for the proof of Theorem 3.4 (and hence for Theorem 3.6) we constructed the graph $\Gamma$ on the conjugates of the involution $x$ with edges $\{x, y\}^G$. We proceeded by showing that $\Gamma$ is locally $K(n, 2)$ and concluded that $\Gamma$ is isomorphic to the Kneser graph $K(n+2, 2)$. Finally, we used that the automorphism group of $K(n+2, 2)$ is the symmetric group $\mathrm{Sym}_{n+2}$. Determining the center of $G$ and counting the order of $G$ we were then able to deduce that $G$ is isomorphic to $\mathrm{Sym}_{n+2}$. We now present an alternative approach to the last part of this proof which does not require knowledge about the automorphism group of $\Gamma$ but instead makes use of the fact that the symmetric group is a Coxeter group.

Let the assumptions of Theorem 3.6 hold and suppose we already showed that the graph $\Gamma$ with vertices $x^G$ and edges $\{x, y\}^G$ is locally the Weyl graph $\mathbb{W}(A_n)$ and hence that $\Gamma$ is isomorphic to $\mathbb{W}(A_{n+2})$. If $s_1, s_2$ are adjacent in $\Gamma$ then they commute and hence $(s_1 s_2)^2 = 1$. On the other hand, let $s_1$ and $s_2$ be nonadjacent in $\Gamma$. Since $n+3 \geqslant 5$ we find a vertex $s$ in $\Gamma \cong \mathbb{W}(A_{n+2})$ such that $s_1, s_2 \in s^\perp$. When establishing the local structure of $\Gamma$ we showed that this is the case if and only if $s_1$ and $s_2$ are both reflections in $J^g \cong W(A_n)$ where $g \in G$ is such that $s = x^g$. Since the product of two noncommuting reflections in $W(A_n)$ has order 3 we find that $(s_1 s_2)^3 = 1$. Note that by definition of $\mathbb{W}(A_{n+2})$ the graph $\Gamma$ contains an induced subgraph $\Lambda$ isomorphic to the complement of the Coxeter graph of type $A_{n+2}$. We just showed that the order of a product $s_1 s_2$ for $s_1, s_2 \in \Lambda$ is 2 or 3 only depending on whether $s_1$ and $s_2$ are adjacent or not. Consequently, the involutions corresponding to the vertices of $\Lambda$ satisfy the defining relations of a Coxeter group of type $A_{n+2}$. Thus $G$ contains a quotient of $W(A_{n+2})$. Because the only nontrivial proper quotient of $W(A_{n+2})$ is isomorphic to $\mathbb{Z}/2$ we see that $G$ actually contains a copy of $W(A_{n+2})$. Since $G$ acts transitively on $\Gamma$ we have indeed $G \cong W(A_{n+2})$ by Lemma 1.15 .

### 3.2.2 Recognizing $W(F_4)$

Recall that the local recognition result of $W(A_n)$ stated in Theorem 3.6 was based on the local recognition of the corresponding Weyl graph $\mathbb{W}(A_n)$. In the sequel, we wish to apply our local recognition results for $\mathbb{W}(F_4)$ in a similar way to recognize the group $W(F_4)$. There are at least two obstacles here. First of all, $\mathbb{W}(F_4)$ is a bichromatic graph having short and long vertices which need to be distinguished. Secondly, the graph $\mathbb{W}(F_4)$ is not locally recognizable without further conditions. We hope to be able to deal with this second obstacle by imposing conditions for a group theoretical application analogous to the conditions we need to apply the graph theoretical recognition result Theorem 2.45.

**Remark 3.9.** To better understand how to deal with the issue that $\mathbb{W}(F_4)$ has vertices of two different types let us reflect for a moment on the assumptions used for proving Theorem 3.6. The vertices of the graph $\Gamma$ defined in the proof of Theorem 3.6 are the conjugates of an involution $x$, and the edges of $\Gamma$ are constructed from a second element $y$ commuting with $x$. The structure of the centralizer of $x$ allows us to describe the neighbors of $x$ as the $J$-conjugates of $y$. In order to analyze adjacency among the neighbors of $x$ we need at least one edge in $x^\perp$. This edge is provided by the element $z$ which by assumption is adjacent to both $x$ and $y$. If we dropped the condition ensuring existence of a triangle $\{x, y, z\}$ then the graph $\Gamma$ might locally be a coclique.

This discussion motivates the following approach to handling bichromatic graphs in a similar fashion. Since there are short as well as long vertices we need two involutions $x$, $y$ such that the conjugates of $x$ correspond to the short and the conjugates of $y$ correspond to the long vertices. To define edges we need a further short as well as a further long vertex. We also need triangles of every type, that is a triangle of only short vertices another one of two short and one long vertex and so on. To this end, we are prepared to need at least six vertices.

If $x$ is a short root reflection of $W(F_4)$ then its centralizer is $\langle x \rangle \times W(B_3)$. Note that the fact that the centralizer is a Coxeter group itself is not surprising in view of Remark 1.83. Accordingly, we take a closer look at $W(B_3)$ in the following example. Of course, analogous statements are true for $W(C_3)$ which occurs in the centralizer of a long root reflection.

**Example 3.10.** Recall the description of the root system $\Phi(B_3)$ given in Example 1.85 as well as the description of the Weyl graph $\mathbb{W}(B_3)$ in Example 2.21. $W(B_3)$ contains the three short root reflections $x_1, x_2, x_3$ and 6 long root reflections $y_{i,j}$, $i \neq j \in [3]$. The corresponding Weyl graph is

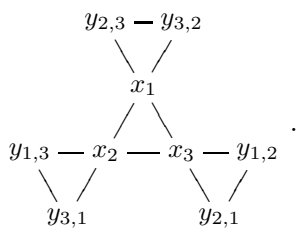

By definition, the order of the product $x_i\, y_{j,k}$ is 4 if $x_i$ and $y_{j,k}$ don't commute or, equivalently, if $i \in \{j, k\}$. We therefore deduce from Proposition 1.18 that

- $y_{j,k}{}^{x_j} = y_{j,k}{}^{x_k} = y_{k,j}$,
- $x_j{}^{y_{j,k}} = x_j{}^{y_{k,j}} = x_k$.

We can rephrase this for instance in the following way. If $x$ and $y$ are a noncommuting short and long root reflection then $y^x$ is the unique long root reflection commuting with $y$ and $x^y$ is the unique short root reflection besides $x$ not commuting with $y$.

Consider two nonadjacent long vertices, say $y_{1,2}$ and $y_{1,3}$. Then $(y_{1,2}\, y_{1,3})^3 = 1$. Using Proposition 1.18 we deduce that $y_{1,2}{}^{y_{1,3}} = y_{1,3}{}^{y_{1,2}}$ does not commute with $y_{1,2}$ or $y_{1,3}$. Hence $y_{1,2}{}^{y_{1,3}}$ equals either $y_{2,3}$ or $y_{3,2}$.

Let $W(B_3)$ act on $\mathbb{W}(B_3)$ by conjugation. The above considerations show that the product $x_1\, x_2\, x_3$ acts trivially. On the other hand, we observe that $W(B_3)$ acts transitively on ordered pairs of nonadjacent vertices of the same type.

**Remark 3.11.** Observe that the automorphism group of $\mathbb{W}(B_3)$ is isomorphic to the wreath product $\mathbb{Z}/2 \wr \mathrm{Sym}_3$ which by Example 1.77 is isomorphic to $W(B_3)$. Notice that this action of $W(B_3)$ on $\mathbb{W}(B_3)$ is necessarily distinct from the action of $W(B_3)$ on $\mathbb{W}(B_3)$ by conjugation which was discussed in Example 3.10. This somewhat subtle point can be observed in Proposition A.8 and Proposition A.9 which show that the graphs $\Gamma_{32a}$ and $\Gamma_{32b}$ are maximally transitive on neighbors while the graphs $\Gamma_{24a}$ and $\Gamma_{24b}$ are not – even though for all four graphs the automorphisms stabilizing a vertex are isomorphic to $W(B_3)$.

Before attempting to recognize a group $G$ as the Coxeter group $W(F_4)$ we present conditions under which we can construct a graph from $G$ that is locally like $\mathbb{W}(F_4)$.

**Proposition 3.12.** *Let $G$ be a group with nonconjugate involutions $x, y$ such that*

- $C_G(x) = \langle x \rangle \times J$ with $J \cong W(B_3)$,
- $C_G(y) = \langle y \rangle \times K$ with $K \cong W(C_3)$,
- *$x$ (respectively $y$) is a short (respectively long) root reflection in $K$ (respectively $J$),*
- *$J \cap K$ contains involutions $x^u$, $y^v$ where $u, v \in G$ such that $x^u$ (respectively $y^v$) is a short (respectively long) root reflection in $K$ as well as in $J$, and*
- *$J \cap J^u$ (respectively $K \cap K^v$) contains an involution that is a short (respectively long) root reflection in both $J$ and $J^u$ (respectively $K$ and $K^v$).*

*Let $\Gamma$ be the bichromatic graph with short vertices the conjugates of $x$ and with long vertices the conjugates of $y$, and edges*

$$E(\Gamma) = \{x, x^u\}^G \cup \{y, y^v\}^G \cup \{x, y\}^G.$$

*Then $\Gamma$ is locally like $\mathbb{W}(F_4)$. Furthermore, if $G = \langle J, K \rangle$ then $\Gamma$ is connected.*

**Proof.** Set $x_1 = x^u$ and $y_1 = y^v$. Further, let $x_2$ (respectively $y_2$) be the involution that is a short (respectively long) root reflection in both $J$ and $J^u$ (respectively $K$ and $K^v$). We have $x_1, x_2, y, y_1 \in J \cong W(B_3)$ where $x_1$ is a short root reflection, and $y$, $y_1$ two commuting long root reflections. Exploiting that $x_1$ commutes with both $y$ and $y_1$, it follows from Example 3.10 that $y_1 = y^{x_2}$. Likewise, $x_1 = x^{y_2}$.

Let $\Gamma$ be the bichromatic graph with short vertices the conjugates of $x$ and with long vertices the conjugates of $y$, and edges

$$E(\Gamma) = \{x, x_1\}^G \cup \{y, y_1\}^G \cup \{x, y\}^G.$$

Note that two short (respectively two long, respectively one short and one long) vertices $a, b$ are adjacent in $\Gamma$ if and only if $(a, b) = (x, x_1)^g$ (respectively $(a, b) = (y, y_1)^g$, respectively $(a, b) = (x, y)^g$) for some $g \in G$. It follows from Proposition 3.3 that $\Gamma$ is connected if $G = \langle J, K \rangle$. The long neighbors of $x$ are the $J$-conjugates of $y$ since

$$x^\perp = \{y^g \colon g \in C_G(x)\} = \{y^g \colon g \in J\}.$$

Likewise the short neighbors of $x$ are the $J$-conjugates of $x_1$. In other words, the long (respectively short) neighbors of $x$ are the long (respectively short) root reflections of $J \cong W(B_3)$. By assumption, $x \perp x_2$, $x_1 \perp x_2$ and likewise $y \perp y_2$, $y_1 \perp y_2$ in $\Gamma$. Since $(x, y)^{x_2 y_2} = (x, y_1)^{y_2} = (x_1, y_1)$, the vertices $x_1$ and $y_1$ are adjacent as well.

In particular, $(x, x_1, x_2)$, $(x, x_1, y)$, $(x, y, y_1)$, $(y, y_1, y_2)$ are ordered triangles of $\Gamma$. We claim that $G$ acts transitively on oriented triangles of the same type. Indeed, let $(a, b, c)$ be an oriented triangle of $\Gamma$ where $a, b, c$ are short vertices. Let $g \in G$ such that $(a, b) = (x, x_1)^g$ and set $d = x_2{}^g$. The vertices $b, c, d$ are short neighbors of $a$, and hence short root reflections in $J^g$. By Example 3.10, we find $h \in J^g$ such that $(b, c) = (b, d)^h$. Therefore $(a, b, c) = (x, x_1, x_2)^{gh}$. Now, let $(a, b, c)$ be an oriented triangle of $\Gamma$ where $a$ is a short and $b, c$ are long vertices. Let $g \in G$ such that $(a, b) = (x, y)^g$. Set $d = y_1{}^g$. $b, c, d$ are long neighbors of $a$, and therefore long root reflections of $J^g$. By assumption, $b$ commutes with $c$ as well with $d$. Consequently, $c = d$, see Example 3.10. Thus $(a, b, c) = (x, y, y_1)^g$. Analogously for triangles of the other types.

We observed that the long (respectively short) neighbors of $x$ are the long (respectively short) root reflections of $J$. Two long neighbors $a, b$ of $x$ are adjacent if and only if $(x, a, b)$ is a triangle in $\Gamma$. By transitivity on oriented triangles of the same type, this is the case if and only if we find $g \in G$ such that $(x, a, b) = (x, y, y_1)^g$, or, equivalently, if and only if we find $g \in J$ such that $(a, b) = (y, y_1)^g$. Two long root reflections $a, b \in J$ are conjugate to the two commuting long root reflections $y, y_1$ if and only if they commute themselves. Thus the long induced subgraph of $x^\perp$ is isomorphic to the long induced subgraph of $\mathbb{W}(B_3)$. Likewise, let $a$ be a short and $b$ a long neighbor of $x$. Again, $a$ and $b$ are adjacent if and only if $(x, a, b)$ is a triangle which is equivalent to finding $g \in J$ such that $(a, b) = (x_1, y)^g$. By Example 3.10, this is the case if and only if $a$ and $b$ commute. Finally, let $a, b$ be two short neighbors of $x$. By Example 3.10, we always find $g \in J$ such that $(a, b) = (x_1, x_2)^g$. Hence, $a$ and $b$ are adjacent. Since $\Gamma$ was constructed to be transitive on short vertices this completes the proof that $\Gamma$ has short local graph $\mathbb{W}(B_3)$. Likewise, we find that the long local graph of $\Gamma$ is $\mathbb{W}(C_3)$. □

We need to add additional assumptions to those given in Proposition 3.12 to deduce that the graph $\Gamma$ defined in Proposition 3.12 is not only locally like $\mathbb{W}(F_4)$ but actually is isomorphic to either $\mathbb{W}(F_4)$ or $\Gamma_{24b}$. In both cases we are then able to deduce that $G$ is isomorphic to the Weyl group $W(F_4)$.

**Theorem 3.13.** *Let $G$ be a group with nonconjugate involutions $x, y$ such that*

- $C_G(x) = \langle x \rangle \times J$ with $J \cong W(B_3)$,
- $C_G(y) = \langle y \rangle \times K$ with $K \cong W(C_3)$,
- *$x$ (respectively $y$) is a short (respectively long) root reflection in $K$ (respectively $J$),*

- $J \cap K$ contains involutions $x_1$, $y_1$ such that $x_1$ (respectively $y_1$) is a short (respectively long) root reflection in $K$ as well as in $J$, and

- there are a long root reflection $y_0 \ne y$, $y_1$ in $J$ and a short root reflection $x_0 \ne x, x_1$ in $K$ such that $x_0$ and $y_0$ commute.

If $G = \langle J, K \rangle$ then $G \cong W(F_4)$.

**Proof.** Observe that $y_0$ does not commute with $x_1$ and that $x_0$ does not commute with $y_1$. Set $x_2 = x_1^{y_0}$ as well as $y_2 = y_1^{x_0}$. Consequently, $x_2 \ne x_1$ is a short root reflection in $J$, and $y_2 \ne y_1$ is a long root reflection in $K$. Summarizing, the elements $x_1$, $x_2$ are short root reflections in $J \cong \mathbb{W}(B_3)$ and $y$, $y_1$, $y_0$ are long root reflections in $J$. Likewise, $y_1$, $y_2$ are long root reflections in $K \cong \mathbb{W}(C_3)$ and $x$, $x_1$, $x_0$ are short root reflections in $J$. Accordingly, the Weyl graphs of $J$ and $K$ are as follows.

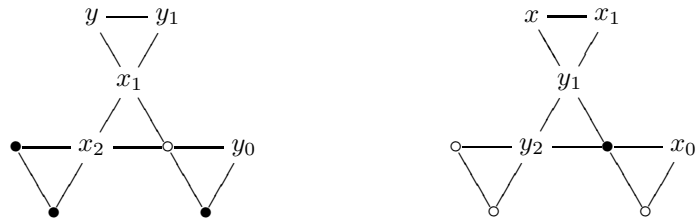

Using Example 3.10, it follows that $x_1 = x^{y_2}$ and $y_1 = y^{y_2}$. Moreover, $(x_0 \, x_1)^3 = 1$ or, equivalently, $x_0^{x_1} = x_1^{x_0}$, see Proposition 1.18, and likewise $y_0^{y_1} = y_1^{y_0}$. Using that $x$ and $y_0$ commute, and that

$$(x, y_0)^{x_0 \, y_1 \, x_0 \, y_0} = (x^{x_2 \, y_0}, y_1^{y_0 \, x_0 \, y_0}) = (x_1^{y_0}, y_1^{x_0}) = (x_2, y_2),$$

we infer that $x_2$ and $y_2$ commute as well. Consequently, $y_2 = y_2^{x_2}$ is a long root reflection in both $K$ and $K^{x_2}$. Analogously, $x_2$ is a short root reflection in both $J$ and $J^{y_2}$. The assumptions of Proposition 3.12 are therefore fulfilled. Hence the graph $\Gamma$ as defined in Proposition 3.12 is connected and locally like $\mathbb{W}(F_4)$.

Recall that $x = x_1^{x_0 \, y_1 \, x_0}$, $y = y_1^{y_0 \, x_1 \, y_0}$, and $x_0^{x_1} = x_1^{x_0}$, $y_0^{y_1} = y_1^{y_0}$. Therefore,

$$(x_0, y_0)^{x_1 \, y_1 \, x_0 \, y_0} = (x_1^{x_0 \, y_1 \, x_0 \, y_0}, y_1^{y_0 \, x_1 \, y_0 \, x_0}) = (x^{y_0}, y^{x_0}) = (x, y)$$

which shows that $x_0$ and $y_0$ are adjacent in $\Gamma$. Notice, for instance using Example 3.10, that $y_0^{x_1 \, y_1}$ is a long root reflection of $J$ not commuting with $y$ and $y_0$. By assumption, $x$ and $y_0$ are adjacent. Since

$$(x_0, y_0^{x_1 \, y_1})^{x_1 \, y_1 \, x_0 \, y_0} = (x, y_0),$$

the vertices $x_0$ and $y_0^{y_1 \, x_1}$ are adjacent as well. We therefore have the following induced subgraph.

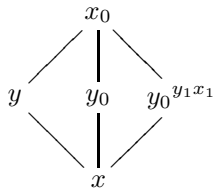

Accordingly, $\mu(x, x_0) = 3$ where $\mu$ is as introduced in Proposition 2.41. Likewise, we find $\mu(y, y_0) = 3$. By construction, $G$ acts transitively on vertices of the same type. Further, recall from Example 3.10 that $J \cong W(B_3)$ acts transitively on ordered pairs of nonadjacent vertices of the same type contained in $x^\perp$. Therefore $G$ acts transitively on oriented 3-paths of the same type.

It follows from Proposition 2.56 that $\mu = 3$. Theorem 2.45 therefore applies, and $\Gamma$ is isomorphic to $\mathbb{W}(F_4)$ or to $\Gamma_{24b}$. In either case, $\mathrm{Aut}(\Gamma)$ is isomorphic to $W(F_4)/Z$.

Note that $J$ is generated by $x_1$, $y_1$, $y_0$, and likewise $K$ is generated by $x_1$, $y_1$, $x_0$. The conjugates of $x$ and $y$ therefore generate $G = \langle J, K \rangle$ and Proposition 3.3 implies that $G/Z(G)$ acts faithfully on $\Gamma$. $Z(G) \leqslant C_G(x)$ and hence $Z(G) \leqslant Z(C_G(x)) = \langle x \rangle \times Z(J)$ by Proposition 1.9. Since $x \notin Z(G)$ we find that $|Z(G)| \leqslant |Z(J)| = 2$. While the stabilizer of a vertex in $\mathrm{Aut}(\Gamma)$ has order 48, the stabilizer in $G$ has order 96. We therefore conclude that $|Z(G)| = 2$. Since $G/Z(G)$ acts transitively on $\Gamma$, Lemma 1.15 implies that $G/Z(G)$ is isomorphic to $\mathrm{Aut}(\Gamma)$. Thus

$$G/Z(G) \cong W(F_4)/Z.$$

With $|G| = 1152$ in mind, the subsequent remark shows that in fact $G \cong W(F_4)$. □

**Remark 3.14.** There is a another approach to proving Theorem 3.13 which does not rely on our previous graph-theoretical recognition results. Suppose that the assumptions of Theorem 3.13 are satisfied. By finding appropriate generators and relations for $G$ we will prove that $G$ is a quotient of $W(F_4)$. Recall from the proof of Theorem 3.13 that $J$ is generated by $x_1$, $y_1$, $y_0$, and likewise $K$ is generated by $x_1$, $y_1$, $x_0$. Accordingly,

$$J = \langle x_1, y_0, y_1 : x_1^2 = y_0^2 = y_1^2 = (x_1 y_0)^4 = (y_0 y_1)^3 = (y_1 x_1)^2 \rangle,$$

and

$$K = \langle y_1, x_0, x_1 : y_1^2 = x_0^2 = x_1^2 = (y_1 x_0)^4 = (x_0 x_1)^3 = (x_1 y_1)^2 \rangle.$$

Since $G = \langle J, K \rangle$ we have $G = \langle x_0, x_1, y_0, y_1 \rangle$, and all the above relations hold. Recall from the proof of Theorem 3.13 that $x = x_1^{x_0 y_1 x_0}$ and $y = y_1^{y_0 x_1 y_0}$. By assumption, $y_0$ commutes with $x$, and $x_0$ commutes with $y$. Further, we assumed that $x_0$ and $y_0$ commute. Summarizing, this accounts for the following additional relations.

$$\begin{aligned}
(x_1^{x_0 y_1 x_0} y_0)^2 &= 1, \\
(x_0 y_1^{y_0 x_1 y_0})^2 &= 1, \\
(x_0 y_0)^2 &= 1.
\end{aligned}$$

The free group generated by four elements $x_0$, $x_1$, $y_0$, $y_1$ together with these relations (actually we can omit one of the first two) and the relations stated above for $J$ and $K$ is isomorphic to $W(F_4)$. For a proof of this statement using GAP see Proposition A.11.

## 3.3 Outlook: Applications to Chevalley groups

### 3.3.1 A characterization of Chevalley groups

We use a characterization based on [Pha70] to introduce Chevalley groups in a way similar to how we introduced Coxeter groups. Recall that we introduced Coxeter groups in Definition 1.60 as groups generated by involutions subject to particularly simple relations. Likewise, we will define (universal) Chevalley groups as groups generated by subgroups isomorphic to $\mathrm{SL}(2,q)$ subject to somewhat analogous relations.

**Example 3.15.** Consider the group $\mathrm{SL}(n+1,q)$ which we identify with $(n+1) \times (n+1)$ matrices over $\mathbb{F}_q$ with determinant 1. For $i \in [n]$ denote with $X_i$ the subgroup corresponding to matrices with block structure

$$\begin{pmatrix} I_{i-1} & & \\ & A & \\ & & I_{n-i} \end{pmatrix}$$

where $I_r$ denotes the $r \times r$ identity matrix, blank entries denote zeros, and $A$ denotes a $2 \times 2$ matrix with determinant 1. Accordingly, $X_i$ is isomorphic to $\mathrm{SL}(2,q)$. Let $T_i$ be the diagonal subgroup of $X_i$. One checks that $\mathrm{SL}(n+1,q)$ is generated by the subgroups $X_i$ and that for all $i \neq j \in [n]$

(a) $\langle X_i, X_j \rangle \cong \begin{cases} \mathrm{SL}(3,q) & \text{if } |i-j|=1, \\ \mathrm{SL}(2,q)^2 & \text{otherwise.} \end{cases}$

(b) $[T_i, T_j] = 1.$

Kok-Wee Phan proves in [Pha70] that $\mathrm{SL}(n+1,q)$ is characterized by these properties for odd $q > 4$, namely that every group generated by subgroups $X_i \cong \mathrm{SL}(2,q)$, $i \in [n]$, satisfying these properties is a quotient of $\mathrm{SL}(n+1,q)$. A simplified proof which relies on the Curtis-Tits theorem was subsequently given by James E. Humphreys in [Hum72]. As Humphreys remarks, his proof generalizes to analogous characterizations for other Chevalley groups.

We define universal Chevalley groups based on characterizations analogous to the one pointed out in Example 3.15. Well-definedness as well as the fact that what we define agrees with the usual notion of a universal Chevalley group is equivalent to the Curtis-Tits theorem in a version described in [Gra08, 4.1.3]. It is based on the results of [Pha70], [Hum72], [Tim04] and [Dun05]. Recall that we classified crystallographic Dynkin diagrams in Theorem 1.82.

**Definition 3.16.** *Let $M$ be a connected crystallographic Dynkin diagram of rank at least three. We say that a group $G$ admits a* Curtis-Tits system *of type $M$ if $G$ is generated by subgroups $X_i \cong \mathrm{SL}(2,q)$, $i \in M$, such that for all vertices $i \neq j$ of $M$*

(a) $\langle X_i, X_j \rangle \cong \begin{cases} \mathrm{SL}(3,q) & \text{if } i \perp j \text{ in } M \text{ with label } 3, \\ \mathrm{Sp}(4,q) & \text{if } i \perp j \text{ in } M \text{ with label } 4, \\ \mathrm{SL}(2,q)^2 & \text{otherwise,} \end{cases}$

(b) $[T_i, T_j] = 1,$

where $T_i$ denotes the diagonal subgroup of $X_i$.

**Definition 3.17.** *Let M be a connected crystallographic Dynkin diagram of rank at least three. The* universal Chevalley group *$M(q)$ is the unique smallest group such that every group admitting a Curtis-Tits system of type M is a quotient of $M(q)$.*

**Example 3.18.** Phan's result stated in Example 3.15 shows that the universal Chevalley group $A_n(q)$ is isomorphic to $\mathrm{SL}(n+1, q)$ for odd $q > 4$.

**Remark 3.19.** Twisted Chevalley groups can be characterized in a similar way. Basically, the generating subgroups isomorphic to $\mathrm{SL}(2, q)$ are replaced by subgroups isomorphic to $\mathrm{SU}(2, q)$, and Curtis-Tits systems are usually replaced by Phan systems. For definitions and details we refer to [Gra08].

### 3.3.2 Classification of the finite simple groups

Consider a group $G$. If we have a normal subgroup $N \trianglelefteq G$ then we can rebuild $G$ from $N$ and $G/N$ together with some information on how to glue these two pieces together. It is therefore natural to look at the building blocks of this decomposition procedure. These are the *simple* groups, namely those groups which contain no normal subgroups besides $\{1\}$ and $G$. The classification of all finite simple groups is considered one of the greatest achievements of all of mathematics. However, "the existing proof of the classification of the finite simple groups runs to somewhere between 10,000 and 15,000 journal pages, spread across some 500 separate articles by more than 100 mathematicians" as it is put in [GLS94]. We state the classification theorem but only comment on the nature of the groups appearing without describing them in detail. For details and definitions we refer to [GLS94].

**Theorem 3.20.** *Every finite simple group is isomorphic to*
- *a cyclic group $\mathbb{Z}/p$ for $p$ prime, or*
- *an alternating group $\mathrm{Alt}_n$ for $n \geqslant 5$, or*
- *a finite simple group of Lie type, or*
- *one of 26 sporadic groups.* □

We are familiar with the cyclic and alternating groups, and we understand the 26 sporadic groups as exceptional phenomena. On the other hand, the bulk of the finite simple groups falls into one of 16 infinite families of the groups of Lie type. These include groups obtained from universal Chevalley groups as defined in Definition 3.17 in most cases by just factoring out the center.

**Example 3.21.** Consider the universal Chevalley group $\mathrm{SL}(n+1, q)$ of type $A_n$. This group is not simple in general because its center, made up by multiples of the identity matrix, is nontrivial whenever $\mathbb{Z}/q$ contains a nontrivial element of order $n+1$. The quotient $\mathrm{PSL}(n+1, q)$, however, is indeed simple if $(n, q) \neq (1, 2), (1, 3)$. A proof and more details can be found for instance in [Rot95, Chapter 8].

**Remark 3.22.** We hint at the nature of the 16 families of groups of Lie type. First, there are the nine families obtained from the universal Chevalley groups of type $A_n$, $B_n$, $C_n$, $D_n$, $E_6$, $E_7$, $E_8$, $F_4$ and $G_2$. Naming differs a lot in literature, and these simple groups are often called Chevalley groups themselves. Additionally, there are families corresponding to the twisted Chevalley groups of the four types $^2A_n$, $^2D_n$, $^3D_4$, $^2E_6$ as well as the three types $^2B_2$, $^2F_4$ and $^2G_2$. The former are often referred to as *Steinberg groups*, and the latter as *Suzuki-Ree groups*. The existence of twisted Chevalley groups can be motivated by recalling the role of diagram automorphisms in proving the isomorphisms $\mathbb{W}(B_n) \cong \mathbb{W}(D_n) \rtimes \mathrm{Sym}_2$ and $\mathbb{W}(F_4) \cong \mathbb{W}(D_4) \rtimes \mathrm{Sym}_3$ stated in Lemma 1.91 and Lemma 1.92. It is such diagram automorphisms that play a role in the definition of Chevalley groups of twisted type. The left superscript corresponds to the order of the used automorphism. We finally remark that the construction of groups of the type $^2B_2$, $^2F_4$ respectively $^2G_2$ only works over fields with order an odd power of 2 respectively 3. A similar motivation can be found in [Coh08].

It is easy to show that the abelian finite simple groups are precisely the cyclic groups of prime order. We are thus left with the classification of the nonabelian finite simple groups. The following remarkable result by Walter Feit and John G. Thompson provides the basis for a strategy.

**Theorem 3.23.** ([FT63]) *All nonabelian finite simple groups have even order.* □

By Cauchy's theorem (which states that a group of order divided by a prime $p$ contains an element of order $p$) this implies that every nonabelian simple group contains an involution. A basic strategy for the classification suggested by Richard Brauer is therefore the following. Assume there was a nonabelian simple group $G$ not among the simple groups listed in Theorem 3.20. According to Theorem 3.23 $G$ contains an involution. Together with further information about the centralizer of this involution and how it interacts with centralizers of other involutions one hopes to identify $G$ with a group that is known to be simple or not simple. In either case, this will be the desired contradiction.

For this reason, recognition results analogous to Theorem 3.6 and Theorem 3.13, which applied to Coxeter groups of type $A_n$ and $F_4$ respectively, are desirable for Chevalley groups.

### 3.3.3 Recognition of Chevalley groups

Recall that the local recognition of Kneser graphs $K(n,2)$ implied the group theoretical recognition result stated in Theorem 3.6 for Coxeter groups of type $A_n$. It may therefore be hoped that Theorem 2.16 which provides a graph theoretical recognition result for a $q$-analog of Kneser graphs implies a recognition result for Chevalley groups of type $A_n$. This is indeed the case as is proved by Ralf Gramlich in his thesis [Gra02]. To state the result we need the notion of a fundamental $\mathrm{SL}(2, q)$ subgroup which will replace the role of a reflection in the statement of Theorem 3.6.

Let $G$ be a subgroup of $\mathrm{GL}(n, \mathbb{F})$, and $V = \mathbb{F}^n$. We define the *commutator*

$$[G, V] \triangleq \{gv - v \colon g \in G, v \in V\}$$

as well as the *centralizer*

$$C_V(G) \triangleq \{v \in V : (\forall g \in G) \ gv = v\}.$$

**Definition 3.24.** *A subgroup $F \leqslant \mathrm{GL}(n, \mathbb{F})$ is said to be a* fundamental $\mathrm{SL}(2, \mathbb{F})$ *subgroup if $F \cong \mathrm{SL}(2, \mathbb{F})$ and*

$$\dim [F, \mathbb{F}^n] = 2, \quad \dim C_{\mathbb{F}^n}(F) = n - 2.$$

**Remark 3.25.** For motivational purposes, recall and compare the definitions of Coxeter groups and universal Chevalley groups. Very roughly, the role of generating involutions of the Coxeter group resembles the role of the generating $\mathrm{SL}(2, q)$ subgroups, and both Coxeter groups and Chevalley groups are described by how two of these generating elements combine. For this reason, one is motivated to think of universal Chevalley groups as $q$-analogs of Coxeter groups. In fact, there exist for instance relations between the orders of these groups analogous to the relation between the orders of $\mathrm{GL}(n, q)$ and $\mathrm{Sym}_n$ pointed out in Remark 1.6. The notion of a fundamental $\mathrm{SL}(2, q)$ subgroup can be interpreted as a $q$-analog of the notion of a reflection.

**Remark 3.26.** There is a one-to-one correspondence between fundamental $\mathrm{SL}(2, q)$ subgroups of $A_n(q) = \mathrm{SL}(n+1, q)$ and nonintersecting line-hyperline pairs of $\mathbb{P}_n(\mathbb{F}_q)$ given by the mapping $F \mapsto ([F, \mathbb{F}_q^n], C_{\mathbb{F}_q^n}(F))$ sending a fundamental $\mathrm{SL}(2, q)$ subgroup to the pair of its commutator and centralizer, see [Gra04]. Moreover, two fundamental $\mathrm{SL}(2, q)$ subgroups $F_1$ and $F_2$ commute if and only if the commutator of each is contained in the centralizer of the other. In the same spirit as the local recognition of the Weyl graphs $\mathbb{W}(A_n)$ is employed to prove Theorem 3.6 in order to recognize Coxeter groups $W(A_n)$, the local recognition of line-hyperline graphs of $\mathbb{P}_n(\mathbb{F}_q)$, see Theorem 2.16, is used to prove the following recognition theorem for the Chevalley groups $A_n(q) = \mathrm{SL}(n+1, q)$.

**Theorem 3.27.** ([Gra04, Theorem 4]) *Let $n \geqslant 6$, $q$ an odd prime power, and let $G$ be a group containing subgroups $X, Y$ isomorphic to $\mathrm{SL}(2, q)$. Denote with $x$ and $y$ the central involutions of $X$ and $Y$. Suppose that*

- *$J \trianglelefteq C_G(x)$ with $J \cong A_n(q)$ and $X \leqslant C_G(J)$,*
- *$K \trianglelefteq C_G(y)$ with $K \cong A_n(q)$ and $Y \leqslant C_G(K)$,*
- *$X$ is a fundamental $\mathrm{SL}(2, q)$ subgroup of $K$,*
- *$Y$ is a fundamental $\mathrm{SL}(2, q)$ subgroup of $J$,*
- *$J \cap K$ contains an involution that is the central involution of a fundamental $\mathrm{SL}(2, q)$ subgroup of both $J$ and $K$.*

*If $G = \langle J, K \rangle$, then $G/Z(G) \cong A_{n+2}(q)/Z$.* □

The bound $n \geqslant 6$ in Theorem 3.27 is optimal as is illustrated by Theorem 3.29. A proof of the special case $n \geqslant 7$ can also be found in [Gra02, 2.5.3]. For a similar result in the case of the twisted Chevalley groups $^2A_n(q)$ see [Gra02, 4.5.6].

**Remark 3.28.** Note that Theorem 3.27 is stated slightly more general then [Gra04, Theorem 4] in that the groups $J, K$ are only required to be normal in the centralizers $C_G(x), C_G(y)$. Recall how we derived Theorem 3.5 from Theorem 3.4 by making use of the fact that the symmetric group $\text{Sym}_n$ has no outer automorphisms when $n \neq 6$. Similarly, the restatement of Theorem 3.27 makes use of the fact that the outer automorphisms of $\text{SL}(n, q)$ are given as products of inner, field, diagonal and diagram automorphisms. The automorphisms of $\text{SL}(n, q)$ were first determined in [SvdW28].

Recently, Ralf Gramlich and Kristina Altmann proved a local recognition result analogous to Theorem 3.27 for the Chevalley groups of type $A_7$ and $E_6$ (as well as for their twisted counterparts) based on results of [Alt07].

**Theorem 3.29. ([Gra08, 7.2.1])** *Let $q$ be an odd prime power, and let $G$ be a group containing subgroups $X, Y$ isomorphic to $\text{SL}(2, q)$. Denote with $x$ and $y$ the central involutions of $X$ and $Y$. Suppose that*

- $J \trianglelefteq C_G(x)$ *with* $J \cong A_5(q)$ *and* $X \leqslant C_G(J)$,
- $K \trianglelefteq C_G(y)$ *with* $K \cong A_5(q)$ *and* $Y \leqslant C_G(K)$,
- $X$ *is a fundamental* $\text{SL}(2, q)$ *subgroup of* $K$,
- $Y$ *is a fundamental* $\text{SL}(2, q)$ *subgroup of* $J$,
- $J \cap K$ *contains an involution that is the central involution of a fundamental $\text{SL}(2, q)$ subgroup of both $J$ and $K$.*

*If $G = \langle J, K \rangle$, then $G/Z(G) \cong A_7(q)/Z$ or $G/Z(G) \cong E_6(q)/Z$.* □

**Remark 3.30.** An outline of the proof can be found in [Gra08]. The basic idea is to look at the commuting graph on the conjugates of $x$ and to find a subgraph which is locally $\mathbb{W}(A_5)$. With the aid of Theorem 2.10 this subgraph can be shown to be isomorphic to either $\mathbb{W}(A_7)$ or $\mathbb{W}(E_6)$. In each case, one is then able to construct a Curtis-Tits system inside $G$. This is done similar to the recognition of $W(A_n)$ outlined in Remark 3.8.

It would be desirable to have a similar result for recognizing the Chevalley group of type $F_4$. However, the strategy employed in proving Theorem 3.29 or its twisted counterpart, see [Gra08], needs to be modified. Hopefully, the analysis of this thesis can provide some starting points as well as hint at possible obstacles involved in finding a recognition result for type $F_4$.

# Appendix A  Computer code

## A.1  Computations in SAGE

### A.1.1  The symplectic graphs

The graphs $\mathcal{S}p_2(2n)$ and $\mathcal{N}\mathcal{S}p^\varepsilon(2n)$ have been introduced in Definition 1.41. They can be implemented in SAGE as follows.

```
def Sp2n(n):
    V = VectorSpace(GF(2), n)
    B = lambda x,y: sum([x[k]*y[k+(-1)**(k%2)] for k in range(n)])
    return Graph([range(1,2**n),
        lambda x,y: x!=y and B(V[x],V[y])==0])
```

```
def NSp2n(n, e):
    V = VectorSpace(GF(2), n)
    B = lambda x,y: sum([x[k]*y[k+(-1)**(k%2)] for k in range(n)])
    Q = lambda x: (sum([x[2*k]*x[2*k+1] for k in range(n/2)])
        + (e==-1 and x[0]**2+x[1]**2))
    return Graph([[ x for x in range(1,2**n) if Q(V[x])==1 ],
        lambda x,y: x!=y and B(V[x],V[y])==0])
```

### A.1.2  The Weyl graphs

We implement the Weyl graphs $\mathbb{W}(A_n)$, $\mathbb{W}(B_n)$, $\mathbb{W}(D_n)$, $\mathbb{W}(F_4)$, $\mathbb{W}(E_6)$, $\mathbb{W}(E_7)$ and $\mathbb{W}(E_8)$ in SAGE based on the descriptions of the corresponding root systems given in Example 1.84, Example 1.85, Example 1.87, Example 1.89 and Example 1.88.

```
def WeylGraphA(n):
    B = (QQ^(n+1)).basis()
    # (n+1)n/2 roots e_i-e_j in R^(n+1)
    R = [ B[i]-B[j] for j in range(n+1) for i in range(j) ]
    return Graph([range(len(R)), lambda x,y: R[x]*R[y]==0])
```

```
def WeylGraphB(n):
    B = (QQ^n).basis()
    # n short roots e_i
    R = [ B[i] for i in range(n) ]
    # n(n-1) long roots e_i+-e_j
    R += [ B[i]+c*B[j] for j in range(n) for i in range(j)
            for c in (-1,1) ]
    return Graph([range(len(R)), lambda x,y: R[x]*R[y]==0])

def WeylGraphD(n):
    B = (QQ^n).basis()
    # n(n-1) roots e_i+-e_j
    R = [ B[i]+c*B[j] for j in range(n) for i in range(j)
            for c in (-1,1) ]
    return Graph([range(len(R)), lambda x,y: R[x]*R[y]==0])

def WeylGraphF4():
    B = (QQ^4).basis()
    # 4+8 short roots e_i, 1/2(e_1+-e_2+-e_3+-e_4)
    R = [ B[i] for i in range(4) ]
    R += [ 1/2*(sum([c[i]*B[i] for i in range(4)]))
            for c in tuples((-1,1), 4) if c[0]==1 ]
    # 12 long roots e_i+-e_j
    R += [ B[i]+c*B[j] for j in range(4) for i in range(j)
            for c in (-1,1) ]
    return Graph([range(len(R)), lambda x,y: R[x]*R[y]==0])

def WeylGraphE8():
    B = (QQ^8).basis()
    # 56 roots e_i+-e_j
    R = [ B[i]+c*B[j] for j in range(8) for i in range(j)
            for c in (-1,1) ]
    # 64 roots of form 1/2 Sum +-e_i
    R += [ 1/2*sum((c[i]*B[i] for i in range(8)))
            for c in tuples((-1,1), 8) if c[0]==1 and prod(c)==1 ]
    return Graph([range(len(R)), lambda x,y: R[x]*R[y]==0])

def WeylGraphE7():
    B = (QQ^8).basis()
    # 30+1 roots e_i+-e_j, i,j<=6, e_7-e_8
    R = [ B[i]+c*B[j] for j in range(6) for i in range(j) for
            c in (-1,1) ]
    R += [ B[6]-B[7] ]
    # 32 roots of form 1/2 Sum +-e_i
    R += [ 1/2*sum((c[i]*B[i] for i in range(8)))
            for c in tuples((-1,1), 8) if c[0]==1
            and prod(c)==1 and c[6]==-c[7] ]
    return Graph([range(len(R)), lambda x,y: R[x]*R[y]==0])
```

```
def WeylGraphE6():
    B = (QQ^8).basis()
    # 20 roots e_i+-e_j, i,j<=5
    R = [ B[i]+c*B[j] for j in range(5) for i in range(j)
            for c in (-1,1) ]
    # 16 roots of form 1/2 Sum +-e_i
    R += [ 1/2*sum((c[i]*B[i] for i in range(8)))
            for c in tuples((-1,1), 8) if c[0]==1
            and prod(c)==1 and c[5]==c[6]==-c[7] ]
    return Graph([range(len(R)), lambda x,y: R[x]*R[y]==0])
```

We can now verify Proposition 2.22 in SAGE.

**Proposition A.1.** *We have the following isomorphisms.*

- $\mathbb{W}(E_6) \cong \mathcal{NS}p^-(6)$,
- $\mathbb{W}(E_7) \cong \mathcal{S}p_2(6)$,
- $\mathbb{W}(E_8) \cong \mathcal{NS}p^+(8)$.

**Proof.** Using the previous implementations of the involved graphs these statements are proved by the following calculations in SAGE.

```
WeylGraphE6().is_isomorphic(NSp2n(6,-1))
>> True
WeylGraphE7().is_isomorphic(Sp2n(6))
>> True
WeylGraphE8().is_isomorphic(NSp2n(8,1))
>> True
```
□

Consider the reflection graphs $\mathbb{W}(H_3)$ and $\mathbb{W}(H_4)$ defined in Remark 2.27.

**Proposition A.2.** *The graph $\mathbb{W}(H_3)$ is isomorphic to $5 \cdot K_3$, and $\mathbb{W}(H_4)$ is a connected graph on $60$ vertices that is locally $\mathbb{W}(H_3)$.*

**Proof.** The first claim is checked as follows.

```
H = gap.FreeGroup(3)
H3 = H / gap.List([H.1^2, H.2^2, H.3^2,
    (H.1*H.2)^5, (H.2*H.3)^3, (H.1*H.3)^2])
WH3 = Graph([H3.ConjugacyClass(H3.1).List(),
    lambda x,y: x*y==y*x])
WH3.order()
>> 15
forall(WH3.connected_components_subgraphs(),
    lambda H: H.is_isomorphic(graphs.CompleteGraph(3)))
>> (True, None)
```

Moreover, the second claim is shown in the following way.

```
H = gap.FreeGroup(4)
H4 = H / gap.List([H.1^2, H.2^2, H.3^2, H.4^2,
    (H.1*H.2)^5, (H.2*H.3)^3, (H.3*H.4)^3,
    (H.1*H.3)^2, (H.1*H.4)^2, (H.2*H.4)^2])
WH4 = Graph([H4.ConjugacyClass(H4.1).List(),
    lambda x,y: x*y==y*x])
WH4.order()
>> 60
WH4.is_connected()
>> True
WH4.subgraph(WH4.neighbors(WH4.vertices()[0])).is_isomorphic(WH3)
>> True
```

□

### A.1.3 Graphs locally like $\mathbb{W}(F_4)$

In the course of the proof of Proposition 2.46 we constructed the two graphs $\Gamma_{24a}$ and $\Gamma_{24b}$ on 24 vertices that are locally like $\mathbb{W}(F_4)$. Copying this construction we implement both graphs in SAGE. Note that the implementations differ only in the last two adjacencies in each of the last four rows.

```
G24a = Graph({
    'x1': ['x2', 'x3', 'y12', 'y21', 'y13', 'y31', 'y14', 'y41'],
    'x2': ['x3', 'x4', 'y12', 'y21', 'y23', 'y32', 'y24', 'y42'],
    'x3': ['x4', 'x1', 'y13', 'y31', 'y23', 'y32', 'y34', 'y43'],
    'x4': ['x1', 'x2', 'y14', 'y41', 'y24', 'y42', 'y34', 'y43'],
    'y12': ['y21', 'y34', 'y43'],
    'y21': ['y12', 'y34', 'y43'],
    'y34': ['y12', 'y21', 'y43'],
    'y13': ['y31', 'y24', 'y42'],
    'y31': ['y13', 'y24', 'y42'],
    'y24': ['y13', 'y31', 'y42'],
    'y14': ['y41', 'y23', 'y32'],
    'y41': ['y14', 'y23', 'y32'],
    'y23': ['y14', 'y41', 'y32'],
    'x5': ['x6', 'x7', 'y12', 'y34', 'y13', 'y24', 'y14', 'y23'],
    'x6': ['x7', 'x8', 'y12', 'y34', 'y42', 'y31', 'y32', 'y41'],
    'x7': ['x8', 'x5', 'y43', 'y21', 'y13', 'y24', 'y32', 'y41'],
    'x8': ['x5', 'x6', 'y43', 'y21', 'y42', 'y31', 'y14', 'y23'],
    'x09': ['x10', 'x11', 'y12', 'y43', 'y13', 'y42', 'y14', 'y32'],
    'x10': ['x11', 'x12', 'y12', 'y43', 'y24', 'y31', 'y23', 'y41'],
    'x11': ['x12', 'x09', 'y34', 'y21', 'y13', 'y42', 'y23', 'y41'],
    'x12': ['x09', 'x10', 'y34', 'y21', 'y24', 'y31', 'y14', 'y32']
})
```

```
G24b = Graph({
    'x1': ['x2', 'x3', 'y12', 'y21', 'y13', 'y31', 'y14', 'y41'],
    'x2': ['x3', 'x4', 'y12', 'y21', 'y23', 'y32', 'y24', 'y42'],
    'x3': ['x4', 'x1', 'y13', 'y31', 'y23', 'y32', 'y34', 'y43'],
    'x4': ['x1', 'x2', 'y14', 'y41', 'y24', 'y42', 'y34', 'y43'],
    'y12': ['y21', 'y34', 'y43'],
    'y21': ['y12', 'y34', 'y43'],
    'y34': ['y12', 'y21', 'y43'],
    'y13': ['y31', 'y24', 'y42'],
    'y31': ['y13', 'y24', 'y42'],
    'y24': ['y13', 'y31', 'y42'],
    'y14': ['y41', 'y23', 'y32'],
    'y41': ['y14', 'y23', 'y32'],
    'y23': ['y14', 'y41', 'y32'],
    'x5': ['x6', 'x7', 'y12', 'y34', 'y13', 'y24', 'y14', 'y23'],
    'x6': ['x7', 'x8', 'y12', 'y34', 'y42', 'y31', 'y32', 'y41'],
    'x7': ['x8', 'x5', 'y43', 'y21', 'y13', 'y24', 'y32', 'y41'],
    'x8': ['x5', 'x6', 'y43', 'y21', 'y42', 'y31', 'y14', 'y23'],
    'x09': ['x10', 'x11', 'y12', 'y43', 'y13', 'y42', 'y41', 'y23'],
    'x10': ['x11', 'x12', 'y12', 'y43', 'y24', 'y31', 'y32', 'y14'],
    'x11': ['x12', 'x09', 'y34', 'y21', 'y13', 'y42', 'y32', 'y14'],
    'x12': ['x09', 'x10', 'y34', 'y21', 'y24', 'y31', 'y41', 'y23']
})
```

**Proposition A.3.** *The graphs $\Gamma_{24a}$ and $\Gamma_{24b}$ are nonisomorphic. $\Gamma_{24a}$ is isomorphic to the Weyl graph $\mathbb{W}(F_4)$.*

**Proof.** We verify these claims in SAGE.

```
G24a.is_isomorphic(G24b)
>> False
G24a.is_isomorphic(WeylGraphF4())
>> True
```

□

**Proposition A.4.** *The bichromatic graphs $\Gamma_{24a}$ and $\Gamma_{24b}$ have isomorphic automorphism groups.*

**Proof.** In order to obtain in SAGE the automorphism groups of the bichromatic graphs $\Gamma_{24a}$ and $\Gamma_{24b}$ we use the partition of vertices into short and long ones based on the fact that short vertices have been labeled $x_i$ and long vertices been labeled $y_{i,j}$.

```
A1 = G24a.automorphism_group(partition=[[v for v in G24a.vertices()
    if v[0]==t] for t in ('x', 'y')])
A2 = G24b.automorphism_group(partition=[[v for v in G24b.vertices()
    if v[0]==t] for t in ('x', 'y')])
A1.group_id() == A2.group_id() == [576, 8654]
>> True
```

□

**Proposition A.5.** *The bichromatic graphs $\Gamma_{24a}$ and $\Gamma_{24b}$ are transitive on vertices of the same type.*

**Proof.** The claim follows by showing that there are exactly two orbits of vertices under the action of the automorphism group.

```
A = G24a.automorphism_group(partition=[[v for v in G24a.vertices()
    if v[0]==t] for t in ('x', 'y')])
len(A.orbits())
>> 2
```

The same code fragment shows that $\Gamma_{24b}$ is transitive on vertices of the same type. □

During the proof of Theorem 2.62 we constructed the graphs $\Gamma_{32a}$ and $\Gamma_{32b}$ which are locally like $\mathbb{W}(F_4)$. In SAGE these constructions read as follows.

```
G32a = Graph({
    'x1': ['x2', 'x3', 'x4', 'y1', 'y2', 'y3', 'y4', 'y5', 'y6'],
    'x2': ['x3', 'y1', 'y2', 'y31', 'y35', 'y51', 'y53'],
    'x3': ['x4', 'y3', 'y4', 'y13', 'y15', 'y51', 'y53'],
    'x4': ['x2', 'y5', 'y6', 'y13', 'y15', 'y31', 'y35'],
    'y1': ['y2', 'y13', 'y15', 'z13', 'z14', 'z15', 'z16'],
    'y2': ['y1', 'y13', 'y15', 'z23', 'z24', 'z25', 'z26'],
    'y3': ['y4', 'y31', 'y35', 'z13', 'z23', 'z35', 'z36'],
    'y4': ['y3', 'y31', 'y35', 'z14', 'z24', 'z45', 'z46'],
    'y5': ['y6', 'y51', 'y53', 'z15', 'z25', 'z35', 'z45'],
    'y6': ['y5', 'y51', 'y53', 'z16', 'z26', 'z36', 'z46'],
    'z13': ['z14', 'z23', 'z24'],
    'z14': ['z23', 'z24'],
    'z23': ['z24'],
    'z15': ['z16', 'z25', 'z26'],
    'z16': ['z25', 'z26'],
    'z25': ['z26'],
    'z35': ['z36', 'z45', 'z46'],
    'z36': ['z45', 'z46'],
    'z45': ['z46'],
    'y13': ['y15', 'z13', 'z14', 'z25', 'z26'],
    'y15': ['y13', 'z15', 'z16', 'z23', 'z24'],
    'y31': ['y35', 'z13', 'z23', 'z45', 'z46'],
    'y35': ['y31', 'z35', 'z36', 'z14', 'z24'],
    'y51': ['y53', 'z15', 'z25', 'z36', 'z46'],
    'y53': ['y51', 'z35', 'z45', 'z16', 'z26'],
    'w1': ['w2', 'w3', 'z13', 'z24', 'z15', 'z26', 'z35', 'z46'],
    'w2': ['w3', 'w4', 'z13', 'z24', 'z16', 'z25', 'z36', 'z45'],
    'w3': ['w4', 'w1', 'z14', 'z23', 'z15', 'z26', 'z36', 'z45'],
    'w4': ['w1', 'w2', 'z14', 'z23', 'z16', 'z25', 'z35', 'z46']
})
```

```
G32b = Graph({
    'x1': ['x2', 'x3', 'x4', 'y1', 'y2', 'y3', 'y4', 'y5', 'y6'],
    'x2': ['x3', 'y1', 'y2', 'y31', 'y35', 'y51', 'y53'],
    'x3': ['x4', 'y3', 'y4', 'y13', 'y15', 'y51', 'y53'],
    'x4': ['x2', 'y5', 'y6', 'y13', 'y15', 'y31', 'y35'],
    'y1': ['y2', 'y13', 'y15', 'z13', 'z14', 'z15', 'z16'],
    'y2': ['y1', 'y13', 'y15', 'z23', 'z24', 'z25', 'z26'],
    'y3': ['y4', 'y31', 'y35', 'z13', 'z23', 'z35', 'z36'],
    'y4': ['y3', 'y31', 'y35', 'z14', 'z24', 'z45', 'z46'],
    'y5': ['y6', 'y51', 'y53', 'z15', 'z25', 'z35', 'z45'],
    'y6': ['y5', 'y51', 'y53', 'z16', 'z26', 'z36', 'z46'],
    'z13': ['z14', 'z23', 'z24'],
    'z14': ['z23', 'z24'],
    'z23': ['z24'],
    'z15': ['z16', 'z25', 'z26'],
    'z16': ['z25', 'z26'],
    'z25': ['z26'],
    'z35': ['z36', 'z45', 'z46'],
    'z36': ['z45', 'z46'],
    'z45': ['z46'],
    'y13': ['y15', 'z13', 'z14', 'z25', 'z26'],
    'y15': ['y13', 'z15', 'z16', 'z23', 'z24'],
    'y31': ['y35', 'z13', 'z23', 'z45', 'z46'],
    'y35': ['y31', 'z35', 'z36', 'z14', 'z24'],
    'y51': ['y53', 'z15', 'z25', 'z36', 'z46'],
    'y53': ['y51', 'z35', 'z45', 'z16', 'z26'],
    'w1': ['w2', 'w3', 'z13', 'z24', 'z15', 'z26', 'z36', 'z45'],
    'w2': ['w3', 'w4', 'z13', 'z24', 'z16', 'z25', 'z35', 'z46'],
    'w3': ['w4', 'w1', 'z14', 'z23', 'z15', 'z26', 'z35', 'z46'],
    'w4': ['w1', 'w2', 'z14', 'z23', 'z16', 'z25', 'z36', 'z45']
})
```

**Proposition A.6.** *The graphs $\Gamma_{32a}$ and $\Gamma_{32b}$ are nonisomorphic.*

**Proof.** We verify this claim in SAGE.

```
G32a.is_isomorphic(G32b)
>> False
```
□

**Proposition A.7.** *The bichromatic graphs $\Gamma_{32a}$ and $\Gamma_{32b}$ are transitive on vertices of the same type.*

**Proof.** The claim follows by showing that there are exactly two orbits of vertices under the action of the automorphism group.

```
A = G32a.automorphism_group(partition=
    [[v for v in G32a.vertices() if v[0] in t]
    for t in (('x', 'z'), ('y', 'w'))])
len(A.orbits())
>> 2
```

The same code fragment shows that $\Gamma_{32b}$ is transitive on vertices of the same type. $\square$

**Proposition A.8.** *The graphs $\Gamma_{32a}$ and $\Gamma_{32b}$ are maximally transitive on neighbors.*

**Proof.** Let $\Gamma$ be a bichromatic graph that is locally like $\mathbb{W}(F_4)$. Let $x \in \Gamma$ be a short a vertex. By assumption, $x^\perp$ is isomorphic to the Weyl graph $\mathbb{W}(B_3)$ whose automorphism group is isomorphic to $W(B_3)$ as observed in Remark 3.11. To show that $\Gamma$ is maximally transitive on neighbors it is therefore sufficient to verify that the stabilizer of $x$ is isomorphic to $W(B_3)$ and that it acts faithfully on $x^\perp$.

```
[A, T] = G32a.automorphism_group(partition=
    [[v for v in G32a.vertices() if v[0] in t]
    for t in (('x', 'z'), ('y', 'w'))], translation=True)
C = gap.Stabilizer(A._gap_(), T['x1'])
C.IdGroup()
>> [48, 48]
gap.Stabilizer(C, [T['x2'],T['x3'],
    T['y1'],T['y3'],T['y5']], gap.OnTuples).Order()
>> 1
```

The above verification works the same for any long vertex. Since $\Gamma_{32a}$ is transitive on vertices of the same type as was shown in Proposition A.7 we find that $\Gamma_{32a}$ is maximally transitive on neighbors. Analogously for $\Gamma_{32b}$. $\square$

The graphs $\Gamma_{24a}$ and $\Gamma_{24b}$ are not maximally transitive on neighbors even though the stabilizer of a vertex is isomorphic to $W(B_3)$. This illustrates the point made in Remark 3.11 that the conjugation action of $W(B_3)$ on $\mathbb{W}(B_3)$ is not faithful.

**Proposition A.9.** *The graphs $\Gamma_{24a}$ and $\Gamma_{24b}$ are not maximally transitive on neighbors.*

**Proof.** We proceed as in the proof of Proposition A.8 by showing that the stabilizer of a vertex $x$ is isomorphic to $W(B_3)$. In the case of the graphs $\Gamma_{24a}$ and $\Gamma_{24b}$, however, the action induced by the stabilizer of $x$ is not faithful.

```
[A, T] = G24a.automorphism_group(partition=
    [[v for v in G24a.vertices() if v[0]==t]
    for t in ('x', 'y')], translation=True)
C = gap.Stabilizer(A._gap_(), T['x1'])
>> [48, 48]
gap.Stabilizer(C, [T['x2'],T['x3'],
    T['y12'],T['y13'],T['y14']], gap.OnTuples).Order()
>> 2
```

Hence $\Gamma_{24a}$ is not maximally transitive on neighbors. Analogously for $\Gamma_{24b}$. □

## A.2 Computations in GAP

### A.2.1 Automorphisms of $\mathbb{W}(F_4)$

**Proposition A.10.** *The automorphism group of the bichromatic graph $\mathbb{W}(F_4)$ is isomorphic to $W(F_4)/Z$.*

**Proof.** We calculate in GAP.

```
LoadPackage("grape");;
F := FreeGroup(4);
F4 := F / [F.1^2, F.2^2, F.3^2, F.4^2, (F.1*F.2)^3, (F.3*F.4)^3,
    (F.2*F.3)^4, (F.1*F.3)^2, (F.1*F.4)^2, (F.2*F.4)^2];
WF4 := Graph(F4, [F4.1, F4.4], OnPoints,
    function(x,y) return x<>y and x*y=y*x; end);
A := Stabilizer(AutomorphismGroup(WF4), Filtered(Vertices(WF4),
    x->IsConjugate(F4, VertexName(WF4, x), F4.1)), OnSets);
IdGroup(A) = IdGroup(F4 / Centre(F4));
>> true
```
□

### A.2.2 Presenting $\mathbb{W}(F_4)$

**Proposition A.11.** *The group presented by*

$$\begin{aligned}
\langle x_0, x_1, y_0, y_1 \ : \ & x_0^2,\ x_1^2,\ y_0^2,\ y_1^2, \\
& (x_1\,y_0)^4,\ (y_0\,y_1)^3,\ (y_1\,x_1)^2,\ (y_1\,x_0)^4,\ (x_0\,x_1)^3, \\
& (x_0\,y_0)^2, \\
& (x_1^{x_0\,y_1\,x_0}\,y_0)^2, \\
& (x_0\,y_1^{y_0\,x_1\,y_0})^2 \rangle
\end{aligned}$$

*is isomorphic to $W(F_4)$.*

**Proof.** This is proved by the following calculation in GAP.

```
F := FreeGroup("x0", "x1", "y0", "y1");
AssignGeneratorVariables(F);
G := F / [ x0^2, x1^2, y0^2, y1^2,
        (x1*y0)^4, (y0*y1)^3, (y1*x1)^2, (y1*x0)^4, (x0*x1)^3,
        (x0*y0)^2,
        (x1^(x0*y1*x0) * y0)^2,
        (x0 * y1^(y0*x1*y0))^2 ];
IdGroup(G);
>> [ 1152, 157478 ]
```
□

# Index

# Thesis Declaration

I hereby declare that the whole of this diploma thesis is my own work and that it does not contain material previously published or written by another person unless otherwise referenced or acknowledged.

This work is submitted to Technische Universität Darmstadt in partial fulfillment of the requirements for being awarded a diploma in Mathematics. I declare that it has not been submitted in whole, or in part, for any other degree.

*Darmstadt, 31 Mar 2008*

Armin Straub